\documentclass[reqno]{amsart}

\usepackage{amsfonts}
\usepackage{amscd}
\usepackage{amsbsy}
\usepackage{amsxtra}
\usepackage{amssymb}
\usepackage{epsfig}
\usepackage{epic,eepic}
\usepackage{graphicx}
\usepackage[all]{xy}
\usepackage{xypic}
\usepackage{psfrag}
\usepackage{amsmath}
\usepackage{amsthm}
\usepackage{setspace}
\usepackage{hyperref}
\usepackage{url}
\usepackage{here}
\usepackage{todonotes}
\usepackage{dsfont}
\newlength{\halfbls}
\usepackage{tikz}
\usetikzlibrary{calc}
\usetikzlibrary{decorations.pathreplacing,decorations.markings}
\usetikzlibrary{arrows}

\DeclareRobustCommand{\SkipTocEntry}[9]{}

   \newcommand{\CC}{\mathbb{C}}
  \newcommand{\HH}{\mathbb{H}}

 \newcommand{\PP}{\mathbb{P}}  \newcommand{\NN}{\mathbb{N}}
\newcommand{\QQ}{\mathbb{Q}} \newcommand{\RR}{\mathbb{R}}  

\newcommand{\ZZ}{\mathbb{Z}}

    \newcommand{\Aut}{{\rm Aut}}

\newcommand{\we}{{\rm wt}}

\newcommand{\Cov}{{\rm Cov}}
\newcommand{\Hur}{{\rm Hur}}

 \newcommand{\Stab}{{\rm Stab}} 

\newcommand{\frakz}{\mathfrak{z}}

  \newcommand{\cPP}{{\mathcal P}}  
 \newcommand{\cFF}{{\mathcal F}} \newcommand{\cEE}{{\mathcal E}}

\newcommand{\cQQ}{{\mathcal Q}}

\newcommand{\cQ}{{\mathcal Q}}

\newcommand{\bfa}{{\bf a}}
\newcommand{\bfb}{{\bf b}}

\newcommand{\bfe}{{\bf e}}

\newcommand{\bfk}{{\bf k}}
\newcommand{\bfl}{{\boldsymbol{\ell}}}
\newcommand{\bfm}{{\bf m}}
\newcommand{\bfn}{{\bf n}}

\newcommand{\bfs}{{\bf s}}
\newcommand{\bft}{{\bf t}}

\newcommand{\bfw}{{\bf w}}

\newcommand{\ual}{{\boldsymbol{\alpha}}}

\newcommand{\uga}{{\boldsymbol{\gamma}}}
\newcommand{\usi}{{\boldsymbol{\sigma}}}

\newcommand{\bfdelta}{{\boldsymbol{\delta}}}

\newcommand{\Hmu}{\Pi}

\newcommand{\QM}{{\rm QM}}

   \newcommand{\N}{{\rm N}} 
\newcommand{\Part}{{\mathcal P}}

\newcommand{\odd}{{\rm odd}}
\newcommand{\reg}{{\rm reg}}
\newcommand{\even}{{\rm even}}

\DeclareMathOperator{\vol}{vol}

\newcommand{\mm}{{\rm{\bf{m}}}}
\newcommand{\zz}{{\rm{\bf{z}}}}

\newcommand{\Ev}{\operatorname{Ev}}
\newcommand{\ev}{\operatorname{ev}}

    \newcommand{\ve}{{\varepsilon}}
\newcommand{\ol}{\overline}   \newcommand{\olp}{\overline{p}}
\newcommand{\olal}{\overline{\alpha}}
\newcommand{\ul}{\underline}

\newtheorem{Defi}{Definition}[section]  
    \newtheorem{Prop}[Defi]{Proposition}
\newtheorem{Lemma}[Defi]{Lemma}    \newtheorem{Cor}[Defi]{Corollary}
\newtheorem{Thm}[Defi]{Theorem}    %\newtheorem*{Anm}{Anmerkung}  


\def\={\;=\;}  \def\+{\,+\,}       \def\h{\tfrac12}  
          
 \def\t#1{\tilde{#1}}

  \newcommand\C{\mathbb C}  \newcommand\Z{\mathbb Z}  \newcommand\Q{\mathbb Q} 
\renewcommand\N{\mathbb N}  \renewcommand\P{\mathbb P}   

\def\a{\alpha} \def\b{\beta}     \def\v{\varepsilon}     \def\l{\lambda}
          \def\2{\pi_2}
\def\G{\Gamma}           
\def\zz{\mathbf z}

\def\SL#1{{\rm SL}(2,#1)}
\def\SL#1{{\rm SL}(2,#1)}       
\def\sm#1#2#3#4{\bigl(\smallmatrix#1&#2\\#3&#4\endsmallmatrix\bigr)}

\def\be{\begin{equation}}   \def\ee{\end{equation}}     \def\bes{\begin{equation*}}    \def\ees{\end{equation*}}
\def\ba{\be\begin{aligned}} \def\ea{\end{aligned}\ee}   \def\bas{\bes\begin{aligned}}  \def\eas{\end{aligned}\ees}

\newcommand{\omoduli}[1][g]{{\Omega \mathcal M}_{#1}}

\newcommand{\B}{B} %pillow
\newcommand{\pari}{\operatorname{par}} % parity of the width of the special cylinders

\newcommand{\PC}{{\rm PC}} 

\newcommand{\EO}{E^0} %edges touching 0
\newcommand{\EOl}{E^0_\ell} %loops touching 0
\newcommand{\EP}{E^+} %oriented edges, not touching 0
\newcommand{\EPl}{E^+_\ell} %oriented loops, not touching 0
\newcommand{\EM}{E^-} % non-oriented edges
\newcommand{\OG}{\overline\Gamma} %reduced graph

\def\bw#1{\bigl\langle#1\bigr\rangle_w}   \def\sbw#1{\langle#1\rangle_w}

\newcounter{savedtocdepth}
\newcommand*{\SaveTocDepth}[1]{%
  \addtocontents{toc}{%
    \protect\setcounter{savedtocdepth}{\protect\value{tocdepth}}%
    \protect\setcounter{tocdepth}{#1}%
  }%
}

\definecolor{cqcqcq}{rgb}{0.75,0.75,0.75}
\definecolor{qqzzqq}{rgb}{0,0.6,0}
\definecolor{ccttqq}{rgb}{0.8,0.2,0}
\definecolor{ccqqcc}{rgb}{0.8,0,0.8}
\definecolor{ffcczz}{rgb}{1,0.8,0.6}
\definecolor{qqttff}{rgb}{0,0.2,1}
\definecolor{uququq}{rgb}{0.25,0.25,0.25}
\definecolor{ffqqqq}{rgb}{1,0,0}
\title[Quasimodularity of pillowcase covers]
{Pillowcase covers: Counting Feynman-like graphs associated with quadratic
differentials}

\author{Elise Goujard}
\address{
  Institut de Math{\'e}matiques de Bordeaux,
  Universit{\'e} de Bordeaux, 
351, cours de la Lib{\'e}ration,
F-33405 Talence \\
}
\thanks{Reaseach of the first author was partially supported by a public grant as part of the FMJH.}
\email{elise.goujard@math.u-bordeaux.fr}

\author{Martin M\"oller}
\thanks{Research  of  the second author is partially supported  
  by the DFG-project MO 1884/1-1 and by the LOEWE-Schwerpunkt
``Uniformisierte Strukturen in Arithmetik und Geometrie''.}
\address{
Institut f\"ur Mathematik, Goethe--Universit\"at Frankfurt,
Robert-Mayer-Str.~6--8,
60325 Frankfurt am Main, Germany\\
}
\email{moeller@math.uni-frankfurt.de}

\begin{document}
\bibliographystyle{halpha}

\begin{abstract}
We prove the quasimodularity of generating functions for counting
pillowcase  covers, with and without Siegel-Veech weight. Similar to
prior work on torus covers, the proof is based on analyzing decompositions
of half-translation surfaces into horizontal cylinders. It provides an
alternative proof of the  quasimodularity results of Eskin-Okounkov and 
a practical method to compute area Siegel-Veech constants.
\par
A main new technical tool is a quasi-polynomiality result for
$2$-orbifold Hurwitz numbers with completed cycles.
\end{abstract}

\maketitle

\tableofcontents

\noindent
\SaveTocDepth{1}

%%%%%%%%%%%%%%%%%%%%%%%%%%%%%%%%%%%%%%%%%%%%%%%%%%%%%%%%%%%
%\section*{Introduction}

%%%%%%%%%%%%%%%%%%%%%%%%%%%%%%%%%%%%%%%%%%%%%%%%%%%%%%%%%%%
\section{Introduction}
%%%%%%%%%%%%%%%%%%%%%%%%%%%%%%%%%%%%%%%%%%%%%%%%%%%%%%%%%%%

Mirror symmetry for elliptic curves can be phrased in its tropical
version by stating that the Hurwitz number counting covers of
elliptic curves can be computed as Feynman integrals and that the
corresponding generating functions are quasimodular forms
(\cite{BBBM}, \cite{GMab}, \cite{Dijk}, \cite{KanZag}, \cite{eo}). A
Feynman integral is physics inspired terminology for an integral
over a product of derivatives of propagators (i.e.\ Weierstrass
$\wp$-functions), the form of the product being encoded by a
(Feynman) graph.
\par
The goal of this paper is to show that the mirror symmetry story
has a complete analog in the scope of pillowcase covers, covers
of the projective line with profile $(2,\ldots,2)$ over $3$ points,
with profile $(\nu,2,\ldots,2)$ over a special point and possibly
a fixed finite number of even order branch points elsewhere (see
Section~\ref{sec:countingPillow}).
These pillowcase covers arise naturally in the volume computation
for strata of quadratic differentials.
\par
Our generalized Feynman graphs have a special vertex~$0$ corresponding
to the branch point with profile $(\nu,2,\ldots,2)$ and, most important,
come with an orientation of the half-edges, so that the edge contribution
to the Feynman integrand is $\wp(z_i \pm z_j)$ according to whether
the half-edges are inconsistently or consistently oriented along the edge.
Those generalized Feynman integrals are quasimodular forms, now
for the subgroup $\Gamma_0(2)$, rather than for $\SL\ZZ$ in the
case of torus covers. The argument is a rather straightforward generalization
of the torus cover case, see Theorem~\ref{thm:coeff0QM} and
Theorem~\ref{thm:QMforS}, and compare to \cite[Section~5 and~6]{GMab}.
The two papers are intentionally parallel whenever possible, to facilitate
comparison. In particular, both papers start with a correspondence
theorem (Proposition~\ref{prop:corr}) that can certainly be rephrased in terms
of covers of tropical curves (\cite[Section~8]{GMab}, \cite{BBBM}).
\par
To arrive from there at our main goals, we need moreover a structure
theorem for the algebra of shifted quasi-polynomials and a polynomiality
theorem for orbifold double Hurwitz numbers, see the end of the introduction.
Altogether, we can first give another proof of the following theorem of
Eskin-Okounkov (\cite{eopillow}).
We let $N^\circ(\Hmu) = \sum N^\circ_d(\Hmu) q^d$ be the generating series
of torus covers with branching profile~$\Hmu$.
\par
\begin{Thm}  \label{intro:count} (= Corollary~\ref{cor:noname})
For any ramification profile~$\Hmu$ the counting 
function $N^\circ(\Hmu)$ for connected pillowcase covers of profile~$\Hmu$ 
is a quasimodular form for the group~$\Gamma_0(2)$ of mixed weight less or
equal to $|\Hmu| + \ell(\Hmu)$.
\end{Thm}
\par
\medskip
The starting point of this paper was to obtain the following version for
a weighted count, motivated by the computation of area-Siegel-Veech
constants (see Section~\ref{sec:GraphSumisQM} for a brief introduction, see 
\cite{emz} and \cite{ekz} for more background.).
\par
\begin{Thm} \label{intro:SV}[= Corollary~\ref{Cor:SVmain}]
For any ramification profile $\Hmu$ and any odd integer $p \geq -1$
the generating series $c^\circ_{p}(\Hmu)$ for counting connected
pillowcase covers with $p$-Siegel-Veech weight is a quasimodular form
for the group~$\Gamma_0(2)$ of mixed weight at most $|\Hmu| + \ell(\Hmu)+p+1$.
\end{Thm}
\par
To explain the use of this result, we compare the knowledge about
strata of the moduli space of abelian differentials $\omoduli[g](\mu)$
and quadratic differentials $\cQQ(\mu)$ with respect to Masur-Veech
volumes and Siegel-Veech constants at the time of writing.
In the abelian case our understanding is nearly complete. Siegel-Veech
constants can be computed recursively by computing ratios of
Masur-Veech volumes of boundary strata (\cite{emz}). These volumes
can be computed efficiently by counting torus covers and closed
formulas derived from this (\cite{eo}, \cite{CMZ}). The volumes have
an interpretation as intersection numbers of tautological classes
(\cite{SauvagetMinimal}, \cite{CMS}) and the formulas are well-understood,
so as to give large genus asymptotics in all detail
(\cite{CMZ}, \cite{Aggarwal}, \cite{CMS}). 
\par
For the moduli space of quadratic differentials much less is known, except
for strata of genus zero surfaces whose volumes are
explicitly computable (\cite{aez}). Siegel-Veech constants are also related
to Masur-Veech volumes by a recursive procedure (\cite{masurzorich},
\cite{GoujardSV}). But these volumes are much harder to evaluate for
higher genus, despite the work of \cite{eopillow}, and some hints being
given in \cite{engel}. The behavior of the large genus asymptotics is
conjectured in \cite{DGZZ} for the sequence of principal strata. Only for
the principal strata an interpretation as intersection number is
known (\cite{DGZZ}).
\par
\medskip
In the current status of knowledge, Theorem~\ref{intro:count}
and Theorem~\ref{intro:SV} provide (besides structural insight) a
somewhat reasonable practical way to compute volumes and Siegel-Veech constants
for strata of quadratic differentials by computing the coefficients
for sufficiently many small~$d$ in order to determine the quasimodular
form uniquely and then using the growth rate of the coefficients. 
This procedure is explained, along with the technical steps of the proof,
in an example in Section~\ref{sec:Example}. Algorithms that compute
volumes and Siegel-Veech constants for quadratic differentials as
efficiently as in the abelian case still have to be found.
\par
\medskip
Finally, we explain the 'local' polynomiality results that are the
intrinsic reason for quasimodularity. In the case of torus covers, double
Hurwitz numbers arise naturally by slicing the torus. These numbers are
polynomials if one uses completed cycles at each slice, see e.g.\ \cite{SSZ}.
Together with a theorem that shows that certain graph sums
with polynomial local contributions are quasimodular and a graph
combination argument to pass to completed cycles~$p_\ell$ we obtained
in \cite{GMab}  quasimodularity for torus covers.
\par
In the case of pillowcase covers, slicing the pillowcase, some
$2$-orbifold Hurwitz numbers arise naturally at special slices. We show
in Theorem~\ref{thm:pbarpoly} that these 2-orbifold Hurwitz numbers are 
quasi-polynomials (rather than just piece-wise quasi-polynomials)
if the $2$-orbifold carries only products of the completed cycles $\olp_k$ in the
algebra of shifted symmetric quasi-polynomials (see Section~\ref{sec:vertex}).
Note that quasi-polynomiality of $2$-orbi\-fold Hurwitz numbers fails
even for the completed cycles~$p_\ell$. The {\em quasi}-polynomiality is the
cause of quasimodularity of the associated generating series
for the subgroup $\Gamma_0(2)$ rather than the full group $\SL\ZZ$.
\par
{\bf Acknowledgements:} The authors are very grateful to Dmitri
Zvonkine for suggesting the form of the one-sided pillowcase
vertex operator. We thank Alex Eskin for sharing a manuscript of
an old project with Andrei Okounkov that also discussed local surfaces
and global graphs. Both authors acknowledge the hospitality
of the {\em Max-Planck Institute for Mathematics (MPIM, Bonn)}, where
much of this work was done.

%%%%%%%%%%%%%%%%%%%%%%%%%%%%%%%%%%%%%%%%%%%%%%%%%%%%%%%%%%%%%%%%%%%%%%
\section{Counting covers of the pillow by global graphs}
\label{sec:countingPillow}
%%%%%%%%%%%%%%%%%%%%%%%%%%%%%%%%%%%%%%%%%%%%%%%%%%%%%%%%%%%%%%%%%%%%%%

The goal of this section is the basic correspondence theorem
Proposition~\ref{prop:corr} and its variants. It allows to count covers of
the pillow by counting graphs with various additional decorations.
As in the abelian case, the correspondence theorem works only if we count
coverings without unramified components. We thus start with standard remarks
on the passage between the various ways of imposing connectivity in the
counting problems.

%%%%%%%%%%%%%%%%%
\subsection{Covers of the pillow and their Hurwitz tuples} 
\label{sec:PillowHurwiz}
%%%%%%%%%%%%%%%%%

We give a short introduction to Hurwitz spaces of covers of the pillow
$B \cong\mathbb{C}P^1 $ and recall some basic notions needed in the sequel.
We provide the pillow with the flat metric that identifies~$B$ with
two squares of side length~$1/2$ glued back to back. We will denote the
corners of the pillow by $P_1,\ldots,P_4$.
\par
A {\em pillowcase cover} is a cover of degree~$2d$ of~$B$ fully branched with~$d$
transpositions over three corners of the pillow, with all odd order branching
stocked together with transpositions over the remaining corner of the pillow,
and with all other even order branch points at arbitrary points different
from the corners.
\par
Let $\Hmu = (\mu^{(1)}, \cdots, \mu^{(n+4)})$ consist of the following
types of partitions. We impose that $\mu^{(1)}=(\nu,2^{d-\vert\nu\vert/2})$ where
$\nu $ is a partition of an even number into odd parts, we require
that $\mu^{(2)}=\mu^{(3)}=\mu^{(4)}=(2^{d})$ and finally that
$\mu^{(i+4)}=(\mu_i,1^{2d-\mu_i})$ with~$\mu_i$ a cycle. We call $\Hmu$ a 
{\em ramification profile} and we define~$g$ by the relation 
$$\ell(\mu)+\ell(\nu)-\vert\mu\vert-\vert\nu\vert/2 \= 2-2g\,, $$
where $\ell(\cdot)$ denotes the length of a partition and $|\cdot|$ the
size of a partition. We write $\Hmu_\emptyset$ for the profile with $n=4$
and $\mu^{(1)} = (2^d)$.
\par
Let $H_{d}(\Hmu)$ (or just $H$ if the parameters are fixed) denote the
$n$-dimensional {\em Hurwitz space} of degree $2d$, genus $g$, 
coverings $p: X \to \PP^1$ of a curve of genus zero with $n+4$ branch points
and ramification profile $\Hmu$, i.e.\ we require that over the $i$-th branch
point~$P_i$ there are $\ell(\mu^{(i)})$ ramification points, of ramification 
orders respectively $\mu^{(i)}_j$.
\par
\begin{figure}[h]
\begin{centering}
\begin{tikzpicture}
\tikzset{
    firstarrow/.style={postaction={decorate},
        decoration={markings,mark=at position .95 with
        {\arrow[scale=1.8,line width=.7pt]{>}}}} }
\tikzset{
    secondarrow/.style={postaction={decorate},
        decoration={markings,mark=at position .5 with
        {\arrow[scale=1.8,line width=.7pt]{>}}}} }

\def\P{2.2pt} % pointwidth 

% Linien
\draw (0,0) node(1){} -- (0,4) node(2){} -- (4.25,4) node(3){} -- (8.5,4) node(4){} -- (8.5,0)
node(5){} -- (4.25,0) node(6){} -- cycle;
\draw (4.25,4) -- (4.25,0);
\fill (1) circle (\P) 
      (2) circle (\P)
      (3) circle (\P)
      (4) circle (\P)
      (5) circle (\P)
      (6) circle (\P);
\fill (0.3,0.3) circle (1.2pt);
\draw (6) -- (1) [secondarrow];
\draw (3) -- (2) [secondarrow];
\draw (3) -- (4) [secondarrow];
\draw (6) -- (5) [secondarrow];

% [Gebilde 1] a1
\draw (0.3,0.3) arc (330:263:.3);
\draw (0.3,0.3) arc (-220:-158:.3) [firstarrow];

% [Gebilde 2] a2
\draw (0.3,0.3) .. controls (0.3,2.9) and (0.3,3.5) .. (0,3.7) [firstarrow];
\draw (0.3,0.3) .. controls (0.4,2.3) and (0.4,3.7) .. (0.25,4);

% [Gebilde 3] a3
\draw (0.3,0.3) .. controls (1,3.5) and (2.2,3.7) .. (3.3,3.8);
\draw (0.3,0.3) .. controls (1,3.5) and (2.2,3.6) .. (3.8,3.7);
\draw (3.3,3.8) .. controls (3.4,3.8) and (3.78,3.9) .. (3.8,4) [firstarrow];
\draw (3.8,3.7) .. controls (4.65,3.65) and (4.7,3.9) .. (4.7,4);

% [Gebilde 4] nicht bezeichnet
\draw (0.3,0.3) .. controls (2,2.34) and (2.19,2.47) .. (2.28,2.85) [firstarrow]; % lange Kurven zu Gebilde
\draw (0.3,0.3) .. controls (2,2.25) and (2.33,2.63) .. (2.52,2.7); 
% Gebilde
\draw plot [smooth, tension=1] coordinates { (2.28,2.85) (2.46,3.27) (2.9,3.35) (2.9,2.9) (2.52,2.7)};  % mittlere Kurve
\fill (2.63,3.02) circle (\P); %Punkt

% [Gebilde 5] gamma 4
\draw (0.3,0.3) .. controls (2.25,1.6) and (2.4,1.7) .. (2.6,2.03) [firstarrow]; % lange Kurven zu Gebilde
\draw (0.3,0.3) .. controls (2.25,1.6) and (2.4,1.73) .. (2.73,1.8);
% Gebilde
\draw plot [smooth, tension=1] coordinates { (2.6,2.03) (2.92,2.46) (3.45,2.4) (3.28,1.91) (2.73,1.8)};  % mittlere Kurve
\fill (3.05,2.15) circle (\P); %Punkt

% [Gebilde 6] nicht bezeichnet
\draw (0.3,0.3) .. controls (3,1) and (3.1,1.2) .. (3.2,1.3) [firstarrow]; % lange Kurven zu Gebilde
\draw (0.3,0.3) .. controls (3,1) and (3.1,1.1) .. (3.3,1.09);
% Gebilde
\draw plot [smooth, tension=1] coordinates { (3.2,1.3) (3.5,1.6) (3.9,1.5) (3.75,1.13) (3.3,1.09)};  % mittlere Kurve
\fill (3.58,1.33) circle (\P); % Punkt

% [Gebilde 7] a4
\draw (0.3,0.3) .. controls (4.25,0.3) and (4.35,0.25) .. (4.5,0) [secondarrow];
\draw (0.3,0.3) .. controls (3.9,0.1) and (4,0.08) .. (4,0);

% Skalierung
\coordinate (S1) at (0,2);   
\coordinate (S2) at (8.5,2);   
\draw[rotate around={45:(S1)}] (-.2,2) -- (.2,2);
\draw[rotate around={45:(S2)}] (8.3,2) -- (8.7,2);
% Kurven in rechter oberer und unterer Ecke 
\draw (8,0) arc (180:90:.5) [secondarrow];
\draw (8.5,3.5) arc (270:180:.5) [secondarrow];

% Beschriftung
\tikzstyle{every node}=[font=\normalsize] 
\node [below] at (1) {$P_1$};
\node [below] at (6) {$P_4$};
\node [below] at (5) {$P_1$};
\node [above] at (4) {$P_2$};
\node [above] at (3) {$P_3$};
\node [above] at (2) {$P_2$};
\node (P) at (0.15,0.4) {$P$};
\node (P5) at (4.2,1.45) {$P_5$};
\node (P6) at (3.76,2.3) {$P_6$};
\node (P7) at (3.2,3.17) {$P_7$};
\node (a1) at (0.6,0.12) {$\alpha_1$};
\node (a2) at (0.1,1) {$\alpha_2$};
\node (a3) at (0.74,2.3) {$\alpha_3$};
\node (a4) at (3.3,0.4) {$\alpha_4$};
\node (g1) at (2.5,1.1) {$\gamma_1$};
\node (g2) at (2,1.7) {$\gamma_2$};
\node (g3) at (1.55,2.1) {$\gamma_3$};
\end{tikzpicture} 
\end{centering}
\caption{Standard presentation of $\pi_1(\PP^1\setminus \{P_1,\ldots,P_{n+4}\})$ 
} \label{cap:PP1gens}
\end{figure}
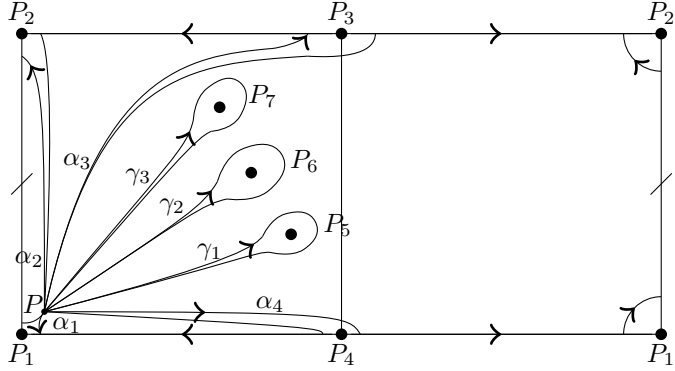

Let $\rho: \pi_1(\PP^1 \setminus \{P_1,\ldots,P_{n+4}\}, \ZZ) 
\to S_{2d}$ be the monodromy representation in the symmetric group of $2d$ elements
associated with a covering in $H$. We use the convention that loops (and elements 
of the symmetric group) are composed from right to left. The elements
$(\alpha_1, \alpha_2, \alpha_3, \alpha_4, \gamma_1, \cdots, \gamma_n)$ as in the 
left picture of Figure~\ref{cap:PP1gens}
generate the fundamental group  $\pi_1(\PP^1 \setminus 
\{P_1,\ldots,P_{n+4}\}, \ZZ)$ with the relation 
\ba\label{eq:FRel}
\alpha_1\alpha_4\gamma_1\dots\gamma_n=\alpha_2^{-1}\alpha_3^{-1}
\ea
%% fig_branchpointpos.tex
%\documentclass{standalone}
%
%\usepackage{tikz}
%\usepackage{pgfplots}
%\usepackage{MnSymbol}
%%\usepackage{amssymb}
%\usepackage{amsmath}
%\usetikzlibrary{calc}
%\usetikzlibrary{decorations.pathreplacing,decorations.markings}
%
%\usetikzlibrary{angles}
%%\usetikzlibrary{intersections}
%\usetikzlibrary{patterns}
%\usetikzlibrary{arrows.meta}
%\usetikzlibrary{decorations.markings, decorations.shapes}
%%\usetikzlibrary{decorations.shapes}
%\usetikzlibrary{arrows}
%\usepackage{adjustbox}
%\usetikzlibrary{positioning}
%
%\usepackage{ifthen}
%
%\newcommand*{\calX}{\mathcal X}
%%\newcommand*{\calY}{{\mathcal Y}
%\newcommand{\bfeta}{{\boldsymbol{\eta}}}
%
%\begin{document}
%\end{document}
Given such a homomorphism $\rho$, we let $\ual_i = \rho(\alpha_i)$, 
and $\uga_i = \rho(\gamma_i)$ and call the tuple
\be \label{eq:HT}
h =(\ual_1, \cdots, \ual_4, \uga_1, \cdots, \uga_n) \in (S_{2d})^{n+4}.
\ee
the {\em Hurwitz tuple} corresponding to $\rho$ and the choice of generators.
Conversely, a Hurwitz tuple as in \eqref{eq:HT} satisfying \eqref{eq:FRel} 
and generating a transitive subgroup of $S_{2d}$
defines a homomorphism $\rho$ and thus a covering $p$. 
We denote by $\Hur^0_d(\Hmu)$ the set of all such  Hurwitz tuples, the upper zero
reflecting that we count connected coverings only. The set of all Hurwitz
tuples (i.e.\ without the requirement of a transitive subgroup) is denoted
by $\Hur_d(\Hmu)$. 
As important technical intermediate notion we need covers
{\em without unramified components}, i.e.\ covers $p: X \to \PP^1$ that
do not have a connected component~$X'$ such that $p|_{X'} = \pi_T \circ p'$
factors into an unramified covering~$p'$ and the torus double covering
$\pi_T: E \to \PP^1$ branched at $P_1,\ldots,P_4$. We let $\ual_1^0 \in S_{2d}$
be the permutation with all transpositions of~$\ual_1$ replaced by the
the identity. In terms of Hurwitz tuples,
we define the tuples without unramified components equivalently as 
\bas
\Hur'_d(\Hmu) &\= \{ h \in \Hur_d(\Hmu)\,:\,
\text{$\langle \ual_1^0, \uga_1, \cdots, \uga_n \rangle$
acts non-trivially on every $\langle h \rangle$-orbit} \} \,.
\eas
The corresponding countings of covers (as usual with weight $1/\Aut(p)$) differ from the
cardinalities of these sets of Hurwitz tuples by the simultaneous conjugation of the 
Hurwitz tuple, hence by a factor of $d!$. Consequently, we let
\be \label{eq:NfromCov}
N_d(\Hmu) \= \frac{|\Hur_d(\Hmu)|}{d!}, \quad
N'_d(\Hmu) \= \frac{|\Hur'_d(\Hmu)|}{d!}, \quad
N^0_d(\Hmu) \= \frac{|\Hur^0_d(\Hmu)|}{d!}\,,
\ee
and package these data into the generating series 
\be \label{eq:genserN}
N(\Hmu) \= \sum_{d=0}^\infty N_d(\Hmu) q^d, \quad N'(\Hmu) \= \sum_{d=0}^\infty 
N_d'(\Hmu) q^d, \quad N^0(\Hmu) \= \sum_{d=0}^\infty N^0_d(\Hmu) q^d \,.
\ee
The connected components of a covering induce a partition of the branch
points of $\ual_1^0$ and $\uga_1,\ldots,\uga_{n}$. This implies that
\begin{equation} \label{eq:NNN}
N'(\Hmu) \= N(\Hmu)/N(\Hmu_\emptyset)\,.
\end{equation}
Similarly, the inclusion-exclusion expression for counting unramified
covers in terms of covers without unramified components carries over
from the case of torus covers (e.g.\ \cite[Proposition~2.1]{GMab}).
\par

%%%%%%%%%%%%%%%%%
\subsection{Covers of the projective line with three marked points} 
\label{sec:CovPP1}
%%%%%%%%%%%%%%%%%

We need coverings of the projective line branched over three points
with two types of parametrizations. As in the case of torus coverings
(see \cite[Section~2.2]{GMab} for more details and remarks on numbered
vs.\ unnumbered enumeration) we define 
\bas
\Cov(\bfw^-, \bfw^+,\mu) &\= \bigl\{(\pi: S \to \PP^1,\sigma_0,\sigma_\infty)\, :  \, \deg(\pi) = 
\sum w_i^+ = \sum w_i^-\,, \\
& \phantom{\= \bigl\{\pi:} \pi^{-1}(1) = [\mu], \,\,\pi^{-1}(0) = \bfw^-,\,\,\
\pi^{-1}(\infty) = \bfw^+ \bigr\}
\eas
to be the set of coverings of~$\PP^1$ with fixed profile~$\mu$
over~$1$ with profile  $\bfw^- = (w_1^-,\ldots,w^-_{n^-})$ and 
$\bfw^+ =  (w_1^+,\ldots,w^+_{n^+})$ over~$0$ and~$\infty$ respectively, and
where $\sigma_0$ and $\sigma_\infty$ are labelings of the branch points over~$0$
and~$\infty$. We usually consider $\bfw^-$ and $\bfw^+$ as 'input' and 'output'
tuples of variables. We denote by
\be \label{eq:defAmu}
 A(\bfw^-, \bfw^+,\mu) \= \sum_{\pi \in \Cov(\bfw^-, \bfw^+,\mu)} \frac{1}{\Aut(\pi)}
\ee
the automorphism-weighted count of these numbers and refer to this quantity
as {\em triple Hurwitz numbers} (although some authors e.g.\ \cite{SSZ}
call them double Hurwitz numbers referring to two sets $\bfw^\pm$ of
free variables).
\par
The second type of covering has only one set of variables and a product
of transpositions of the point at~$\infty$. That is, we define
\bas
\Cov_2(\bfw,\nu) &\= \bigl\{(\pi: S \to \PP^1,\sigma_0)\, :  \, \deg(\pi) = 
\sum w_i\,,  \pi^{-1}(1) = [\nu,2^{(\deg(\pi)-|\nu)/2}]\,, \\
& \phantom{\= \bigl\{\pi:} \,\,
\pi^{-1}(0) = \bfw,\,\,\ \pi^{-1}(\infty) = [2^{\deg(\pi)/2}] \bigr\}
\eas
be the set of coverings with fixed profile over~$1$ and~$\infty$ (but stabilized
by transpositions rather than by adding ones as usual!) and with variable profile
$\bfw = (w_1,\ldots,w_k)$ over~$0$. Finally, we let
\be  \label{eq:defA2mu}
 A_2(\bfw,\nu) \= \sum_{\pi \in \Cov_2(\bfw,\nu)} \frac{1}{\Aut(\pi)}\,.
 \ee
and we refer to them as {\em simple Hurwitz numbers with $2$-stabilization}.
\par
As usual, all these notions have their respective variants for coverings
without unramified components (decorated by a prime) and for connected
coverings (decorated by an upper zero).

%%%%%%%%%%%%%%%%%
\subsection{Global graphs and cylinder decompositions} \label{sec:globalgraph}
%%%%%%%%%%%%%%%%%

We normalize the pillow to be the quotient orbifold $B = E_i/\pm$ where
$E_i = \CC/\ZZ[i]$ is the rectangular torus provided with the unique up
to scale quadratic differential~$q_B$ such that $\pi_T^* q_B$ is holomorphic on~$E$.
The pillow comes with the distinguished points $P_1,\ldots,P_4$ that are
the images of $0, \tfrac{i}2, \tfrac{i+1}2, \tfrac{1}2 \in E_i$ respectively.
For the remaining branch points we usually use in the sequel the {\em branch 
point normalization} that the $i$-th branch point~$P_i$ is the $\pi_T$-image of
a point with coordinates 
$z_i = x_i + \sqrt{-1} \ve_i$ with $0\leq \ve_5 < \ve_6 < \cdots < \ve_{n+4} <1/2$
and any $x_i\in[0,1)$. 
\par
The horizontal foliation with respect to~$q_B$ on~$B$ and thus on every pillowcase
cover $p: X \to B$ with respect to~$p^*q_B$ is periodic. There are two possible
variants to encode the covering by a graph and local data: First, we might
use that the complement of the leaves through all the preimages of the $P_i$
consists of cylinders only. Second, one can use that the cylinders can be continued
across the leaves 'at height~$\tfrac12$', i.e.\ the keeping the leaves through
$P_2$ and~$P_3$ still gives a cylinder decomposition. This can be pushed
further by realizing that the leaves 'at height~$0$' joining two simple
transposition preimages of $P_1$ and $P_4$ can also be added to the cylinders.
In this paper we use the second viewpoint throughout, i.e.\ extending cylinders
as much as possible over fake saddle connections. Said differently, we
mark~$X$ only at the points where~$q$ has a zero or pole, not at the
preimages of the~$P_i$ where $q$ is regular and only remove saddle connections
between those marked points to get a horizontal cylinder decomposition.
\par

\begin{figure}
\mbox{\begin{tikzpicture}[line cap=round,line join=round,>=triangle 45,x=0.8cm,y=0.8cm]
\draw [color=cqcqcq,dash pattern=on 2pt off 2pt, xstep=0.8cm,ystep=0.8cm] (-3.42,-4.68) grid (7,4.64);
\clip(-3.42,-4.68) rectangle (7,4.64);
\draw (2,4)-- (-2,4);
\draw (-0.11,4) -- (0,4.14);
\draw (-0.11,4) -- (0,3.87);
\draw (2,4)-- (6,4);
\draw (4.11,4) -- (4,3.87);
\draw (4.11,4) -- (4,4.14);
\draw (-2,4)-- (-1.5,2.5);
\draw (-1.5,2.5)-- (-3,0);
\draw (-1,0)-- (-3,0);
\draw (-2.21,0) -- (-2.11,0.14);
\draw (-2.21,0) -- (-2.11,-0.14);
\draw (-2,0) -- (-1.9,0.14);
\draw (-2,0) -- (-1.9,-0.14);
\draw (-1,0)-- (1,0);
\draw (0.21,0) -- (0.11,-0.14);
\draw (0.21,0) -- (0.11,0.14);
\draw (0,0) -- (-0.11,-0.14);
\draw (0,0) -- (-0.11,0.14);
\draw (1,0)-- (1.5,-1.5);
\draw (3.5,-1.5)-- (1.5,-1.5);
\draw (2.4,-1.5) -- (2.5,-1.37);
\draw (2.4,-1.5) -- (2.5,-1.64);
\draw (2.61,-1.5) -- (2.71,-1.37);
\draw (2.61,-1.5) -- (2.71,-1.64);
\draw (2.19,-1.5) -- (2.29,-1.37);
\draw (2.19,-1.5) -- (2.29,-1.64);
\draw (4.5,2.5)-- (6.5,2.5);
\draw (5.61,2.5) -- (5.5,2.37);
\draw (5.61,2.5) -- (5.5,2.64);
\draw (5.4,2.5) -- (5.29,2.37);
\draw (5.4,2.5) -- (5.29,2.64);
\draw (5.82,2.5) -- (5.71,2.37);
\draw (5.82,2.5) -- (5.71,2.64);
\draw (3.5,-1.5)-- (3,0);
\draw (3,0)-- (4.5,2.5);
\draw (6.5,2.5)-- (6,4);
\draw (0,-4)-- (0,-3);
\draw (0,-3)-- (1,-3);
\draw (1,-3)-- (1,-4);
\draw (1,-4)-- (0,-4);
\draw [->] (0.5,-1) -- (0.5,-2.5);
\draw [dash pattern=on 5pt off 5pt] (-1.5,2.5)-- (4.5,2.5);
\draw [dash pattern=on 5pt off 5pt] (1,0)-- (3,0);
\draw [shift={(0.5,-2.3)}] plot[domain=4.32:5.11,variable=\t]({1*1.3*cos(\t r)+0*1.3*sin(\t r)},{0*1.3*cos(\t r)+1*1.3*sin(\t r)});
\draw [shift={(0.5,-4.7)},dash pattern=on 2pt off 2pt]  plot[domain=1.18:1.97,variable=\t]({1*1.3*cos(\t r)+0*1.3*sin(\t r)},{0*1.3*cos(\t r)+1*1.3*sin(\t r)});
\draw (1.94,-0.54) node[anchor=north west] {$ C_1 $};
\draw (0.62,1.48) node[anchor=north west] {$ C_2 $};
\draw (1.82,3.46) node[anchor=north west] {$ C_3 $};
\begin{scriptsize}
%losange
%\fill [color=black] (-3,0) ++(-2.0pt,0 pt) -- ++(2.0pt,2.0pt)--++(2.0pt,-2.0pt)--++(-2.0pt,-2.0pt)--++(-2.0pt,2.0pt);
%\fill [color=black,shift={(-1,0)}] (0,0) ++(0 pt,3.0pt) -- ++(2.6pt,-4.5pt)--++(-5.2pt,0 pt) -- ++(2.6pt,4.5pt);
%\fill [color=black] (1,0) ++(-2.0pt,0 pt) -- ++(2.0pt,2.0pt)--++(2.0pt,-2.0pt)--++(-2.0pt,-2.0pt)--++(-2.0pt,2.0pt);
%\fill [color=black] (3,0) ++(-2.0pt,0 pt) -- ++(2.0pt,2.0pt)--++(2.0pt,-2.0pt)--++(-2.0pt,-2.0pt)--++(-2.0pt,2.0pt);
\fill [color=black] (-3,0) ++(-3.0pt,0 pt) -- ++(3.0pt,3.0pt)--++(3.0pt,-3.0pt)--++(-3.0pt,-3.0pt)--++(-3.0pt,3.0pt);
\fill [color=black,shift={(-1,0)}] (0,0) ++(0 pt,3.0pt) -- ++(2.6pt,-4.5pt)--++(-5.2pt,0 pt) -- ++(2.6pt,4.5pt);
\fill [color=black] (1,0) ++(-3.0pt,0 pt) -- ++(3.0pt,3.0pt)--++(3.0pt,-3.0pt)--++(-3.0pt,-3.0pt)--++(-3.0pt,3.0pt);
\fill [color=black] (3,0) ++(-3.0pt,0 pt) -- ++(3.0pt,3.0pt)--++(3.0pt,-3.0pt)--++(-3.0pt,-3.0pt)--++(-3.0pt,3.0pt);
\draw [color=black] (1.5,-1.5) circle (2.0pt);
\draw [color=black] (3.5,-1.5) circle (2.0pt);
\draw [color=black] (-1.5,2.5) circle (2.0pt);
\draw [color=black] (4.5,2.5) circle (2.0pt);
\draw [color=black] (6.5,2.5) circle (2.0pt);
\fill [color=black,shift={(-2,4)},rotate=180] (0,0) ++(0 pt,3.0pt) -- ++(2.6pt,-4.5pt)--++(-5.2pt,0 pt) -- ++(2.6pt,4.5pt);
\fill [color=black,shift={(6,4)},rotate=180] (0,0) ++(0 pt,3.0pt) -- ++(2.6pt,-4.5pt)--++(-5.2pt,0 pt) -- ++(2.6pt,4.5pt);
\fill [color=black,shift={(2,4)},rotate=270] (0,0) ++(0 pt,3.0pt) -- ++(2.6pt,-4.5pt)--++(-5.2pt,0 pt) -- ++(2.6pt,4.5pt);
\fill [color=black] (0,-3) circle (1.0pt);
%nouveau carré
\fill [color=black,shift={(0,-4)},rotate=45] (0,0) ++(-3.0pt,0 pt) -- ++(3.0pt,3.0pt)--++(3.0pt,-3.0pt)--++(-3.0pt,-3.0pt)--++(-3.0pt,3.0pt);
%ancien rond
%\fill [color=black] (0,-4) circle (2.0pt);
\fill [color=black] (1,-4) circle (1.0pt);
\fill [color=black] (1,-3) circle (1.0pt);
\draw [color=black] (0.5,-3.6) circle (2.0pt);
\end{scriptsize}
\end{tikzpicture}}
\mbox{\begin{tikzpicture}[line cap=round,line join=round,>=triangle 45,x=0.8cm,y=0.8cm]
\clip(-3.16,-3.2) rectangle (1.34,6.08);
\draw (-1,3)-- (-1,0);
\draw [shift={(0.58,1.5)}] plot[domain=2.38:3.9,variable=\t]({1*2.18*cos(\t r)+0*2.18*sin(\t r)},{0*2.18*cos(\t r)+1*2.18*sin(\t r)});
\draw [shift={(-2.58,1.5)}] plot[domain=-0.76:0.76,variable=\t]({1*2.18*cos(\t r)+0*2.18*sin(\t r)},{0*2.18*cos(\t r)+1*2.18*sin(\t r)});
\draw (-2.25,1.6) node[anchor=north west] {$ C_1 $};
\draw (-1.10,1.6) node[anchor=north west] {$C_2$};
\draw (-0.32,1.6) node[anchor=north west] {$C_3$};
\draw [shift={(-1.5,5.2)}] plot[domain=4.32:5.11,variable=\t]({1*1.3*cos(\t r)+0*1.3*sin(\t r)},{0*1.3*cos(\t r)+1*1.3*sin(\t r)});
\draw [shift={(-1.5,2.8)},dash pattern=on 2pt off 2pt]  plot[domain=1.18:1.97,variable=\t]({1*1.3*cos(\t r)+0*1.3*sin(\t r)},{0*1.3*cos(\t r)+1*1.3*sin(\t r)});
\draw [shift={(-0.5,5.2)}] plot[domain=4.32:5.11,variable=\t]({1*1.3*cos(\t r)+0*1.3*sin(\t r)},{0*1.3*cos(\t r)+1*1.3*sin(\t r)});
\draw [shift={(-0.5,2.8)},dash pattern=on 2pt off 2pt]  plot[domain=1.18:1.97,variable=\t]({1*1.3*cos(\t r)+0*1.3*sin(\t r)},{0*1.3*cos(\t r)+1*1.3*sin(\t r)});
\draw [shift={(-1.6,4.7)}] plot[domain=4.32:5.11,variable=\t]({1*1.3*cos(\t r)+0*1.3*sin(\t r)},{0*1.3*cos(\t r)+1*1.3*sin(\t r)});
\draw [shift={(-1.6,2.3)},dash pattern=on 2pt off 2pt]  plot[domain=1.18:1.97,variable=\t]({1*1.3*cos(\t r)+0*1.3*sin(\t r)},{0*1.3*cos(\t r)+1*1.3*sin(\t r)});
\draw [shift={(-0.4,4.7)}] plot[domain=4.32:5.11,variable=\t]({1*1.3*cos(\t r)+0*1.3*sin(\t r)},{0*1.3*cos(\t r)+1*1.3*sin(\t r)});
\draw [shift={(-0.4,2.3)},dash pattern=on 2pt off 2pt]  plot[domain=1.18:1.97,variable=\t]({1*1.3*cos(\t r)+0*1.3*sin(\t r)},{0*1.3*cos(\t r)+1*1.3*sin(\t r)});
\draw [shift={(-1,7.76)}] plot[domain=4.41:5.01,variable=\t]({1*3.41*cos(\t r)+0*3.41*sin(\t r)},{0*3.41*cos(\t r)+1*3.41*sin(\t r)});
\draw [shift={(-1,1.24)}] plot[domain=1.27:1.87,variable=\t]({1*3.41*cos(\t r)+0*3.41*sin(\t r)},{0*3.41*cos(\t r)+1*3.41*sin(\t r)});
\draw (-2,4.5)-- (-2,4);
\draw (-2,4)-- (-2.1,3.5);
\draw (-1,4)-- (-0.9,3.5);
\draw (-1,4)-- (-1.1,3.5);
\draw (0,4)-- (0.1,3.5);
\draw (0,4)-- (0,4.5);
\draw (-2.76,4.24) node[anchor=north west] {$S_1$};
\draw (-1.16,5.02) node[anchor=north west] {$w_3$};
\draw (-1.82,3.4) node[anchor=north west] {$w_2$};
\draw (-0.44,3.34) node[anchor=north west] {$w_1$};
\draw (-2.5,-1)-- (-1.5,-1);
\draw (-1.5,-1)-- (-1.5,-1.5);
\draw (-2.5,-1)-- (-2.5,-1.5);
\draw [shift={(-2,-0.3)}] plot[domain=4.32:5.11,variable=\t]({1*1.3*cos(\t r)+0*1.3*sin(\t r)},{0*1.3*cos(\t r)+1*1.3*sin(\t r)});
\draw [shift={(-2,-2.7)},dash pattern=on 2pt off 2pt]  plot[domain=1.18:1.97,variable=\t]({1*1.3*cos(\t r)+0*1.3*sin(\t r)},{0*1.3*cos(\t r)+1*1.3*sin(\t r)});
\draw (-1,-1)-- (-1,-0.5);
\draw (-1,-1)-- (0,-1);
\draw (0,-1)-- (0,-1.5);
\draw (1,-1)-- (1,-1.5);
\draw (1,-1)-- (1,-0.5);
\draw [shift={(0.5,-0.3)}] plot[domain=4.32:5.11,variable=\t]({1*1.3*cos(\t r)+0*1.3*sin(\t r)},{0*1.3*cos(\t r)+1*1.3*sin(\t r)});
\draw [shift={(0.5,-2.7)},dash pattern=on 2pt off 2pt]  plot[domain=1.18:1.97,variable=\t]({1*1.3*cos(\t r)+0*1.3*sin(\t r)},{0*1.3*cos(\t r)+1*1.3*sin(\t r)});
\draw [shift={(0.5,0.2)}] plot[domain=4.32:5.11,variable=\t]({1*1.3*cos(\t r)+0*1.3*sin(\t r)},{0*1.3*cos(\t r)+1*1.3*sin(\t r)});
\draw [shift={(0.5,-2.2)},dash pattern=on 2pt off 2pt]  plot[domain=1.18:1.97,variable=\t]({1*1.3*cos(\t r)+0*1.3*sin(\t r)},{0*1.3*cos(\t r)+1*1.3*sin(\t r)});
\draw [shift={(0,2.76)}] plot[domain=4.41:5.01,variable=\t]({1*3.41*cos(\t r)+0*3.41*sin(\t r)},{0*3.41*cos(\t r)+1*3.41*sin(\t r)});
\draw [shift={(0,-3.76)}] plot[domain=1.27:1.87,variable=\t]({1*3.41*cos(\t r)+0*3.41*sin(\t r)},{0*3.41*cos(\t r)+1*3.41*sin(\t r)});
\draw (-3.3,-0.82) node[anchor=north west] {$S_0$};
\draw (-2.1,-1.58) node[anchor=north west] {$w_3$};
\draw (0.02,0.1) node[anchor=north west] {$w_2$};
\draw (0.42,-1.58) node[anchor=north west] {$w_1$};
\begin{scriptsize}
\fill [color=black] (-1,3) circle (2.0pt);
\draw [color=black] (-1,0)-- ++(-2.0pt,-2.0pt) -- ++(4.0pt,4.0pt) ++(-4.0pt,0) -- ++(4.0pt,-4.0pt);
\draw [color=black] (-1,4) circle (2.0pt);
\fill [color=black,shift={(-2.5,-1)},rotate=180] (0,0) ++(0 pt,3.0pt) -- ++(2.6pt,-4.5pt)--++(-5.2pt,0 pt) -- ++(2.6pt,4.5pt);
\fill [color=black,shift={(-1.5,-1)},rotate=270] (0,0) ++(0 pt,3.0pt) -- ++(2.6pt,-4.5pt)--++(-5.2pt,0 pt) -- ++(2.6pt,4.5pt);
\fill [color=black,shift={(-1,-1)}] (0,0) ++(0 pt,3.0pt) -- ++(2.6pt,-4.5pt)--++(-5.2pt,0 pt) -- ++(2.6pt,4.5pt);
\fill [color=black] (0,-1) ++(-3.0pt,0 pt) -- ++(3.0pt,3.0pt)--++(3.0pt,-3.0pt)--++(-3.0pt,-3.0pt)--++(-3.0pt,3.0pt);
\end{scriptsize}
\end{tikzpicture}}
\caption{A pillowcase cover, the global graph and the local surfaces}
\label{fig:patronQ21iii}
\end{figure}
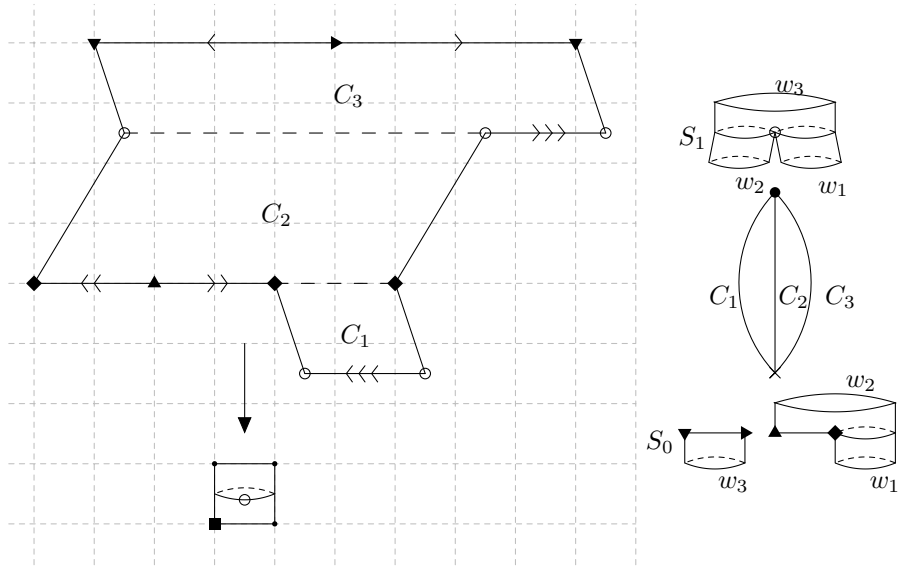

The {\em global graph~$\Gamma$} associated with the pillowcase covering 
surface $(X,q = p^* q_B)$ of ramification profile~$\Hmu$ is the 
graph~$\Gamma$ with $n+1=|\Hmu|-3$ vertices. The vertex with the special
label~$0$ corresponds to the union of leaves through a preimage
of $P_1$ or $P_4$ and the remaining vertices, labeled by $j \in \{1,\ldots,n\}$
correspond to the leaves through $P_{j+4}$.
below. If the partition~$\nu$ is empty, then there is no vertex~$0$. The
edges $E(\Gamma)$
of $\Gamma$ are in bijection with the core curves of the horizontal cylinder
decomposition described in the second viewpoint above.
\par
We illustrate this using Figure~\ref{fig:patronQ21iii} that gives
a covering in the stratum $\cQ(2,1,-1^3)$, in other terms it has a ramification profile given by $n=1$, $\nu=(3,1,1,1)$ and $\mu_5=2$ in the notations of
Section~\ref{sec:PillowHurwiz}. (For more background on strata of quadratic
differentials, including the notation $\cQ(2,1,-1^3)$, see for
example \cite{zorich06}.) The small triangles (with different
orientations) are the three simple poles, the diamond indicates the simple
zero. These points map to the black square on the pillow. The white circle
indicates the double zero, mapping the white circle on the pillow.
This corresponds to the point~$P_5$ while $P_1,\ldots,P_4$ are the corners
of the pillow. 
\par
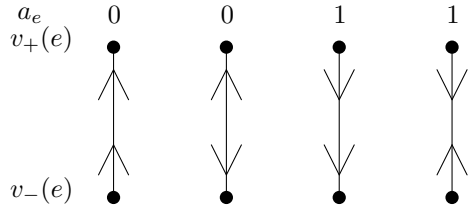
\begin{figure}
\begin{center}
\begin{tikzpicture}[line cap=round,line join=round,>=triangle 45,x=1.0cm,y=1.0cm]
\clip(-4,2.79) rectangle (4,5.72);
\draw (-2,5)-- (-2,3);
\draw (-0.5,5)-- (-0.5,3);
\draw (1,5)-- (1,3);
\draw (2.5,5)-- (2.5,3);
\draw (-2.2,4.3)-- (-2,4.7);
\draw (-2,4.7)-- (-1.8,4.3);
\draw (-2.2,3.3)-- (-2,3.7);
\draw (-2,3.7)-- (-1.8,3.3);
\draw (-0.7,4.3)-- (-0.5,4.7);
\draw (-0.5,4.7)-- (-0.3,4.3);
\draw (0.8,4.7)-- (1,4.3);
\draw (1,4.3)-- (1.2,4.7);
\draw (2.3,4.7)-- (2.5,4.3);
\draw (2.5,4.3)-- (2.7,4.7);
\draw (-0.7,3.7)-- (-0.5,3.3);
\draw (-0.5,3.3)-- (-0.3,3.7);
\draw (0.8,3.7)-- (1,3.3);
\draw (1,3.3)-- (1.2,3.7);
\draw (2.3,3.3)-- (2.5,3.7);
\draw (2.5,3.7)-- (2.7,3.3);
\draw (-3.5,2.8) node[anchor=south west] {$v_-(e)$};
\draw (-3.5,4.8) node[anchor=south west] {$v_+(e)$};
\draw (-3.4,5.2) node[anchor=south west] {$a_e$};
\draw (-2.2,5.2) node[anchor=south west] {$0$};
\draw (-0.7,5.2) node[anchor=south west] {$0$};
\draw (0.8,5.2) node[anchor=south west] {$1$};
\draw (2.3,5.2) node[anchor=south west] {$1$};
\begin{scriptsize}
\fill [color=black] (-2,5) circle (2.5pt);
\fill [color=black] (-2,3) circle (2.5pt);
\fill [color=black] (-0.5,5) circle (2.5pt);
\fill [color=black] (-0.5,3) circle (2.5pt);
\fill [color=black] (1,5) circle (2.5pt);
\fill [color=black] (1,3) circle (2.5pt);
\fill [color=black] (2.5,5) circle (2.5pt);
\fill [color=black] (2.5,3) circle (2.5pt);
\end{scriptsize}
\end{tikzpicture}
\end{center}
\caption{Orientation of half-edges and height minimum~$a_e$} \label{cap:HEorient}
\end{figure}

We provide~$\Gamma$ with an {\em orientation of its half-edges} as follows.
We provide the pillowcase without the special layers with one of the two choices of an orientation
of the vertical direction, say the upward pointing. 
We orient a half-edge at a vertex~$v$ outward-pointing,
if the orientation of the cylinder pointing towards the boundary
representing the half-edge is consistent with the chosen global (``vertical'')
orientation, and we orient the half-edge inward-pointing
otherwise. In particular, all the half-edges starting at the vertex~$0$
(if it exists) are oriented outward-pointing. Recall that cylinders may
cross the special layers at height $0$ or~$1/2$ any number of times.
The two half-edges corresponding to a cylinder are oriented consistently
(see Figure~\ref{cap:HEorient}, leftmost and third arrow) if and only if
the cylinder crosses the special layers an even number of times.  We refer
to this extra datum as an {\em orientation $G \in \Gamma$}.
\par
To reconstruct a pillowcase covering  from a global graph, we need 
as in the case of torus coverings, two types of extra data that encode the geometry of
the cylinders and the geometry of the local surfaces, respectively. The
first extra datum is the cylinder geometry.
\par
Each cylinder (corresponding to an edge~$e$) has an integral positive 
width $w_e$ and a real positive height~$h_e$. The heights~$h_e$ are not 
arbitrary, but related to the position of the branch points. To define the
space parametrizing possible heights, we need a finer classification of the
set $E(\Gamma)$.
\par
We use an upper index~$0$ to denote edges having the vertex~$0$ as an extremity
and we use a lower index~$\ell$ to denote loop edges, i.e.\ edges linking
a vertex to itself. For all non-loop edges we denote by $v^+(e)$ (resp.\  $v^-(e)$) 
the label of the vertex whose height according to the branch point normalization
is higher (resp.\ lower).  We refer to them as the 'upper' (resp.\ 'lower')
vertex of the edge~$e$. For all inconsistently oriented loop edges we extend this notation into $v^+(e)=v^-(e)=v(e)$. Once we provided~$\Gamma$ with an orientation~$G$
we can distinguish those edges with consistent orientation. We denote this subset
by an upper index~$+$. We can now define the {\em height space} to be 
\be \label{eq:heightspace}
\widetilde{\NN}^{E(G)} \= \Bigl\{ (h_e)_{e \in E(\Gamma)} \,:\,
\begin{cases} 
h_e   \in \NN_{>0} \,  &\text{if} \quad e \in\EOl(G)\cup\EPl(G) \\
h_e - \Delta(e)  \in \NN_{\geq a_e} &\text{if} \quad e \in E(G) \setminus (\EOl(G)\cup\EPl(G))
\end{cases} \Bigr\}\,,
\ee
where $\Delta(e) = \pm \ve_{v^+(e)} \pm \ve_{v^-(e)}$, where the sign in front of each $\ve$ is positive
if and only if the edge at the corresponding vertex is incoming and with $a_e$
depending on the orientation as indicated in Figure~\ref{cap:HEorient}
\par
We claim that the collection of heights the cylinders in a pillowcase covering
belongs to the height space and that conversely each element in the height space
can be realized by such a covering. The integrality of the heights corrected
by $\Delta(e)$ follows directly from the branch point normalization and the
conventions of half-edge markings. It remains to justify the lower bounds one
for the corrected heights. This happens if and only if the cylinder has to
go all the way up to the preimage of the height~$1/2$-line and down again.
Loops based at~$0$ and consistently oriented loops have this property. For the
remaining loops, it depends on the orientation of the half-edges. Note that
%$v^\pm(e)$ is not well-defined for loops, but for inconsistently oriented loops
(second and fourth case in Figure~\ref{cap:HEorient}) the lower bound~$a_e$
is independent of the choice. For non-loop edges the integer~$a_e$ encodes
whether the cylinder has to go around the pillowcase. This completes the
proof of the claim.
\par
The last piece of local information for a cylinder is the {\em twist} 
$t_e \in \ZZ \cap [0,w_e-1]$. The twist depends on the choice 
of a ramification point $P^-(e)$ and $P^+(e)$ in each of the two components 
adjacent to the cylinders and it is defined as the integer part
of the real part $\lfloor  \Re(\int_s \omega) \rfloor$ of the
integral along the unique straight line joining $P^-(e)$ to $P^+(e)$
such that $t_e \in [0,w_e-1]$. The exact values of the twist will hardly 
matter in the sequel. It is important to retain simply that there 
are $w_e$ possibilities for the twist in a given cylinder.

%%%%%%%%%%%%%%%%%
\subsection{The basic correspondence theorem} \label{sec:corrthm}
%%%%%%%%%%%%%%%%%

The second extra datum needed for the correspondence theorem is the local
geometry at the vertices. We let $X^0$ be the complement of the core curves
of the cylinders. We call the union of  connected components of $X^0$
that carry the same label the {\em local surfaces} of $(X,\omega)$. We label these local
surfaces also by an integer in $\{0,1,\ldots,n\}$ according to the ramification 
point they carry. This labeling is well-defined, since $p$ is a cover without
unramified components. The restriction of the cover~$p$ to any local surface
besides the one corresponding to the special vertex is metrically the pullback
of an infinite cylinder branched over one point, as in the case of torus
coverings. We thus encode these local surfaces by elements in 
$\Cov'(\bfw_v^-,\bfw_v^+,\mu_v)$ where $\bfw_v^-$ and $\bfw_v^+$ are the
widths of the incoming and outgoing edges. The restriction of the cover~$p$ to
a neighborhood of the line at height~$0$ is precisely the type of
cover parameterized by an element in $\Cov_2(\bfw,\nu)$, with $\bfw$ the
tuple of widths of the outgoing edges. 
\par
For the following proposition we fix a ramification profile~$\Hmu$ and let $\nu$
resp.\ $\mu_v$ be the component of the tuple~$\Hmu$ that corresponds to the vertex~$v$
under the vertex marking conventions explained in Section~\ref{sec:globalgraph}.
\par
\begin{Prop} \label{prop:corr}
There is a bijective correspondence between
\begin{itemize}
\item[i)] flat surfaces $(X,q)$ with a covering $p:X\to \B$ of degree~$2d$ and
profile~$\Hmu$ of the pillow~$\B$ without unramified components and
with $q=\pi^*q_\B$, and
\item[ii)] isomorphism classes of tuples $(G, (w_e,h_e, t_e)_{e \in E(G)}, (\pi_v)_{v\in V(G)})$ consisting of 
\begin{itemize}
\item[$\bullet$] a global graph~$\Gamma$ with labeled vertices including
a special vertex~$0$ if $\nu$ is non-empty, without isolated vertices, together with an
orientation $G \in \Gamma$ of the half edges such that
%: \begin{itemize} \item
all half edges emerging from the special vertex are outgoing, 
%\item any other vertex has at least one incoming and one outgoing half edge
%\end{itemize}
\item[$\bullet$] a collection of real numbers $(w_e,h_e,t_e)_{e \in E(G)}$ 
representing the width, height and twist of the cylinder corresponding to~$e$. 
The widths $w_e$ are integers, the tuple of heights $(h_e)_{e \in E(G)} \in 
\widetilde{\NN}^{E(G)}$ is in the height space, $t_e \in \ZZ \cap [0,w_e-1]$ 
and these numbers satisfy 
\be \label{eq:whd}
2\sum_{e \in E(G)} w_eh_e \= 2d\,,
\ee
\item[$\bullet$] a collection of $\PP^1$-coverings $(\pi_v)_{v\in V(G) \setminus \{0\}}
\in \Cov'(\bfw_v^-,\bfw_v^+,\mu_v)$ without unramified components %, with $\pi_v 
 where $\bfw_v^-$ is the tuple of widths at the incoming edges at $v$, $\bfw_v^+$ is 
the tuple of widths at the outgoing edges at~$v$, and $\mu_v$ is the ramification
profile given by the labels at the vertex~$v$.
\item[$\bullet$] and a $\PP^1$-covering $\pi_0 \in \Cov'_2(\bfw_0,\nu)$ 
where $\bfw_0$ is the tuple of widths at the outgoing edges at $v=0$
and $\nu$ is the ramification profile given by the labels at the vertex~$v=0$.
\end{itemize}
up to the action of the group $\Aut(\Gamma)$ of automorphisms of the
labeled graph~$\Gamma$. 
\end{itemize}
\end{Prop}
\par
\begin{proof}
With the setup and the orientation of half-edges adapted to pillowcase
coverings, the proof proceeds now exactly as in the case of torus covers,
see \cite[Proposition~2.4]{GMab}.
\end{proof}

%%%%%%%%%%%%%%%%%
\subsection{Variants of the correspondence theorem} \label{sec:corrvar}
%%%%%%%%%%%%%%%%%

For the proof of the main theorem we will also need variants of
the correspondence theorem that arise from counting covers by graphs
while declaring a subset of points $P_i$ for $i \in S \subset \{5,\ldots n+4\}$
to be part of the layer of the special vertex. For extreme cases $S = \emptyset$
we are back in the situation of the previous situation while for
$S=\{5,\ldots n+4\}$ the global graph is tautologically just a single vertex
with no edges, decorated by a local Hurwitz number which is just the
global Hurwitz number we are interested in.
\par
For the concrete statement, we start with the {\em branch point normalization.}
We place the points $P_i$ at $z_i = x_i + \sqrt{-1}\ve_i$ where now
$0 < \ve_i < \kappa$ for all $i \in S$ and $\kappa < \ve_i < 1/2$
for all $i \in S^c = \{5,\ldots n+4\} \setminus S$, and moreover within these
constraints strictly increasing with~$i$.
\par
The {\em global graph} associated with $(X,q = p^*q_B)$ is now the following
graph~$\Gamma_S$ with $n+1-|S|$ vertices. The special vertex~$0$ corresponds
to the region $R = \{0 \leq \Im(z) \leq \kappa\}$ and the remaining vertices
are indexed by~$S^c$. Edges correspond the cylinders that are not entirely
contained in a connected component of $p^{-1}(R)$. The notion of an
{\em orientation} $G_S \in \Gamma_S$ carries over verbatim from the above
discussion, and the same holds for the height space $\widetilde{\NN}^{E_S(G)}$,
declaring $\ve_i = 0$ for $i \in S$.
\par
The simplification in the graph is accounted for by a more complex
Hurwitz number at the special vertex. We extend the definition
of Hurwitz numbers with $2$-stabilization by 
\bas
\Cov_2(\bfw,\{\mu_i\}_{i \in S}, \nu) \=
\bigl\{(\pi: S \to \PP^1,\sigma_0)\, :  \, \deg(\pi) = 
\sum w_i\,, \pi^{-1}(0) = \bfw, \\
%& \phantom{\= \bigl\{\pi:} \,\,
\pi^{-1}(1) = [\nu,2^{(\deg(\pi)-|\nu)/2}]\,, 
\pi^{-1}(a_i) = [\mu_i] \, (i \in S),\,\,\ \pi^{-1}(\infty) = [2^{\deg(\pi)/2}] \bigr\}
\eas
for some points $a_i \not\in \{0,1,\infty\}$, and set 
\be  \label{eq:defA2mugen}
A_2(\bfw,\{\mu_i\}_{i \in S},\nu) \= \sum_{\pi \in \Cov_2(\bfw,\{\mu_{i}\}_{i \in S},\nu)}
\frac{1}{\Aut(\pi)}\,.
\ee
The same proof as above now yields the following proposition.
\par
\begin{Prop} \label{prop:corrvariant}
There is a bijective correspondence between 
\begin{itemize}
\item[i)] flat surfaces $(X,q)$ with a covering $p:X\to \B$ of degree~$2d$ and
profile~$\Hmu$ of the pillow~$\B$ without unramified components and
with $q=\pi^*q_\B$, and
\item[ii)] isomorphism classes of tuples
$(G_S, (w_e,h_e, t_e)_{e \in E(G_S)}, (\pi_v)_{v\in V(G_S)})$ as
in Proposition~\ref{prop:corr}, with $\pi_0 \in \Cov'_2(\bfw_0,\nu)$
replaced by $\pi_0 \in \Cov'_2(\bfw_0,\{\mu_i\}_{i \in S},\nu)$,
up to the action of the group $\Aut(\Gamma)$ of automorphisms of the
labeled graph~$\Gamma$. 
\end{itemize}
\end{Prop}
\par

%%%%%%%%%%%%%%%%%%%%%%%%%%%%%%%%%%%%%%%%%%%%%%%%%%%%%%%%%%%%%%%%%%%%
\section{Fermionic Fock space and its balanced subspace}
\label{sec:vertex}
%%%%%%%%%%%%%%%%%%%%%%%%%%%%%%%%%%%%%%%%%%%%%%%%%%%%%%%%%%%%%%%%%%%%

In this section we briefly recall the necessary background
material about fermionic Fock space and the balanced subspace
for the evaluation of the $w$-brackets that compute the generating
functions of pillowcase covers. (See also \cite{rzsurvey}, \cite{OkPand},
\cite{eopillow},  for this formalism)
The new result here is
Theorem~\ref{thm:gnugens} stating that
the generalized shifted symmetric functions $f_\ell$ and $g_\nu$
are rich enough to generate the algebra $\ol{\Lambda}$.
\par
Recall the definition of the normalized characters
$f_{\mu}(\lambda) \= \frakz_{\mu} \chi^\lambda(\mu)/\dim \chi^\lambda$,
where  $\frakz(\mu) =  \prod_{m=1}^{\infty} m^{r_m(\mu)} 
\prod_{m=1}^\infty r_m(\mu)! 
\=  \prod_{i=1}^{\ell(\mu)} \mu_i \prod_{m=1}^\infty r_m(\mu)!$
denotes the order of the centralizer of the partition 
$\mu = 1^{r_1}2^{r_2}3^{r_3}\cdots $.
We also write $f_\ell$ for the special case that $\sigma$ is
a $\ell$-cycle. The {\em algebra of shifted symmetric polynomials} is
defined as  $\Lambda^* = \varprojlim \Lambda^*(n)$, where $\Lambda^*(n)$
is the algebra of symmetric polynomials in the $n$ variables
\hbox{$\l_1-1$,\,}$\ldots,\l_n-n$. The functions 
\be \label{eq:def2pk}  \widetilde{P}_\ell(\l) \= \sum_{i=1}^\infty
\left( (\l_i -i +\h)^\ell - (-i + \h)^\ell \right)  \quad \text{and}\quad
\widetilde{P}_\mu \=\prod_i \widetilde{P}_{\mu_i} \ee
belong to $\Lambda^*$. We add constant terms corresponding to
regularizations to these functions to obtain
\be \label{eq:smallpk}
p_\ell(\lambda) \= P_\ell(\lambda) \+ (1-2^{-\ell})\,\zeta(-\ell).
\ee
Here $\ell! \beta_{\ell+1} = (1-2^{-\ell})\zeta(-\ell)$ with~$\beta_k$
defined by $B(z) \,:=\, \frac{z/2}{\sinh(z/2)} = \sum_{k=0}^\infty \b_k\,z^k$.
Recall the first basic structure result.
\par
\begin{Thm}[\cite{KerOls}] \label{thm:KO}
The algebra $\Lambda^*$ is freely generated by all the $p_\ell$ (or
equivalently, by the $P_\ell$) with~$\ell\ge1$.
The functions $f_\mu$ belong to $\Lambda^*$. More precisely, 
as $\mu$ ranges over all partitions, these  functions $f_\mu$ form
a basis of~$\Lambda^*$.
\end{Thm}
\par
Let $f:\P\to\QQ$ be an arbitrary function on the set $\P$ of all
partitions. That is, we define (following~\cite{eopillow})
the {\em $w$-brackets} 
\be \label{defqbrac} 
\sbw f \= \frac{\sum_{\l\in\P} w(\l) f(\l)\,q^{|\l|}}{\sum_{\l\in\P} w(\l)
q^{|\l|}}\;\,\in\;\QQ[[q]]\,,
\ee
where the difference to the $q$-brackets used to discuss torus
coverings is the weight function
\be \label{eq:defw}
   \sqrt{w(\lambda)}\=\frac{\dim(\lambda)}{|\lambda|!}
   f_{2, \dots, 2}(\lambda)^2, \quad w(\lambda) \= \sqrt{w(\lambda)}^2\,.
   \ee
The main reason for introducing $w$-brackets is the
expression 
\be \label{eq:Covwbracket}
N'(\Hmu) = \sbw{g_\nu f_{\mu_5}\cdots f_{\mu_{n+4}}}
\ee
for the connected Hurwitz numbers. This follows directly from
the classical Burnside formula (see \cite[Section~2]{eo}).
\par
The algebra $\Lambda^*$ is enlarged to the algebra 
\be
\ol{\Lambda} \= \QQ[p_\ell, \olp_k (k,\ell \geq 0)]
\ee
of {\em shifted symmetric quasi-polynomials}, where 
\ba
\olp_k(\lambda) \= \sum_{i \geq 0} \Bigl((-1)^{\lambda_i-i+1}
(\lambda_i -i + \tfrac12)^\ell
- (-1)^{-i+1}(-i + \tfrac12)^\ell \Bigr) \,+\, \gamma_k\,,
\ea
and where the constants $\gamma_i$ are zero for $i$~odd, $\gamma_0 = 1/2$, 
$\gamma_2 = -1/8$, $\gamma_4 = 5/32$ 
and in general defined by the expansion
$ C(z) = 1/(e^{z/2} + e^{-z/2}) = \sum_{k \geq 0}  \gamma_k z^k/k!$.
We provide the algebra $\ol{\Lambda}$ with a grading by defining 
the generators to have
\bes \we(p_\ell) \= \ell + 1 \quad \text{and} \quad \we(\olp_k) \= k\,. \ees
\par
The main reason for introducing $\Lambda^*$ is the following result.
The reason for introducing $g_\nu$ will become clear
by~\ref{eq:Covwbracket} in Section~\ref{sec:HurwitzGraph}.
In this section $\nu$ is always a partition consisting of an even number
of odd parts.
\par
\begin{Thm}[{\cite[Theorem~2]{eopillow}}] \label{thm:gnuinbarL}
There is a function $g_\nu \in \ol{\Lambda}$ of (mixed) weight 
less or equal to $|\nu|/2$ such that
\be \label{eq:defgnu}
g_\nu(\lambda) = \frac{f_{(\nu,2,2,\ldots)}(\lambda)}
{f_{(2,2,\ldots)}(\lambda)}\, \quad \text{for $\lambda$ balanced}.
\ee
\end{Thm}
\par
Our goal is the following converse, for which we define
$\we(g_\nu) = |\nu|/2$.
\par
\begin{Thm} \label{thm:gnugens}
The elements $g_\nu$ generate $\ol{\Lambda}$ as a graded $\Lambda^*$-module
i.e.\ the subspace of $\ol{\Lambda}$ of weight less or equal to~$n$ is
generated by expressions $h g_\nu$ for $h \in \Lambda^*$ with
$\we(h) + \we(g_\nu) \leq n$.
\end{Thm}
\par
Of course, the elements $g_\nu$ do not form a basis as
there are many more $g_\nu$ than products of $\olp_k$ for a given weight. 
We recall the main steps of the proof of Theorem~\ref{thm:gnuinbarL},
since we need them for Theorem~\ref{thm:gnugens}.
\par

%%%%%%%%%%%%%%%%
\subsection{Fermionic Fock space}
%%%%%%%%%%%%%%%%
Let $\Lambda_0^{\frac{\infty}2}V$ be the {\em charge zero subspace of the
  half-infinite wedge} or {\em Fermionic Fock space} over the
countably-infinite-dimensional vector space $V$.
We denote the basis elements of~$V$ by underlined half-integers.
A orthonormal basis of $\Lambda^{\frac{\infty}2}V$ is given by the elements
\bes
v_\lambda  \= \ul{\xi_1} \wedge \ul{\xi_2} \wedge \ul{\xi_3} \cdots, 
\qquad \xi_i = \lambda_i -i + \tfrac12\,.
\ees
indexed by partitions $\lambda = (\lambda_1 \geq \lambda_2 \geq \cdots)$.
The basic operators on the half-infinite wedge is for any $k \in \ZZ + \tfrac12$
the {\em creation operator} $\psi_k(v) = \ul{k} \wedge v$ and its
adjoint $\psi^*$, the {\em annihilation operator}. For any function~$f$ on
the real line we define the (unregularized) operators
\bes \widetilde{\cEE}_k[f] \= \sum_{m \in \ZZ + \tfrac12} f(m) \colon\! 
\psi_{m-k} \psi^*_m \colon,
\ees
where the colons denote the normally ordered product.
We use the convention that~$x$ is the default variable on the real line.
Consequently, if $T$ is a term in~$x$ (typically a polynomial, or a character
like $(-1)^x$ times a polynomial) we write $\widetilde{\cEE}_k(T)$ as
shorthand for $\widetilde{\cEE}_k[x \mapsto T(x)]$. The regularized
operators for exponential arguments are defined by
\bes \cEE_k[e^{zx}] \,:=\, \cEE_k[x \mapsto e^{zx}]
\= \widetilde{\cEE}_k[e^{zx}] + \delta_{0,k}  \frac{B(z)}z
\ees
and in the presence of a character $(-1)^x$ by
\bes \cEE_k[e^{(z+\pi i)x}] \,:=\, \cEE_k[x \mapsto e^{(z+\pi i)x}] \=
\widetilde{\cEE}_k[x \mapsto e^{(z+\pi i)x}]  - i\delta_{k,0} C(z)\,,
\ees
where
\bes
 C(z) \,:= \,
\frac{1}{e^{z/2} + e^{-z/2}} \= \sum_{k \geq 0}  \gamma_k \frac{z^k}{k!}\,.
\ees
We frequently need three special cases of these operators. First
\bes \alpha_{-n} \= \cEE_{-n}[1] \= \sum_{m \in \ZZ + \tfrac12} \colon\!
\psi_{m+n} \psi^*_m \colon\,,
\ees
whose adjoint is denoted by $\alpha_{n} = \alpha_{-n}^*$. The
\emph{Murnaghan-Nakayama rule} says that
\be \label{eq:chiviaalpha}
\prod_i \alpha_{-\mu_i} v_\emptyset \= \sum_{\lambda} \chi^\lambda(\mu)v_\lambda\,.
\ee
Second, the expansion of the (regularized) formal power series 
\bes
 \cEE_0(z) := \cEE_o[x \mapsto e^{zx}] \= 
 \frac1z + \sum_{\ell \geq 1} \,\cPP_\ell \,\frac{z^\ell}{\ell!}
 \ees
and the expansion of 
\bes
i\cEE_o[x \mapsto e^{(z+\pi i)x}]  \= \sum_{\ell \geq 0} \frac{z^k}{k!} 
\ol{\cPP}_k
\ees
gives operators $\cPP_\ell$ and $\ol{\cPP}_k$ with the property
\bes
\cPP_\ell v_\lambda \= p_\ell(\lambda) v_\lambda , \quad
\ol{\cPP}_k v_\lambda \= \olp_k(\lambda) v_\lambda\,.
\ees
Finally, note that the unregularized $\cEE_0(z)$-operator admits the
useful formula
\be \label{eq:CEEexpr}
\widetilde{\cEE}_0(z)  \= [y^0] \psi(e^zy) \psi^*(y)\,,
\ee 
where the interior expression can be checked by the
commutator lemma for vertex operators (\cite[Lemma~7.1]{rzsurvey}
or \cite[Section~14]{KacInf}) to be
\be \label{eq:psiviaextract}
\psi(xy) \psi^*(y) \= \frac{1}{x^{1/2} - x^{-1/2}}
\exp\Bigl(\sum_{n>0} \frac{(xy)^n - y^n}{n} \alpha_{-n} \Bigr)
\exp\Bigl(\sum_{n>0} \frac{y^{-n} - (xy)^{-n}}{n} \alpha_{n} \Bigr).
\ee 

%%%%%%%%%%%%%%%%
\subsection{The balanced subspace}
%%%%%%%%%%%%%%%%

Note that the definition of the algebra~$\Lambda$ excludes the
operator~$p_0$, the charge operator, since it is equal to zero
on~$\Lambda_0^{\frac{\infty}2}V$. Similarly the definition of the
algebra~$\ol{\Lambda}$ excludes the operator
$$\olp_0(\lambda) \ = \frac 12 \+ 
\sum_{i \geq 0} \Bigl((-1)^{\lambda_i-i+1} - (-1)^{-i+1} \Bigr)\,.$$ 
A partition~$\lambda$ is called {\em balanced} if among the $\lambda_i-i+1$
for $\lambda_i \geq 0$ there are as many odd as even numbers, i.e.\ if
and only if $\olp_0(\lambda) \ = \tfrac 12$.
Every partition~$\lambda$ determines (by sorting the $\lambda_i -i +1$
into even and odd) two partitions~$\alpha$ and~$\beta$, 
called the {\em $2$-quotients}, such that 
$$
\{\lambda_i-i+\tfrac12 \} \=  \{2(\alpha_i-i+\tfrac12) + \olp_0(\lambda)\}
\cup  \{2(\beta_i-i+\tfrac12) - \olp_0(\lambda)\}\,.$$
We let $\Lambda^{\rm bal}V$ denote the {\em balanced subspace} 
of $\Lambda_0^{\frac{\infty}2}V$, i.e.\ the subspace spanned by
the $v_\lambda$ for~$\lambda$~balanced. It inherits from 
$\Lambda_0^{\frac{\infty}2}V$ the grading by eigenspace of the energy
operator, i.e.\
$\Lambda^{\rm bal}V = \oplus_{d \geq 0, {\rm even}} \Lambda^{\rm bal}V_d$.
We use the shorthand notation
\bes 
|[\rho; \ol{\rho}]\rangle \= \frac{1}{\frakz(\rho)\frakz({\ol{\rho}})}
\prod_i \alpha_{-\rho_i} \prod_j \olal_{-\ol{\rho_j}} \, v_\emptyset
\ees
for the following reason (\cite{eo}).
\par
\begin{Prop} \label{prop:brhobasis}
For $\rho = (\rho_i)$ and $\ol{\rho} = (\ol{\rho_i})$ running over all
partitions with entries in $2\ZZ$, the elements $|[\rho; \ol{\rho}]\rangle$
form an orthogonal basis of $\Lambda^{\rm bal}V$.
\end{Prop}
\par
\begin{proof}[Proof of Theorem~\ref{thm:gnuinbarL}]
By~\eqref{eq:chiviaalpha} the content of the theorem is that the
orthogonal projection of $|[\nu, 2^{d-|\nu|/2}; \emptyset]\rangle$ to
the balanced subspace $\Lambda^{\rm bal}V$ is a linear combination
of the projections of
$|\prod_i {\cPP}_{\mu_i} \prod \ol{\cPP}_{\ol{\mu_i}}\, 2^d\rangle$
with $\mu = (\mu_i)_{i \geq 1}$ and $\ol{\mu} = (\ol{\mu_i})_{i \geq 1}$
partitions with $\we(\cPP_\mu) + \we(\ol{\cPP}_{\ol{\mu}}) \leq |\nu|/2$
with coefficients independent of~$d$. For this
purpose one calculates using the commutation laws of the
vertex operators that on the one hand
\be \label{eq:nuagainstbarrho}
\bigl\langle [\rho; \ol{\rho}]\, | \,[\nu, 2^{d-|\nu|/2}; \emptyset]
\bigr\rangle
\= \frac{2^{\ell(\nu) - \ell(\ol{\rho})}}{2^{d-|\nu|/2} (d-|\nu|/2)!
\frakz(\nu)\frakz(\ol{\rho})} C(\nu,\ol{\rho}), \quad \text{if} \quad
\rho = 2^{d-|\nu|/2}\,,
\ee
where $C(\nu,\ol{\rho})$ is the number of ways to assemble the parts
of~$\ol{\rho}$ from the parts of~$\nu$, and zero otherwise. In
particular~$|\nu| = |\ol{\rho}|$ for~\eqref{eq:nuagainstbarrho}
to be non-zero. The squared norms of the element
$|[\rho; \ol{\rho}] \rangle$ for $\rho$ and $\ol{\rho}$ having even parts
only is equal to $1/\frakz(\rho)\frakz(\ol{\rho})$. In particular
the scalar product~\eqref{eq:nuagainstbarrho} divided by
$|||[\rho; \ol{\rho}] \rangle ||^2$ is independent of~$d$.
\par
On the other hand one computes using the commutation laws of the
vertex operators that the brackets 
\be \label{eq:nuagainstmumubar}
D_{[\rho;\ol{\rho}], (\mu,\ol{\mu})}(d) \= \bigl\langle [(\rho, 2^{(d-|\rho|)/2};
\ol{\rho}]\, |\,
\prod_i {\cPP}_{\mu_i} \prod \ol{\cPP}_{\ol{\mu_i}}\,|\, 2^d\bigr\rangle /
||\,[(\rho,2^{(d-|\rho|)/2}); \ol{\rho}] \rangle ||^2
\ee
for $\rho$ an partition with only even parts of length different from two
are non-zero only if
\bes \label{eq:weupperbound}
\Delta(\rho,\ol{\rho},\mu,\ol{\mu})
:= \we(\rho) + \we(\ol{\rho}) - \we({\mu}) + \we(\ol{\mu}) \geq 0\,.
\ees
Since an additional factor $\cPP_1$ in~\eqref{eq:nuagainstmumubar}
gives an additional factor $(d-\tfrac1{24})$, we first consider
$D_{[\rho;\ol{\rho}], (\mu,\ol{\mu})}(d)$ with $\mu$ without a part equal to one.
Then, if~\eqref{eq:weupperbound} is attained, the scalar product
$D_{[\rho;\ol{\rho}], (\mu,\ol{\mu})}(d)$ is non-zero if and only
if $\mu = \rho/2$ and $\ol{\mu} = \ol{\rho}/2$. In general,
$D_{[\rho;\ol{\rho}], (\mu,\ol{\mu})}(d)$
is a polynomial of degree $\Delta(\rho,\ol{\rho},\mu,\ol{\mu})/2$.
This also implies that (for fixed weight of $[\rho; \ol{\rho}]$)
the matrix $D_{[\rho;\ol{\rho}], (\mu,\ol{\mu})}(d)$ (with entries
in $\QQ[d]$) is an invertible matrix~$D$, in fact block triangular.
\par
Using these facts we can write
\be \label{eq:gnuformula}
g_\nu = \sum_{|\ol{\rho}| = |\nu|} \frac{2^{\ell(\nu)-\ell(\ol{\rho})}}{\frakz(\nu)}
C({\nu,\ol{\rho}}) \sum_{\mu, \ol{\mu}}
(D^{-1})_{[\emptyset;\ol{\rho}], (\mu,\ol{\mu})}(p_1 + \tfrac1{24})
p_\mu \ol{p}_{\ol{\mu}}
\ee
where the sum all $(\mu, \ol{\mu})$ with
$|\mu| + \ell(\mu) + |\ol{\mu}| \leq \nu/2$.
\end{proof}
\par
\begin{Lemma} \label{le:nurhomat}
For every fixed $d$ the matrix $C(\nu,\ol{\rho})$,
where $\nu$ is a partition of~$2d$ consisting only of odd parts
and $\ol{\rho}$ is a partition of~$2d$ consisting only of even parts,
has full rank equal to $\P(d)$.
\end{Lemma}
\par
\begin{proof} We order the rows $\ol{\rho}$ lexicographically
and consider the submatrix with columns $\nu$ consisting of the
partitions
$$\nu(\rho) \= (\ol{\rho}_1 -1,1,\ol{\rho}_2-1,1,
\ldots,\ol{\rho}_n-1,1)$$
formed from $ \ol{\rho} = (\ol{\rho}_1 \geq \cdots \geq \ol{\rho}_n)$.
Since $C(\nu(\ol{\rho}),\ol{\rho}') \neq 0$ if and only if
$\ol{\rho}' \leq \ol{\rho}$ lexicographically, the claim follows.
\end{proof}
\par
\begin{proof}[Proof of Theorem~\ref{thm:gnugens}]
Since the matrix $D_{[\rho;\ol{\rho}], (\mu,\ol{\mu})}(d)$ is invertible,
it suffices to prove by induction on the weight that the operators
of contraction against
\bes
b_{\rho,\ol{\rho}} \=
\frac{1}{|| \langle \,[(\rho,2^{(d-|\rho|)/2}); \ol{\rho}] ||^2}
\bigl\langle [(\rho, 2^{(d-|\rho|)/2};
\ol{\rho}]\, |\, 
\ees
are in the $\Lambda^*$-module generated by the $g_\nu$. For $\rho = \emptyset$
this follows from Lemma~\ref{le:nurhomat}, in fact those
$b_{\emptyset, \ol{\rho}}$ can be spanned by $g_\nu$ with constant coefficients.
For $\rho \neq \emptyset$ we use the expression of $b_{\rho,\ol{\rho}}$ 
as linear combination of $p_\mu \olp_{\ol{\mu}}$. The terms with
$\we(p_\mu \olp_{\ol{\mu}}) = \we(b_{\rho,\ol{\rho}})$ are either
$p_{\rho/2} \olp_{\ol{\rho}/2}$ or involve a factor of $p_1^j$ for
some $j>0$ by the properties of the matrix
$D_{[\rho;\ol{\rho}], (\mu,\ol{\mu})}(d)$. Consequently, these terms
are generated by $g_\nu$ as $\Lambda^*$-module by induction hypothesis
and the extra factor $p_1^j$ does not alter this fact.
For the terms with smaller weight the induction hypothesis applies directly.
\end{proof}

%%%%%%%%%%%%%%%%%%%%%%%%%%%%%%%%%%%%%%%%%%%%%%%%%%%%%%%%%%%%%%%%%%%%
\section{Hurwitz numbers and graph sums} \label{sec:HurwitzGraph}
%%%%%%%%%%%%%%%%%%%%%%%%%%%%%%%%%%%%%%%%%%%%%%%%%%%%%%%%%%%%%%%%%%%%

The goal of this section is to use the correspondence theorems to express
any $w$-bracket in terms of auxiliary brackets that directly reflect
the graph sums of the correspondence theorems. The precise form of the
goal, Theorem~\ref{thm:bracketgraphsum}, will involve in the
auxiliary brackets only arguments for which the $A'(\cdot)$-functions will
later be proven to be polynomial.
\par
We first define for any function~$F$ on partitions
\be \label{eq:AwithF}
A(\bfw^-,\bfw^+,F) \= \frac{1}{\prod_i w^-_i\prod_i w^+_i } 
\sum_{|\lambda| = d} \chi_{\bfw^-}^\lambda
\chi_{\bfw^+}^\lambda  F(\lambda)
\ee
and we define the connected variant, denoted by $A'(\bfw^-,\bfw^+,F)$ by
the usual inclusion-exclusion formula (e.g.\ \cite[Equation~(17)
or~(25)]{GMab}).
The reason for this definition is that on one hand
the triple Hurwitz number introduced in~\eqref{eq:defAmu} can be
written using the Burnside Lemma (see e.g.\ \cite[Section~2]{GMab}) as
\bes
A'(\bfw^-,\bfw^+,\mu) \= A'(\bfw^-,\bfw^+,f_\mu)\,.
\ees
On the other hand, we will use that the function with completed
cycles argument
\bes
\ol{A}'(\bfw^-,\bfw^+,\mu) \,:=\, A'\Bigl(\bfw^-,\bfw^+,\frac{P_\mu}
{\prod \mu_i}\Bigr)\,.
\ees
is a polynomial of even degree for $\mu = (\mu_1)$ being a partition
consisting of a single part and for $\mu_1 +1-\ell(\bfw^-)- \ell(\bfw^+)$
even (\cite{SSZ}, rephrased as \cite[Theorem~4.1]{GMab}).
\par
\medskip
A new feature of pillowcase covers is the use
of the one-variable analog
\be \label{eq:A2withF}
A_2(\bfw,F) \= \frac{1}{\prod_i w_i } 
\sum_{|\lambda| = d} \sqrt{w(\lambda)}   F(\lambda)\,,
%(f_{(2,2,\ldots}^\lambda)^2 \chi_{\bfw}^\lambda
\ee
where the second variable has been replaced by the character
for the fixed partition $(2,\ldots,2)$. We define the
connected version $A_2'(\bfw,F)$ by the usual inclusion-exclusion formula.
Again, the reason for this definition is two-fold. By the Burnside
formula the simple Hurwitz numbers with $2$-stabilization
introduced in~\eqref{eq:defA2mu} and generalized in~\eqref{eq:defA2mugen}
can be written as
\be
A_2'(\bfw,\{\mu_i\}_{i \in S},\nu) \= A_2'(\bfw,g_\nu\prod_{i \in S} f_{\mu_i} )\,.
\ee
We study polynomiality properties of $A'_2$ for suitable
$F = \prod \ol{p}_{k_i} \in \ol{\Lambda}$ in detail
in Section~\ref{sec:QPoly}.
\par
\medskip
Let $\Hmu$ be a profile as specified in Section~\ref{sec:PillowHurwiz}.
We decompose the Hurwitz number~$N'(\Hmu)$ according to the
contribution of the global graphs, i.e.\ we write
$$  N'(\Hmu) \= \frac{1}{|\Aut(\Gamma)|} \, \sum_{\Gamma}
N'(\Gamma, \Hmu) \,,$$
where the sum is over all (not necessarily connected) labeled 
graphs~$\Gamma$ with $n = |\Hmu|$ non-special vertices and
possibly a special vertex and where $\Aut(\Gamma)$
are the automorphisms of the graph~$\Gamma$ that respect the
vertex labeling. (Note that $\Gamma$ has neither
a labeling nor an orientation on the edges.) Following the
results in the correspondence theorem we define an
{\em admissible orientation~$G$ of~$\Gamma$} (symbolically
written as $G \in \Gamma$) to be an orientation of the half-edges
of~$\Gamma$ such that all the half-edges at the special vertex~$0$
(if it exists) are outward-pointing. Now the following proposition
is an immediate consequence of the correspondence theorem
Proposition~\ref{prop:corr}.
\par
\begin{Prop} \label{prop:N'exp}
The contributions of individual labeled graphs to $N'(\Hmu)$ can
be expressed in terms of triple Hurwitz numbers as
\be \label{eq:NpGamma1} 
N'(\Gamma, \Hmu) \= \sum_{G \in \Gamma } {N}'(G,\Hmu)\,, 
\ee
where 
 \be \label{eq:NpGamma2}
{N}'(G,\Hmu) \= \!\!\sum_{h\in \widetilde{\NN}^{E(G)}, \atop w\in \Z_+^{E(G)}}\prod_{e\in E(G)}w_eq^{h_e w_e}
\cdot A_2'(\bfw_0,\nu) 
\prod_{v\in V(G) \setminus \{0\}}A'(\bfw_v^-,\bfw_v^+,\mu_v)\,\, \delta(v)
\ee
where $V(G)^* = V(G)\setminus \{0\}$  and where
\be \label{eq:NotDelta}
\delta(v) = \delta\bigl(\sum_{i\in e_+(v)} w_i^+ - \!\sum_{i\in e_-(v)} \!w_i^-\bigr).\ee
\end{Prop}
\par
\medskip
We formalize the type of expression appearing in the previous proposition
by defining {\em auxiliary brackets}
\be \label{eq:defauxbrack}
[F_1,\ldots,F_n;F_0] =  \sum_{\Gamma} \, [F_1,\ldots,F_n;F_0]_\Gamma\,,\quad
[F_1,\ldots,F_n;F_0]_\Gamma = \sum_{G \in \Gamma } \, [F_1,\ldots,F_n;F_0]_G,
\ee
where the sum is over all labeled graphs~$\Gamma$ with $n$~vertices and
over all admissible orientations, respectively, and where 
\bes
[F_1,\ldots,F_n; F_0]_G = \sum_{h\in \widetilde{\NN}^{E(G)}, \atop w\in \Z_+^{E(G)}}\prod_{i\in E(G)}w_iq^{h_iw_i}\cdot A_2'(\bfw_0,F_0)
\prod_{v\in V(G)^*} A'(w_v^-,w_v^+,F_{\#v})\,\, \delta(v)\,.
\ees
Here $\# v$ denotes the label of the vertex~$v$. This notation is designed
so that Proposition~\ref{prop:N'exp} can be restated as
\be \label{eq:fasbracket}
\bw{f_{\mu_1}\cdots f_{\mu_n} g_\nu} \= [f_{\mu_1},\ldots,f_{\mu_n};g_{\nu}]\,.
\ee
More generally, by verbatim the same proof, Proposition~\ref{prop:corrvariant}
can be restated as the generalization
\be \label{eq:fasbracketS}
\bw{f_{\mu_1}\cdots f_{\mu_n} g_\nu} \= [\prod_{i \not\in S} f_{\mu_i};
  \prod_{i \in S} f_{\mu_i}g_{\nu}]
\ee
for any subset $S \subseteq \{1,\ldots,n\}$.
\par
\medskip
We are now ready to formulate the goal of this section in detail.
\par
\begin{Thm} \label{thm:bracketgraphsum} 
The $w$-bracket of any element in $\ol{\Lambda}$ can be expressed
as a finite linear combination of the auxiliary brackets, i.e.\ 
for every $\bfl = (\ell_1,\ldots,\ell_n)$ and every 
$\bfk = (\ell_1,\ldots,k_m)$ there exist $c_{(\bft,\bfs)} \in \QQ$ 
(depending on $(\bfl,\bfk)$)
such that
\be \label{eq:wbrackASaux}
\bw{ \prod_{j=1}^n p_{\ell_j } \, \prod_{i=1}^m \olp_{k_i}}  
\= \sum_{(\bft,\bfs)}  c_{(\bft,\bfs)} \Bigl[p_{t_1},\ldots,p_{t_{\ell(\bft)}};
\prod_{i=1}^{\ell(\bfs)} \olp_{s_i}\Bigr]\,,
\ee
where the sum is over all $(\bft,\bfs)$ with 
$\sum_j (t_j+1) + \sum_i s_i \leq \sum_j (\ell_j +1) + \sum_i k_i$.
\end{Thm}
\par
\begin{proof} The proof is by induction on the weight $w = \sum_{i=1}^m k_i$
of the $\olp_{k_i}$-part of the bracket, the case of weight zero
being trivial (no special vertex, i.e.\ as in the abelian case.)
\par
By Theorem~\ref{thm:gnugens} we can write the left hand side
of~\eqref{eq:wbrackASaux} as a linear combination of 
$\bw{f_{\mu_1}\cdots f_{\mu_n} g_\nu}$ with $\we(g_\nu) \leq w$.
By~\eqref{eq:fasbracket} each such summand is equal to
\be \label{eq:fgback}
[f_{\mu_1},\ldots,f_{\mu_n};g_{\nu}]
\= \sum_{\bfb,\bfa}  [f_{\mu_1},\ldots,f_{\mu_n};
\prod_{j \geq 1} p_{b_j } \, \prod_{i \geq 1} \olp_{a_i}] \,,
\ee
where the sum is over all partitions~$\bfa$ and~$\bfb$,
by Theorem~\ref{thm:gnuinbarL}. In the summands where $\bfb$
is the empty partition, we replace $f_{\mu_i}$ by a linear
combination of products of $p_\ell$ thanks to Theorem~\ref{thm:KO}
and these contributions are of the required form of the right hand
side of~\eqref{eq:wbrackASaux}. In all the summands with $\bfb$
non-empty we use the converse base change of Theorem~\ref{thm:KO}
to write the product of $p_{b_j}$ as a linear combination of a
product of $f_{\mu_j}$. We can now use~\eqref{eq:fasbracketS} from
right to left to express all the terms as a sum $w$-brackets with
$\olp_{k_i}$-part of weight $w - \sum_{j} b_j+1$. Since $\bfb$ is
non-empty, we conclude thanks to the induction hypothesis.
\end{proof}

%%%%%%%%%%%%%%%%%%%%%%%%%%%%%%%%%%%%%%%%%%%%%%%%%%%%%%%%%%%%%%%%%%%%%%%%%%%%%%%%%
\section{Constant coefficients of 
quasi-elliptic functions} 
\label{sec:QMFProp}
%%%%%%%%%%%%%%%%%%%%%%%%%%%%%%%%%%%%%%%%%%%%%%%%%%%%%%%%%%%%%%%%%%%%%%%%%%%%%%%%

In this section we consider the constant coefficient
(in $z_1,\ldots,z_n$) of a function that is quasi-elliptic in these
variables, has a globally a quasimodular transformation behavior and
poles at most at two-torsion translates of the coordinate axes
and diagonals. We show in Theorem~\ref{thm:coeff0QM} that this constant
coefficient is indeed a quasimodular form for the subgroup $\Gamma(2)$
of $\SL\Z$.
\par
%%%%%%%%%%%%%%%%%%%%%
\subsection{Quasimodular forms} \label{sec:QMF}
%%%%%%%%%%%%%%%%%%%%%

A {\em quasimodular form} for the cofinite Fuchsian group $\Gamma \subset 
\SL\RR$ of weight~$k$ is a function $f: \HH \to \CC$ that is holomorphic
on $\HH$ and the cusps of $\Gamma$ and such that there exists and integer~$p$ and
holomorphic functions $f_i: \HH \to \CC$ such that 
$$ (c\tau + d)^{-k} f\Bigl(\frac{a\tau +b }{c\tau +d}\Bigr) \=  
\sum_{i=0}^p f_i(\tau) \Bigl(\frac{c}{c\tau + d} \Bigr)^i\,\quad 
\text{for all} \quad \sm abcd \in \Gamma\,.$$ 
\par
The smallest integer~$p$ with the above property is called the {\em depth}
of the quasimodular form.  By definition, 
quasimodular forms of depth zero are simply modular forms.
The basic examples of quasimodular forms are the Eisenstein series defined by
\[G_{2k}(\tau) \= \frac{(2k-1)!}{2(2\pi i)^{2k}}\sum_{(m,n)\in\Z^2\setminus\{(0,0\}}
\frac{1}{(m+n\tau)^{2k}} \= -\frac{B_{2k}}{4k}+\sum_{n=1}^{\infty}\sigma_{2k-1}(n)q^n\,.\]
Here  $B_l$ is the Bernoulli number, $\sigma_l$ is the divisor sum function
and $q=e^{2\pi i\tau}$. For $k\geq 2$ these are modular forms, while for $k=2$ the
Eisenstein series 
\[G_2(\tau)=-\frac{1}{24}+\sum_{n=1}^{\infty}\sigma_1(n)q^n\] 
is a quasimodular form of weight~$2$ and depth~$1$ for $\SL\ZZ$. By
\cite{KanZag} we can write the ring of quasimodular forms for any group
$\Gamma$ with $\Stab_\infty(\Gamma) = \pm \langle \sm1101 \rangle$
as $\QM(\Gamma) =
\CC[G_2]\otimes M(\Gamma)$ in terms of the ring of modular forms $M(\Gamma)$.
%We will be mainly interested in the congruence groups $\Gamma_0(n)$ as well
%as in
%\be
%\Gamma^C_0(2) \=
%\begin{pmatrix}1 & 0\\ 0 & 2\end{pmatrix}\Gamma_0(2)
%\begin{pmatrix}1 & 0\\0 & 2\end{pmatrix}^{-1} \,\subset \, \Gamma_0(4) 
%\ee
We will be mainly interested in the congruence groups $\Gamma_0(2)$ and
\bes
\Gamma(2)=\begin{pmatrix}2 & 0\\ 0 & 1\end{pmatrix}\Gamma_0(4)
\begin{pmatrix}2 & 0\\0 & 1\end{pmatrix}^{-1} \,\subset \, \Gamma_0(2).
\ees
Since $M(\Gamma_0(2))$ is freely generated by $G_2^{\odd}(\tau) = G_2(\tau)
-2G_2(2\tau)$ and $G_4(2\tau)$ by the transformation formula (e.g.\
\cite[Proposition~4.2.1]{DiaShu}) and the usual dimension formula for
modular forms, we deduce that
\be \label{eq:G02gens}
\QM(\Gamma_0(2)) \,\cong\, \CC[G_2(\tau), G_2(2\tau), G_4(2\tau)]
\ee
is a polynomial ring.
% and also that 
%\be \label{eq:G02Cgens}
%\QM(\Gamma_0^C(2)) \cong \CC[G_2(2\tau), G_2(4\tau), G_4(4\tau)]\,
%\ee
%is a polynomial ring.
 For $\Gamma_0(4)$ we restrict our attention to
the subring of even weight quasimodular forms. Since $M_{2*}(\Gamma_0(4))$
is freely generated by $G_2^{\odd}(\tau)$ and $G_2^{\odd}(2\tau)$ (again
using the transformation formula together with the isomorphism $\Gamma_0(4) 
\cong \Gamma(2)$ given by conjugation with ${\rm diag}(2,1)$),
we deduce that 
\be \label{eq:G2gens}
\QM(\Gamma(2)) \,\cong\, \CC[G_2(\tau/2), G_2(\tau), G_2(2\tau)]\,.
\ee
\par
We use the notation $q=e^{2\pi i \tau}$, hence $q^{1/2}=e^{\pi i \tau}$.
Note that a typical element $\QM(\Gamma(2))$ has a Fourier expansion
in~$q^{1/2}$. The following observation allows us to prove
quasimodularity by the larger group $\Gamma_0(2)$.
\par
\begin{Lemma} \label{le:EisG2} For all $k \in \NN$ 
the Eisenstein series $G_{2k}(\tau/2)$, $G_{2k}(\tau)$ and $G_{2k}(2\tau)$
are quasimodular forms for $\Gamma(2)$. Moreover, any 
even weight quasimodular form for $\Gamma(2)$ whose Fourier expansion 
is a series in~$q$ is in fact a quasimodular form for $\Gamma_0(2)$.
\end{Lemma}
\par
\begin{proof}
The second statement follows immediately from~\eqref{eq:G02gens}
and~\eqref{eq:G2gens}.
\end{proof}

%%%%%%%%%%%%%%%%%%%
\subsection{Quasimodular forms as constant coefficients of  
quasi-elliptic functions}  \label{sec:intEF}
%%%%%%%%%%%%%%%%%%%

We are now ready to state the first main criterion for quasimodularity, 
involving the constant coefficients of some quasi-elliptic functions introduced below. We start with a general remark
on the domains where the expansions are valid.
Suppose that the meromorphic function $f(z_1, z_2,\dots, z_n; \tau)$ 
is periodic under $z_j \mapsto z_j + 1$ for each~$j$ and under 
$\tau \mapsto \tau +1$. We can then write $f(z_1, z_2,\dots, z_o;\tau) 
= \ol{f}(\zeta_1, \ldots, \zeta_n,q)$ where $\zeta_j = e^{2\pi i z_j}$ as above. 
For any permutation $\pi \in S_n$ on we fix the domain
\bes
 \Omega_{\pi} \= |q^{1/2}|<|\zeta_{\pi(i)}|<|\zeta_{\pi(i+1)}|<1 \quad 
\text{for all~$i=1,\ldots,n-1$}\,.
\ees
On such a domain the {\em constant term  with respect to all the $\zeta_i$}
is well-defined. It can be expressed as integral
$$[\zeta_n^0, \ldots,\zeta_1^0]_{\pi} \, \ol{f} \= \frac{1}{(2\pi i)^n}\,
\oint_{\gamma_n}\dots \oint_{\gamma_1}f(z_1, \dots, z_n;\tau)dz_1\dots dz_n
$$
along  the integration paths 
\[ \gamma_j: [0,1]\to\C, \quad  \; t\mapsto iy_j+t\,,\]
where  $0\leq y_{\pi(1)}< y_{\pi(2)}<\dots y_{\pi(n)}< 1/2$. 
We call these our {\em standard integration paths} for the permutation~$\pi$.
If the domain $\Omega_\pi$ is clear from the context we also write
$[\zeta^0]$ or $[\zeta_n^0, \ldots,\zeta_1^0]$ as shorthand for the 
coefficient extraction $[\zeta_n^0, \ldots,\zeta_1^0]_\pi$.
\par
Let $\Delta = \Delta_{\tau}$ be the operator on meromorphic functions defined by
$$ \Delta(f)(z) \= f(z+\tau) - f(z)\,.$$
A meromorphic function~$f$ is called {\em quasi-elliptic} (for the lattice $\ZZ+\tau\ZZ$)
if $f(z+1) = f(z)$
and if there exists some integer~$e$ such that $\Delta^e(f)$ is elliptic.
The minimal such~$e$ is called the {\em order (of quasi-ellipticity)} of~$f$.
\par
We say that a meromorphic function $f: \CC^n \times \HH \to \CC$ 
is quasi-elliptic, if it is quasi-elliptic in each of the first~$n$ variables.
For such a function we write $\bfe = (e_1,\ldots,e_n)$ for the tuple of
orders of quasi-ellipticity in the $n$~variables. Consequently, 
a quasi-elliptic function of order $(0,\ldots,0)$ is simply an
elliptic function.
\par
We write $\Delta_i$ for the operator $\Delta$ acting on the $i$-th variable.
Note that these operators~$\Delta_i$ commute.
Let $T=\{0,1/2, \tau/2, (1+\tau)/2\}$ be the set of 2-torsion points.
\par
The functions we want to take constant coefficients of belong to
the space  in the following definition. It is similar to the
quasi-elliptic quasimodular forms used in \cite[Definition~5.5]{GMab}.
The difference consists of allowing $2$-torsion translates for
the poles and requiring a modular transformation law for a smaller group.
We do not decorate our new definition of $\cQQ_{n,\bfe}^{(k)}$ by
an extra symbol~$2$ to avoid overloading notation.
\par
\begin{Defi}
We define for $n \geq 0$, $k \geq 0$ and $\bfe \geq 0$
the vector space of $\cQQ_{n,\bfe}^{(k)}$ of {\em quasi-elliptic quasimodular forms} for $\Gamma(2)$
to be the space of meromorphic functions~$f$ on 
$\CC^n \times \HH$ in the variables $(z_1,\ldots,z_n; \tau)$ that 
\begin{enumerate}
\item[i)] have poles on $\CC^n$ at most at the $\ZZ + \tau\ZZ$-
translates of the diagonals $z_i = z_j$, $z_i=-z_j$ and the 2-torsion points $z_i\in T$,
\item[ii)] that are quasi-elliptic of order $\bfe$, and
\item[iii)]  that are quasimodular of weight~$k$ for $\Gamma(2)$, i.e.\
$f$~is holomorphic in $\tau$ on $\HH \cup \infty$ and 
there exists some $p \geq 0$  and functions 
$f_i(z_1,\ldots,z_n; \tau)$ that are holomorphic in~$\tau$ and meromorphic in
the~$z_i$ such that
$$ (c\tau + d)^{-k} f\Bigl(\frac{z_1}{c\tau +d},\ldots, \frac{z_n}{c\tau +d}; 
\frac{a\tau +b }{c\tau +d}\Bigr) \=  
\sum_{i=0}^p f_i(z_1,\ldots,z_n; \tau) \Bigl(\frac{c}{c\tau + d} \Bigr)^i$$ 
for all $\sm abcd \in \Gamma(2)$.
\end{enumerate}
\end{Defi}
\par
Examples of such quasi-elliptic quasimodular forms will be constructed
from the {\em propagator}, the shift $P(z;\tau) =
\tfrac{1}{(2\pi i)^2}\wp(z;\tau) + 2G_2(\tau)$ of the Weierstrass
$\wp$-function, and from the shift $Z(z;\tau) = -\zeta(z;\tau)/2\pi i
+2G_2(\tau)2\pi iz$ of the Weierstra\ss\ $\zeta$-function.
The reason for this shift, as well as the Fourier and Laurent series
expansion of these functions is summarized in \cite[Section~5.2]{GMab}.
In particular we will need 
\bas
P_{even}(z;\tau)&\= 2P(2z;2\tau) \quad \text{and}\\
P_{odd}(z;\tau)&\=P(z;\tau)-P_{even}(z;\tau)\,.
\eas
\par
\begin{Prop} \label{prop:PZinQring}
  The functions  $P^{(k)}(z_i-a;\tau)$ where $a\in T$, and each of the functions
$P^{(k)}(2z_i; \tau)$, $P^{(k)}(2z_i; 2\tau)$, $P_{even}^{(k)}(z_i;\tau)$,
$P_{odd}^{(k)}(z_i;\tau)$, $ P^{(k)}(z_i-z_j; \tau)$ and $P^{(k)}(z_i+z_j; \tau)$ 
belong to $\cQQ_{n,\bf 0}^{(k+2)}$.
\par
The functions $Z(z_i-a;\tau)$ for $a\in T$ belong to $\cQQ_{n,\bfe_i}^{(1)}$,
where $\bfe_i=(0,\dots,0,1,0,\dots,0)$, and the functions $Z(z_i-z_j;\tau)$ and
$Z(z_i+z_j;\tau)$ belong to $\cQQ_{n,\bfe_i+\bfe_j}^{(1)}$. 
\end{Prop}
\begin{proof}
Since $P$ is an elliptic meromorphic function with poles at $\Z+\tau\Z$, and
quasimodular in the sense of iii) of weight 2 for $\SL\Z$, more precisely
\bes
(c\tau + d)^{-2} P\Bigl(\frac{z}{c\tau +d}; 
\frac{a\tau +b }{c\tau +d}\Bigr) \=  
P(z; \tau)+ \frac{c}{c\tau + d} \,, 
\ees
we deduce easily the result for the functions derived from $P$. In fact,
the functions $P^{(k)}(z_i-1/2;\tau)$ (hence also $P^{(k)}(2z_i; 2\tau)$,
$P_{even}^{(k)}(z_i;\tau)$ and $P_{odd}^{(k)}(z_i;\tau)$)  are quasimodular for
the bigger group $\Gamma_0(2)$. Moreover, the functions $P^{(k)}(2z_i; \tau)$
are quasimodular for the full group $\SL\Z$.
We proceed similarly for the functions derived from $Z$, which is quasi-elliptic
of order 1, quasimodular of weight one and depth one with $Z_1(z;\tau)=z$.
\end{proof}
\par
\begin{Prop} \label{prop:Qring}
The direct sum
$$ \cQQ_n \= \bigoplus_{k \geq 0} \cQQ_n^{(k)}, \quad \text{where} 
\quad \cQQ_n^{(k)} \= \bigoplus_{\bfe \geq 0}\, \cQQ_{n,\bfe}^{(k)} \,,$$
is a graded ring. The derivatives $\partial/\partial z_i$ map 
$\cQQ_n^{(k)}$ to $\cQQ_n^{(k+1)}$ for all $i=1,\ldots,n$
and the derivative $D_q = q\tfrac{\partial}{\partial q}$ maps
$\cQQ_n^{(k)}$ to $\cQQ_n^{(k+2)}$.
\par
For all $i,j \in \{1,\ldots,n\}$, $i\neq j$, the functions
\be \label{eq:LZ2}
L(z_i;\tau) \= -\frac{1}{2}Z^2(z_i;\tau) + 
\frac{1}{2}P(z_i;\tau)-G_2(\tau)+\frac{1}{12}
\ee as well as
$L(2z_i; \tau), L(2z_i; 2\tau), L(z_i-z_j;\tau)$ and $L(z_i+z_j;\tau)$
belong to $\cQQ_n^{(0)} \oplus \cQQ_n^{(2)}$.
\par
\end{Prop}
\par
\begin{proof} The proof is similar to the $\SL\Z$-case
 (cf.~\cite[Proposition~5.6]{GMab}).
\end{proof}
\par From now on we omit the variable $\tau$ in the notation, if not necessary.
\par
\begin{Prop} \label{prop:Naddbasis}
The vector space $\cQQ_n$ is (additively) generated as $\cQQ_{n-1}$-module  
by the functions $Z^e(z_n-a)$ and $Z^e(z_n-a) P^{(m)}(z_n-a)$
for $a\in T$ together with $Z^e(z_n+z_j)$ and
$Z^e(z_n+z_j) P^{(m)}(z_n+z_j)$,$Z^e(z_n-z_j)$, $Z^e(z_n-z_j) P^{(m)}(z_n-z_j)$
for $j=1,\ldots,n-1$  and for all $e \geq 0$ and $m \geq 0$.
\par
More precisely, if~$f \in \cQQ_n^{(k)}$ then we can write
\bas f(z_1,\ldots,z_n) &\= \sum_{a\in T}\left(\sum_{e,m} A_{e,m,j} Z^e(z_n-a) P^{(m)} (z_n-a)  + 
\sum_{e} B_{a,e} Z^e(z_n-a)\right)\\
&+  \sum_{e,m,i} C_{e,m,j} Z^e(z_n+z_i) P^{(m)} (z_n+z_i)  + 
\sum_{e,i} D_{e,i} Z^e(z_n+z_i)\\
&+\sum_{e,m,i} E_{e,m,i} Z^e(z_n-z_i) P^{(m)} (z_n-z_i)  + 
\sum_{e,i} F_{e,i} Z^e(z_n-z_i) + G
\eas
with $A_{a,e,m}, C_{e,m,i},  E_{e,m,i} \in \cQQ_{n-1}^{(k-e-m+2)}$, $B_{a,e}, D_{e,i},F_{e,i}  \in \cQQ_{n-1}^{(k-e)}$
and $G \in \cQQ_{n-1}^{(k)}$.
\end{Prop}
\par
\begin{proof}
For every~$n$ we argue inductively on the order $e = \min_{j \geq 0}
\{\Delta_n^j(f)\,\, \text{elliptic}\}$ of quasi-ellipticity with respect to 
the last variable. 
Suppose, without loss of generality, that $f \in \cQQ_n^{(k)}$ is homogeneous
of weight~$k$. 
\par
We first treat the case $e=0$. We show that we can write 
\[f=\sum_{a\in T} \left(\sum A_{m,a}P^{(m)}_a + B_{a}Z_{a}\right) + \sum C_{m,i}P^{(m)}_{i,+}+D_i Z_{i,+}+ \sum E_{m,i}P^{(m)}_{i,-}+ F_{i,j}Z_{i,j}+G\]
with $A_{m,a}, C_{m,i}, E_{m_i}\in \cQQ_{n-1}^{(k-m-2)}$, $B_a, D_i, F_{i,j}\in \cQQ_{n-1}^{(k-1)}$ and $G\in\cQQ_{n-1}^{(k)}$,
where $P^{(m)}_{j,a} = P^{(m)}(z_n+a)$ for $a\in T$, 
$P^{(m)}_{i,-} = P^{(m)}(z_n-z_i)$,$P^{(m)}_{i,+} = P^{(m)}(z_n+z_i)$,
for all $m \geq 0$ and all $i=1,\ldots,n-1$ 
$Z_{a}=Z(z_n+a)-Z(z_n+z_{n-1})+Z(z_{n-1})$ for $a\in T$, $Z_{i,+}=Z(z_n+z_i)-Z(z_n-z_i)-2Z(z_i)$, for all $i=1,\ldots,n-1$,
$Z_{i,j} = Z(z_n-z_i) - Z(z_n-z_j)  + Z(z_i-z_j)$ for all $1 \leq i < j \leq n-1$. By
Proposition~\ref{prop:PZinQring}, these functions are clearly in $\cQQ_{n,\bf 0}$.
\par
We proceed by induction on the pole orders, first along the divisors $z_n-a$, then along the divisors $z_n+z_i$, then along the divisors $z_n-z_i$. The rest of the proof is then totally similar to \cite[Proposition~5.4]{GMab}. The residue theorem ensures that we can eliminate the last poles with the functions $Z_{i,j}$ to end the procedure. 
\par
For the case $e>0$, the proof is a straightforward adaptation of the proof of \cite[Proposition~5.7]{GMab}.
\end{proof}
\par
Using this additive basis we can now prove the main result.
\par
\begin{Thm} \label{thm:coeff0QM}
For any permutation $\pi$ the constant term with respect to the domain 
$\Omega_\pi$ of a function in $\cQQ_{n}^{(k)}$ is a quasimodular form for $\Gamma(2)$
of mixed weight $\leq k$.
\end{Thm}
\par
\begin{proof} Again, the proof is exactly the same as in \cite[Theorem~5.8]{GMab}: we reduce the problem to the computation of $[\zeta^0] Z^e(z-a)$ by some transformations preserving the weight of the quasimodular form. The last step of the proof is isolated in the statement of the next proposition.
\end{proof}
\par
\begin{Prop} \label{prop:Zpowint}
The constant coefficient $[\zeta^0] Z^e(z-a)$ for $a\in T$ is
a quasimodular form for $\SL\Z$ of mixed weight less or equal to $e$.
\end{Prop}
\par
\begin{proof} 
  Since $Z(z)$ is $1$-periodic, we clearly have  $[\zeta^0]Z^e(z-1/2)
  =[\zeta^0]Z^e(z)$ and these coefficients are quasimodular forms for $\SL\Z$ of mixed weight less or equal to~$e$, by \cite[Proposition~5.9]{GMab}.
Similarly, $[\zeta^0]Z^e(z-\tau/2 -1/2)=[\zeta^0]Z^e(z-\tau/2)$ so we just
have to compute $[\zeta^0]Z^e(z-\tau/2)$.
\par
Using again the $1$-periodicity of $Z$, we obtain for all $\ell$
\ba \label{eq:none}
[\zeta^0]Z^\ell(z-\tau/2) &\= \int_{\gamma_1} Z^\ell(z-\tau/2)dz
\= \int_{\gamma_2} Z^\ell(z)dz \= \int_{\gamma_3}Z^\ell(-z)dz \\
&\= [\zeta^0]Z^\ell(-z),
\ea
where the integration paths $\gamma_1$, $\gamma_2$, $\gamma_3$ are described
in Figure~\ref{cap:intpaths}.
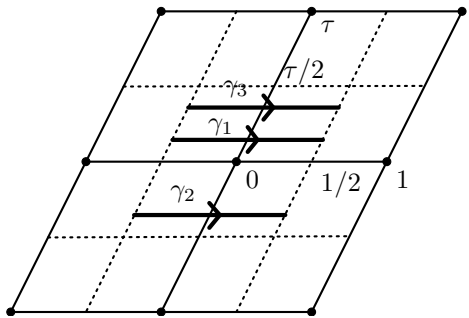
\begin{figure}
\begin{center}
\begin{tikzpicture}[line cap=round,line join=round,>=triangle 45,x=1.0cm,y=1.0cm]
\clip(-3.66,-0.05) rectangle (3.7,4.45);
\draw [line width=0.8pt] (-1.,0.)-- (1.,4.);
\draw [line width=0.8pt] (-3.,0.)-- (-1.,4.);
\draw [line width=0.8pt] (-3.,0.)-- (1.,0.);
\draw [line width=0.8pt] (-2.,2.)-- (2.,2.);
\draw [line width=0.8pt] (-1.,4.)-- (3.,4.);
\draw [line width=0.8pt] (1.,0.)-- (3.,4.);
\draw [line width=0.8pt,dotted] (-1.5,3.)-- (2.504,3.008);
\draw [line width=0.8pt,dotted] (0.,4.)-- (-2.,0.);
\draw [line width=0.8pt,dotted] (0.,0.)-- (2.,4.);
\draw [line width=0.8pt,dotted] (-2.504,0.992)-- (1.504,1.008);
\draw [line width=1.6pt] (-0.85,2.29)-- (1.15,2.29);
\draw [line width=1.6pt] (0.28,2.29) -- (0.15,2.13);
\draw [line width=1.6pt] (0.28,2.29) -- (0.15,2.46);
\draw [line width=1.6pt] (-1.36,1.29)-- (0.65,1.29);
\draw [line width=1.6pt] (-0.22,1.29) -- (-0.36,1.12);
\draw [line width=1.6pt] (-0.22,1.29) -- (-0.36,1.45);
\draw [line width=1.6pt] (1.36,2.72)-- (-0.64,2.72);
\draw [line width=1.6pt] (0.49,2.72) -- (0.36,2.88);
\draw [line width=1.6pt] (0.49,2.72) -- (0.36,2.56);
\draw (0,2) node[anchor=north west] {$0$};
\draw (1,4) node[anchor=north west] {$\tau$};
\draw (1,2) node[anchor=north west] {$1/2$};
\draw (2,2) node[anchor=north west] {$1$};
\draw (0.5,3.5) node[anchor=north west] {$\tau/2$};
\draw (-0.5,2.7) node[anchor=north west] {$\gamma_1$};
\draw (-1,1.8) node[anchor=north west] {$\gamma_2$};
\draw (-0.3,3.2) node[anchor=north west] {$\gamma_3$};
\begin{scriptsize}
\draw [fill=black] (0.,2.) circle (1.5pt);
\draw [fill=black] (-1.,0.) circle (1.5pt);
\draw [fill=black] (1.,4.) circle (1.5pt);
\draw [fill=black] (2.,2.) circle (1.5pt);
\draw [fill=black] (3.,4.) circle (1.5pt);
\draw [fill=black] (1.,0.) circle (1.5pt);
\draw [fill=black] (-3.,0.) circle (1.5pt);
\draw [fill=black] (-2.,2.) circle (1.5pt);
\draw [fill=black] (-1.,4.) circle (1.5pt);
\end{scriptsize}
\end{tikzpicture}
\end{center}
\caption{Integration paths to evaluate $[\zeta^0]Z^e$} \label{cap:intpaths}
\end{figure}
Since $Z$ is odd, the constant coefficients $[\zeta^0]Z^\ell(z-\tau/2)$ are
then given by $(-1)^\ell[\zeta^0]Z^\ell$ (see \cite{GMab} for some explicit values).
\end{proof}
\par
Note that the proof of Theorem~\ref{thm:coeff0QM} provides an effective
algorithm to compute constant coefficients of quasi-elliptic functions. Applications
of this algorithm and explicit computations of quasi-modular forms are
detailed in Section~\ref{sec:Example}.

%%%%%%%%%%%%%%%%%%%%%%%%%%%%%%%%%%%%%%%%%%%%%%%%%%%%%%%%%%%%%%%%%%%%%%%%%%%%%
\section{Quasimodularity of graph sums} 
\label{sec:GraphSumisQM}
%%%%%%%%%%%%%%%%%%%%%%%%%%%%%%%%%%%%%%%%%%%%%%%%%%%%%%%%%%%%%%%%%%%%%%%%%%%%%

The goal of this section is to show the quasimodularity of
graph sums of the form~\eqref{eq:SGammam} below. The motivation for
considering these sums will become apparent in comparison with
the quasi-polynomiality theorem in Section~\ref{sec:QPoly}.
We encourage the reader to look at Section~\ref{sec:Example} simultaneously
with this one, we hope that all notations will become transparent on the example
that we treat in detail.
\par
We will show quasimodularity for graphs that arise as global graphs
as in Section~\ref{sec:globalgraph} with the following extra decoration.
\par
A {\em global graph with distinguished edges} is a graph with vertices with
the labels $1,\dots,n$ and possibly a special vertex~$0$ and a subset~$E^+(\Gamma)$
of $E(\Gamma)$ of distinguished edges such that no extremity of an edge
in~$E^+(\Gamma)$ is the vertex~$0$. We let $V^*(\Gamma) = V(\Gamma) \setminus \{0\}$,
we let $E^0(\Gamma)$ be the edges adjacent to~$0$. Finally we define
$E^*(\Gamma) = E(\Gamma) \setminus E^0(\Gamma)$ and $E^-(\Gamma) = E^*(\Gamma)
\setminus E^+(\Gamma)$. An {\em admissible orientation $G$ of $(\G,E^+)$} is
an orientation of the half-edges of~$\G$, such that
\begin{itemize}
\item all half-edges adjacent to the vertex $0$ are outgoing, and 
\item the orientations of the two half-edges are consistent on marked edges,
and inconsistent on the other edges that are not adjacent  to $v_0$. 
\end{itemize}
We write $G \in (\G,E^+)$ for the specification of an admissible orientation.
Obviously every admissible orientation $G$ of $\G$ (in the sense of
Section~\ref{sec:HurwitzGraph})
is admissible for $(\G,E^+)$ for a uniquely determined subset~$E^+ \subset E(\G)$.
We define the {\em set of parity conditions} to be $\PC(\Gamma) = \{0,1\}^{E^0(\Gamma)}$.
It specifies a congruence class mod~$2$ for the width of each edge adjacent to
the vertex~$0$.
\par
We consider here for fixed $\mm = (m_1,\ldots,m_{E(\Gamma)})$ and fixed $E^+(\Gamma)$
the graph sums
\be \label{eq:SGammam} {S}(\Gamma,E^+,\mm, \pari)\=\sum_{G \in (\Gamma,E^+)}{S}(G,\mm, \pari)\ee
over all admissible orientations $G$ of the half edges of $\Gamma$, where for
$\pari \in \PC(\Gamma)$
\be \label{eq:SGm} {S}(G,\mm, \pari) \=\sum_{h\in\widetilde\NN^{E(G)} \atop \bfw_*\in \N_{>0}^{E(G)^*}}
\sum_{\bfw_0 \in \NN_{>0}^{E^0(\Gamma)} \atop  \bfw_0 \cong \pari\,{\rm mod}\, 2}
\prod_{i\in E(G)}w_i^{m_i+1}q^{h_iw_i}\prod_{v\in V(G)^*}\delta(v)\,. \ee
Here $\widetilde\NN^{E(G)}$ is the height space introduced
in~\eqref{eq:heightspace} and~$\delta(v)$ is as in~\eqref{eq:NotDelta}.
The goal of this section is to show the quasimodularity of these graph sums.
\par
\begin{Thm} \label{thm:QMforS}
For a fixed tuple of non-negative even integers
$\mm = (m_1,\ldots,m_{|E(\Gamma)|})$ and for any
$\pari \in \PC(\Gamma)$ the graph sums $S(\Gamma, E^+,\mm,\pari)$ are
quasimodular forms for the group $\Gamma_0(2)$  of mixed 
weight at most $k(\mm):= \sum_{i} (m_i +2)$.
\end{Thm}
\par
The following is a first simplification step for the computation.
Recall the definition of the offsets $a_e$ from Figure~\ref{cap:HEorient}.
\par
\begin{Lemma}\label{lem:heightspace} Replacing the height space  
$\widetilde\NN^{E(G)}$  by
\be \label{eq:newheightspace}
\widehat{\NN}^{E(G)} \= \Bigl\{ (h'_e)_{e \in E(\Gamma)} \,:\,
\begin{cases} 
h'_e   \in \NN_{>0} \,  &\text{if} \quad e \in\EOl(G)\cup\EPl(G) \\
h'_e   \in \NN_{\geq a_e} \, &\text{otherwise}
%if} \quad e \in E_\ell(G) \setminus(\EOl(G)\cup\EPl(G)) \\
%h_e - \Delta(e)  \in \NN_{\geq a_e} &\text{if} \quad e \in E_\ell(G) \setminus E_\ell(G)
\end{cases} \Bigr\}\,,
\ee
does not change the total sum $S(G,\mm,\pari)$. 
\end{Lemma}
\par
\begin{proof} We apply the linear change of variables
$h'_e=h_e-\Delta(e)$ with $\Delta(e)$ as in the line below~\eqref{eq:heightspace}. 
This maps $\widetilde\NN^{E(G)}$ bijectively onto $\widehat\NN^{E(G)}$. For notational
convenience we set $\Delta(e)=0$ for the remaining edges. Each summand of~\eqref{eq:SGm}
for fixed $(w_1,\dots,  w_{|E(G)|})$ is multiplied under the variable change by
\bes q^{\sum_{e}\Delta(e)w_e}\prod_{v\in V(G)^*}\delta(v) \= 
q^{\sum_v \varepsilon_v\left(\sum_{i\in e^-(v)} w_i - \sum_{i\in e^+(v)} w_i\right)}\prod_{v\in V(G)^*}
\delta(v) \=1\,.
\ees
This implies the claim.
\end{proof}

%%%%%%%%%%%%%%%%%%%%%%%
\subsection{The reduced graph}
%%%%%%%%%%%%%%%%%%%%%%%%

We will simplify the graph sums by isolating the contribution from the loop
edges $\EOl(G)\cup\EPl(G)$. The {\em reduced graph} $\OG$ is obtained from $\G$
by deleting those loops, i.e.\ the loops adjacent to the vertex $v_0$
and the loops among the distinguished edges.
\par
\begin{Lemma} \label{le:loopsplit}
The graph sums $S(\G, \mm, \pari)$ factor as 
 \[S(\G, E^+,\mm, \pari) \=S_{loops}(\G, \mm, \pari)\,S(\OG, E^+,\mm, \pari)\]
 where 
\[S_{loops}(\G, \mm, \pari)\=\sum_{h\in\N_{>0}^{\EOl(\G)\cup\EPl(\G)}} \sum_{\bfw\in\Z_+^{\EPl(\G)}}
\sum_{\bfw_0 \in \NN_{>0}^{E^0(\Gamma)} \atop  \bfw_0 \cong \pari\,{\rm mod}\, 2}
\prod_{i\in\EOl(\G)\cup\EPl(\G)}w_i^{m_i+1}q^{h_iw_i}\,.\]
 \end{Lemma}
\par 
\begin{proof} The length constraints $\delta(v)$ for $v \in V(G)^*$ are unchanged
under removing a loop edge $e \in E^+(G)$ that contributes equally to both incoming
and outgoing weight. The length parameters at the vertex~$0$ is always unconstrained.
This proves the factorization we claim.
\end{proof}
\par
\begin{Lemma} \label{le:loopOK}
If $\mm$ is even, $S_{loops}(\G, \mm, \pari)$ is a quasimodular form for~$\Gamma_0(2)$
of mixed weight $k(\mm)$.
\end{Lemma}
\par
\begin{proof}
The graph sum  $S_{loops}(\G, \mm, \pari)$ is a product of 
\bas
S_m &\= \sum_{w,h,=1}^{\infty} w^{m+1}q^{2wh} &&\= G_{m+2}(q)-G_{m+2}(0)\\
 S_{m, even} &\= \sum_{w,h=1}^{\infty} (2w)^{m+1}q^{(2w)h} &&\=2^{m+1}( G_{m+2}(q^2)-G_{m+2}(0))\\
 S_{m,odd} &\= \sum_{w,h,=1}^{\infty} (2w-1)^{m+1}q^{(2w-1)h} &&\= S_m-S_{m, even}\,.
\eas
All the right hand sides are quasimodular forms for $\Gamma_0(2)$
by \eqref{eq:G02gens}.
\end{proof}

%%%%%%%%%%%%%%%%%%%%%%%%%%%%%%%%%
\subsection{Contour integrals} \label{sec:ContourInt}
%%%%%%%%%%%%%%%%%%%%%%%%%%%%%%%%%

We now write the sum of the reduced graph as contour integral
of suitable derivatives of the following variants of the propagator.
Let $P_{even}(z;\tau) = 2P(2z;2\tau)$ and $P_{odd}(z;\tau) = P(z;\tau)-P_{even}(z;\tau)$. We also use
$P_i$ for $i \in \ZZ/2$ to refer to these two functions. For a reduced
graph~$\OG$, for $\mm$ even, and for a given parity condition $\pari \in \PC(\Gamma)$,
define
\ba \label{eq:defPG}
P_{\OG, E^+,\mm, \pari}(\zz)\=\prod_{i\in \EO(\OG)}
P_{\pari_i}^{(m_i)}(z_{v_+(i)})
\cdot \prod_{i\in \EP(\OG)}
P^{(m_i)}( z_{v_-(i)}-z_{v_+(i)})\cdot\\
\cdot \prod_{i\in \EM(\OG)}
P^{(m_i)}( z_{v_-(i)}+z_{v_+(i)})\,,
\ea
where $v_+(i)$ and $v_-(i)$ are the two ends of the edge $i$, with $v_-$ being the
one of lower index (so for $i\in\EO(\OG)$ necessarily $v_-(i)=0$). Note that since
the graph is reduced, $v_+(i)$ and $v_-(i)$ can be the same, but only if $i\in\EM(\OG)$.
Note that the variable $z_0$ does not appear in the expression~$P_{\OG, \mm, \pari}(\zz)$
at all.
\par
\begin{Prop} \label{prop:S1} For a tuple of non-negative even integers~$\mm_1$
and a parity condition $\pari$ we can express the graph sum as
\ba \label{eq:SasContInt}
{S}(\OG,E^+,\mm, \pari)\= [\zeta_n^0, \dots, \zeta_1^0]\,
P_{\OG, E^+,\mm, \pari}(\zz; \tau)\,,
\ea
where the coefficient extraction is for the expansion on the domain 
 $|q^{1/2}|<|\zeta_i|<|\zeta_{i+1}|<1$ for all~$i$.
\end{Prop}
\par
\begin{proof} The proof is analogous to the proof of \cite[Proposition~6.7]{GMab}
and we suggest to read the two proofs in parallel since we will not reproduce
the bulky main formulas. We rather indicate the main changes.
First note that in the domain specified above the  inequalities
\bes
|q|<|\zeta_i\zeta_j|<1\; \quad \forall i,j, \qquad 
|q|<|q^{1/2}|<|\zeta_i/\zeta_j|<1 \;\quad \forall j>i
\ees
hold, and hence the following Fourier expansions 
\ba \label{eq:PropExpansion}
P_{even}^{(m)}(z;\tau)&\= \sum_{w\geq 1, even} w^{m+1} \left(\zeta^{w}\sum_{h\geq 0} q^{wh} +\zeta^{-w}\sum_{h\geq 1}q^{wh} \right)\\
P_{odd}^{(m)}(z;\tau)&\= \sum_{w\geq 1, odd} w^{m+1} \left(\zeta^{w}\sum_{h\geq 0} q^{wh} +\zeta^{-w}\sum_{h\geq 1}q^{wh} \right)\\
P^{(m)}(z_i+z_j;\tau)&\= \sum_{w\geq 1} w^{m+1} \left((\zeta_i\zeta_j)^{w}\sum_{h\geq 0} q^{wh} +(\zeta_i\zeta_j)^{-w}\sum_{h\geq 1}q^{wh}\right)\; \forall i,j\\
P^{(m)}(z_i-z_j;\tau)&\= \sum_{w\geq 1} w^{m+1} \left((\zeta_i/\zeta_j)^{w}\sum_{h\geq 0} q^{wh} +(\zeta_i/\zeta_j)^{-w}\sum_{h\geq 1}q^{wh}\right)\; \forall j>i
\ea
are valid. For the proof we first consider the factors in~\eqref{eq:defPG} that
involve the edges~$E_1$ adjacent to the vertex~$1$, i.e.\ those involving the
variable~$z_1$. The propagators $P,P_{\rm odd}$ or $P_{\rm even}$ in~\eqref{eq:defPG}
are chosen so that the parity conditions for $w_e$ specified in~\eqref{eq:SGm} hold.
Each of the propagators in~\eqref{eq:PropExpansion} has (for fixed $(w,h)$) two
summands, that we consider as incoming ($\zeta$-exponent~$+w$) or outgoing
($\zeta$-exponent~$-w$). Consequently, expanding the product of propagators
involving the edges in~$E_1$ is a sum over all partitions $E_1 = J_1 \cup K_1$
of the incoming and outgoing terms. The integration with respect to~$z_1$
forces that all contributions vanish except for those where the incoming~$w_e$
are equal to the outgoing~$w_e$. This ensures the appearance of the factor~$\delta(v_1)$
in~\eqref{eq:SGm}. The proof proceeds by similarly considering the vertex~$2$ and
expanding the propagator factors that involve the edges $E_2$ adjacent this
vertex but not already in~$E_2$, which produces a sum over all partitions
$E_2 = J_2 \cup K_2$ according to whether the incoming or outgoing summand
of the propagator has been taken.
\par
The main difference to the abelian case is the consequence of orienting
the half-edges. Suppose that $e$ joins~$v_1$ to~$v_2$. If $e \in E^+(\Gamma)$
then in all admissible orientations $e$ is incoming at $v_1$ and outgoing at $v_2$,
or vice versa. If $e \in J_1$ is incoming at $v_1$, we have to make sure that
the propagator terms have $\zeta_2^{-w}$, i.e.\ we have to use $P^{(m)}(z_1 - z_2)$.
On the other hand if $e \in E^-(\Gamma)$, then in all admissible orientations $e$ is
incoming or outgoing simultaneously at $v_1$ and $v_2$. I.e.\ $\zeta_1$ and
$\zeta_2$ have to appear with the same $w$-exponent, whence the use
of $P^{(m)}(z_1 + z_2)$. The reader can check that this orientation convention
is also consistent for the special vertex~$0$ and that the range of the
sums $h \geq 0$ versus $h \geq 1$ in~\eqref{eq:defPG} is consistent with
the conditions of the height space that appear in the $h$-summation
in~\eqref{eq:SGm}.
\end{proof}
\par
\begin{proof}[Proof of Theorem~\ref{thm:QMforS}] This is now a direct
consequence of Proposition~\ref{le:loopOK} for the loop contribution, of
Proposition~\ref{prop:S1} for the reduced graph and of Theorem~\ref{thm:coeff0QM}
for quasimodularity (for $\Gamma(2)$) of contour integrals, combined with Lemma~\ref{le:EisG2} to get quasimodularity for the bigger group $\Gamma_0(2)$.
\end{proof}

%%%%%%%%%%%%%%%%%%%%%%%%%%%%%%%%%%%%%%%%%%%%%%%%%%%%%%%%%%%%%%%%%%%%%%%%
\section{Quasipolynomiality of $2$-orbifold double Hurwitz numbers} 
\label{sec:QPoly}
%%%%%%%%%%%%%%%%%%%%%%%%%%%%%%%%%%%%%%%%%%%%%%%%%%%%%%%%%%%%%%%%%%%%%%%%

The main result of this section is the quasi-polynomiality of
the simple Hurwitz numbers with $2$-stabilization $A_2'(\bfw,F)$
in the case that $F$ is a product of $\olp_k$. The meaning of
quasi-polynomiality is  that the restriction to a congruence class
mod~$2$ in each variable is a polynomial. The crucial statement
for the quasi-modularity is that these polynomials are global, i.e.\
not piece-wise polynomials depending on a chamber decomposition of
the domain of~$\bfw$. As a first application we combine this with the
correspondence and quasimodularity theorems of the previous section
to give in Corollary~\ref{cor:noname} another proof the Eskin-Okounkov
theorem on the quasimodularity of the number of pillowcase covers. 
\par
%%%%%%%%%%%%%%%%%%%
\subsection{The one-sided pillowcase operator}
%%%%%%%%%%%%%%%%%%%
Our goal here is to write $A_2'(\bfw,F)$ in terms of vertex operators.
For this purpose we define the
\emph{one-sided pillowcase operator}
\bes
\Gamma_{\sqrt{w}} \= \exp\Bigl( \sum_{i>0} \tfrac{\alpha_{-i}^2}{2i}\Bigr)\,.
\ees
\par
\begin{Prop} \label{prop:A2vertex}
The simple Hurwitz numbers with $2$-stabilization can be expressed
using the one-sided pillowcase operator as
\be \label{eq:A2vertex}
A_2(\bfw, F) \=  \cfrac{1}{\prod w_i}\, \Bigl\langle 0 \,|\,
\prod_{i=1}^{\ell(\bfw)}
\alpha_{w_i} \cFF\, \Gamma_{\sqrt{w}} \,|\, 0 \Bigl\rangle
\ee
where $\cFF v_\lambda = F(\lambda) v_\lambda$.
\end{Prop}
\par
\begin{proof} We first observe that
\be \label{eq:sqrtwfirst}
\langle \Gamma_{\sqrt{w}} v_{\emptyset}, v_\lambda \rangle
\= \sum_\nu \chi^\lambda(2\nu)\prod_{i\geq 1} \cfrac{(1/2i)^{r_i(\nu)}}{r_i(\nu)!}\,,
\ee
where $2\nu$ is the partition obtained by repeating twice each row of the
Young diagram of $\nu$, i.e.\ if $\nu = 1^{r_1(\nu)}2^{r_2(\nu)}\cdots$ is
written in terms of the multiplicities of the parts then
$2\nu = 1^{2r_1(\nu)}2^{2r_2(\nu)}\cdots$ . This observation follows from developing
the exponential in $\Gamma_{\sqrt{w}}$ and the Murnaghan-Nakayama rule. It
thus remains to show that
\bes
\sum_\nu \chi_\lambda(2\nu)\prod_{i=1}^{\nu_1} \cfrac{(1/2i)^{r_i(\nu)}}{r_i(\nu)!}
\=f_{2, 2, \dots, 2}(\lambda)^2\cfrac{\dim\lambda}{|\lambda|!}
\=\sqrt{\textrm{\bf w}(\lambda)}
\ees
or equivalently that
\be \label{eq:f2reform}
f_{2, \dots, 2}(\lambda)^2
\=\sum_{\nu} f_{2\nu}(\lambda)\cdot \prod_{i \geq 1} (i^{r_i}(2r_i-1)!!).
\ee
Let $C_2$ be the conjugacy class corresponding to the partition
$(2,\dots, 2)$ and $C(2\nu)$ the conjugacy class of $2\nu$. Define
\[n^{C}_ {C_2^2} \= \# \{ (a,b)\in C_2^2 \,|\, ab\in C\}\] 
and let $e_C=\sum_{g\in C} [g]$ be  the central elements the group ring
$\mathbb Z[G]$. Then  $e_{C_2}e_{C_2}=\sum_C n^{C}_ {C_2^2} e_C$.
Note that $e_C$ acts on any irreducible representation $\lambda$ as
multiplication by the scalar $f_C(\lambda)$. Consequently, the
claim~\eqref{eq:f2reform} is equivalent to 
\[n^{C}_ {C_2^2} \=\begin{cases} \prod_i (i^{r_i}(2r_i-1)!!)
\mbox {\quad  if } C=C(2\nu)\\
0 \mbox{ otherwise\,.}\end{cases}\]
To prove this, note first that the product $\tau = \rho \sigma$ of two
permutations $\rho, \sigma$ in $C_2$ is in $C(2\nu)$. This follows
be recursively proving that $\sigma \tau^{-n}(P) = \tau^{n}(\sigma(P)$,
i.e.\ the cycles starting at $P$ and $\sigma(P)$ have the same length.
To justify the combinatorial factor, assume first that all $r_i=1$.
In order to specify the factorization it is necessary and sufficient
to specify for each~$i$ and for some point in an $i$-cycle
its $\sigma$-image in the other $i$-cycle. The rest of the factorization
is determined by the requirement of profile $(2, \dots, 2)$ and this
initial choice. In the general case $r_i>1$ we moreover have to match
the $2r_i$ cycles of length~$i$ in pairs (which corresponds to the factor
$(2r_i-1)!!$) and make a choice as above for each pair. 
\par
The proposition follows from~\eqref{eq:sqrtwfirst} using again
the Murnaghan-Nakayama rule.
\end{proof}
\par
We next collect some rules to evaluate the right hand side
of~\eqref{eq:A2vertex}. The standard commutation law between creation operators
is
\be\label{eq:comut2}\alpha_n^a\alpha_{-n}^b \=
\sum_{k=0}^{\min(a,b)} \binom{a}{k} \frac{b!}{(b-k)!} n^k \alpha_{-n}^{b-k}
\alpha_n^{a-k}\,.
\ee
Together with $[\alpha_m,\alpha_n] = 0$ for $m \neq -n$ this implies that we
can deal with the $\alpha_{\pm n}$ for all~$n$ separately, i.e.\
\bes
\langle \prod_n f_n(\alpha_n)g_n(\alpha_{-n})v_\emptyset, v_\emptyset\rangle
\= \prod_n\langle f_n(\alpha_n)g_n(\alpha_{-n}) v_\emptyset, v_\emptyset\rangle\,,
\ees
where $f_n(\alpha_n)=\sum c_{j,n}\alpha_n^j$ and $g_n(\alpha_{-n})=
\sum c'_{j,n}\a_{-n}^j$. From the commutation relation 
\bes e^{c_n\alpha_{-n}}e^{d_n\alpha_{n}} \= e^{-nc_nd_n}e^{d_n\alpha_{n}}e^{c_n\alpha_{-n}}\,,
\ees
that follows directly from~\eqref{eq:comut2}, 
we deduce successively the following relations involving the factors
of~$\Gamma_{\sqrt{w}}$. First
\bes
\Pi \,:=\, \prod_{i=1}^{s}\left(e^{c_{i,n}\alpha_{-n}}e^{d_{i,n}\alpha_n}\right)
\= e^{-nA_n}e^{D_n\alpha_n}e^{C_n\alpha_{-n}}
\ees
where $C_n=\sum_{i=1}^s c_{i,n}$, $D_n=\sum_{i=1}^s d_{i,n}$ and
$A_n=\sum_{j=1}^{s} d_{j,n}\sum_{i=1}^{j}c_{i,n}$. Next, from
\bas 
\langle 0|e^{c_n\alpha_{-n}}e^{d_n\alpha_n}e^{\alpha_{-n}^2/2n}|0\rangle & =  e^{nd_n^2/2}\\
\langle 0|\alpha_{n}e^{c_n\alpha_{-n}}e^{d_n\alpha_n}e^{\alpha_{-n}^2/2n}|0\rangle & =  n(d_n+c_n)e^{nd_n^2/2} 
\eas
we obtain
\ba \label{eq:brwithGOp}
\langle 0|\prod_{i=1}^s \left(e^{c_{i,n}\alpha_{-n}}e^{d_{i,n}\alpha_n}\right)e^{\alpha_{-n}^2/2n}|0\rangle & = e^{n\tilde A_n}e^{nD_n^2/2}  \\
 \langle 0|\alpha_{n}\prod_{i=1}^s \left(e^{c_{i,n}\alpha_{-n}}e^{d_{i,n}\alpha_n}\right)e^{\alpha_{-n}^2/2n}|0\rangle & =  n(D_n+C_n)e^{n\widetilde A_n}e^{nD_n^2/2} 
\ea
where $\widetilde A_n= D_nC_n - A_n = \sum_{j=2}^s d_{j,n}\sum_{i=1}^{j-1} c_{i,n}$,
and most generally 
\be \label{eq:genbrwithGOp}
\langle 0|\alpha_{n}^\ell\, \Pi\,
%\left(e^{C_{n}\alpha_{-n}}e^{D_{n}\alpha_n}\right)
e^{\alpha_{-n}^2/2n}|0\rangle \= \sum_{i=0}^{\ell/2}
%OLD, equivalent \frac{\ell!\, n^{\ell-i}(C_n+D_n)^{\ell-2i}  }{(2i)!! \, (\ell-2i)!}
\binom{\ell}{2i} \frac{(2i)!}{2^i(i)!}
n^{\ell-i}(C_n+D_n)^{\ell-2i} e^{n\widetilde A_n}e^{nD_n^2/2}\,.
\ee

%%%%%%%%%%%%%%%%%%%
\subsection{Polynomiality}
%%%%%%%%%%%%%%%%%%%

We can now evaluate the vertex operator expressions and prove the
main result of this section.
\par
\begin{Thm} \label{thm:pbarpoly}
The simple Hurwitz number with $2$-stabilization $A_2'(\bfw,F)$ without
unramified components is a quasi-polynomial if~$F$  is a product of $\olp_k$,
i.e.\ for each coset $\mm = (m_1,\ldots,m_t) \in \{0,1\}^n$ with $\sum m_i$
even there exists a polynomial $R_{F,\mm} \in \QQ[w_1,\ldots,w_t]$ such that
\bes
A_2'(\bfw,F) \= R_{F,\mm}(\bfw) \quad \text{for all}
\quad \bfw \in 2\NN^t + \mm\,.
\ees
\end{Thm}
\par
The proof relies on matching piece-wise polynomial on sectors like
$w_1>w_2$ to form a global polynomial. The parity constraints to match
the piece-wise polynomials do not work out for elements in $\Lambda^*$,
not even for $F = p_\ell$, if there is more than one boundary variable~$w_i$.
In fact for any $w_1,w_2 \in \NN$ with $w_1+w_2$ even
\bes
\frac15\,A_2'((w_1,w_2),p_5) \= \frac78 u^3+ \frac{13}{8}uv^2 - u\,,
\ees
where $u = \min(w_1,w_2)$ and $v = \max(w_1,w_2)$ is only piece-wise
polynomial, while e.g.\ on the coset $\mm = (0,0)$
\bes
\frac14\,A_2'((2w_1,2w_2),\ol{p}_4) \= 10(2w_1)^2 + 10(2w_2)^2 -3
\ees
is globally a polynomial.
\par
The basic source of polynomiality is the following lemma, relevant
for the case of $F = \olp_k$.
\par
\begin{Lemma} \label{eq:1VarPoly}
For each~$k\geq 1$ there is a polynomial $Q_k$ such that for $n \in \NN$
\bes
Q_k(n) \= [y^{-2n}][z^k] D\,, \quad \text{where}\quad
D(y,z) \= \frac{1+y^{-2}e^{-z}}{\sqrt{(1-y^{-2})(1-y^{-2}e^{-2z})}}\,.
\ees
Moreover, $Q_k$ is even for $k$ even and $Q_k$ is odd for $k$ odd
and $Q_k(0)=0$ in both cases.
\end{Lemma}
\par
\begin{proof} We abbreviate $\partial_z = \partial/\partial z$.
Since $[z^k] = k! \partial_z f|_{z=0}$ it suffices
to write $\partial^k_z D(y,z)|_{z=0} =  R_k(y)/(1-y^{-2})^{k+1}$
for some polynomial $R_k(y^{-2})$ of degree $\leq k$ without
constant coefficient. The relation $D(1/y,-z)= D(y,z)$ implies
that $R_k$ is palindromic, i.e.\ in the span of
$y^{-s} + y^{-(2k-s)}$ for $s = 2,4,\ldots,2k-2$. Since by the binomial
theorem
\bes
[y^{-2n}] \frac1{(1-y^{-2})^{k+1}} \= \frac{1}{k!}(n+k-1)\cdots (n+1)n
\ees
agrees with a polynomial for integers $n \geq 1-k$, we obtain the polynomiality
claim. The parity claim follows from $R_k$ being palindromic. 
\end{proof}
\par
\begin{proof}[Proof of Theorem~\ref{thm:pbarpoly}]
We may shift $\olp_k$ by the regularization constant~$\gamma_k$
in order to use~\eqref{eq:CEEexpr} and~\eqref{eq:psiviaextract} and assume that
$F = \prod_{j=1}^s (\olp_{k_j} -\gamma_{k_j})$ for some $k_j$, not necessarily
distinct. Proposition~\ref{prop:A2vertex} now translates into
\be \label{eq:A2eval1}
A_2(\bfw,F) \= [y_1^0\cdots y_s^0][z_1^{k_1} \cdots z_s^{k_s}]
\prod_{n \geq 1} \frac{1}{n^{r_n(\bfw)}} 
\langle 0|\alpha_{n}^{r_n(\bfw)} \Psi_F e^{\alpha_{-n}^2/2n}|0\rangle\,,
\ee
where $r_n(\bfw)$ is the multiplicity of~$n$ in $\bfw$ and where
\bes
\Psi_F = \prod_{j=1}^s  \frac1{e^{z_i/2} + e^{-z_i/2}}
\exp\Bigl( \frac{y_j^n((-e^{z_j})^n -1)}{n} \alpha_{-n} \Bigr)
\exp\Bigl(\frac{y_j^{-n}(1 - (-e^{-z_j})^{n})}{n} \alpha_{n} \Bigr).
\ees
By~\eqref{eq:brwithGOp}, the first factor common to the evaluation of the
brackets~\eqref{eq:A2eval1} for all~$n$ is $e^{nD_n^2/2}$, results in
a product of 
\bes
D_{[j]}\,:=\, 
\exp\Bigl(\sum_{n>0} \frac{y_j^{-2n}(1 - (-e^{-z_j})^{n})^2}{2n} \Bigr)
\=  \frac{1+y_j^{-2}e^{-z_j}}{\sqrt{(1-y_j^{-2})(1-y_j^{-2}e^{-2z_j})}}\,.
\ees
and
\bes
D_{[ij]} \,:=\,  \frac{(1+y_iy_je^{-z_i})(1+y_iy_je^{-z_j})}
{{(1-y_iy_j)(1-y_iy_je^{-z_i-z_j})}}\,.
\ees
The second common factor $e^{n\widetilde A_n}$ for all~$n$ results in a factor
of
\bes
\widetilde{A} \,:=\, \prod_{j=2}^s \prod_{i=1}^{j-1}
\frac
{(1+\frac{y_j}{y_i}e^{-z_i})(1+\frac{y_j}{y_i}e^{z_j})}
{(1- \frac{y_j}{y_i})(1-\frac{y_j}{y_i}e^{z_j-z_i})}
\,. \ees
In order to built an arbitrary covering with boundary lengths~$\bfw$
from a covering without unramified components, we have to choose for
each length~$n$ of the boundary components among the $\ell_n=r_n(\bfw)$
an even number~$2i$ of boundary components that are glued together in pairs
to form cylinders, and the number of such gluing is $\frac{(2i)!}{2^i(i)!}$,
the number of fixed point free involutions. That is, the combinatorial
factor~$\binom{\ell_n}{2i} \frac{(2i)!}{2^i(i)!}$ in front
of the summand in~\eqref{eq:genbrwithGOp} counts precisely these possibilities.
Consequently this formula implies that
\be \label{eq:no76}
A_2'(\bfw,F) \=  [y_1^0\cdots y_s^0][z_1^{k_1} \cdots z_s^{k_s}]
\Biggl(\frac{\widetilde{A} \, \prod_{i<j}D_{[ij]} \prod_j D_{[j]}}
     {\prod_{j} (e^{z_j/2} + e^{-z_j/2})}
\prod_{n: r_n(\bfw) \geq 1}
K_n^{r_n(\bfw)} \Biggr)
\ee
where $K_n = C_n + D_n$ as in the vertex operator manipulations above, i.e.\
\bes
K_n \= \sum_{j=1}^s \frac{y_j^n((-e^{z_j})^n -1)}{n} \,+\,
\frac{y_j^{-n}(1 - (-e^{-z_j})^{n})}{n}
\ees
The claim follows if we can show two statements, first that the
expression~\eqref{eq:no76} is piece-wise polynomial and second that this
expression with each $K_n$ replaced
by $\widetilde{K_n} = n K_n$ is globally a polynomial. The factor
$\prod_{j} (e^{z_j/2} + e^{-z_j/2})$ results just in a shift of $z_j$-degrees
and will be ignored in the sequel. Note that we can write the
last factor in~\eqref{eq:no76} equivalently as
$\prod_{n: r_n(\bfw) \geq 1} K_n^{r_n(\bfw)} = \prod_{i=1}^t K_{w_i}$. 
\par
We start with the case $s=1$, illustrating the main idea. Let $t$
be the number of~$n$ with $r_n(\bfw) \geq 1$, say these are $n_1,\ldots,n_t$.
First, we want to show that
\bes \bfn = (n_1,\ldots,n_t) \,\mapsto\,
[y_1^0] \Bigl([z_1^{k-j}] \prod_{i=1}^t \widetilde{K}_{n_i}\,\cdot \,
[z_1^j] D_{[1]} \Bigr)
\ees
is polynomial in the $n_i$ for each~$j \in [0,k]$ (and zero otherwise). To
evaluate this, we can choose in each $\widetilde{K}_{n_i}$-factor the
$y_i^{n_i}$-term or
the $y_i^{-n_i}$-term and then sum over the contributions of all choices. For
each $\bfdelta  = (\delta_1,\ldots,\delta_t) \in \{\pm 1\}^t$ we consider the
linear form $f_\bfdelta(n_1,\ldots,n_t) = \sum_{i=1}^t \delta_i n_i$. We claim that
already the sum of the contributions of $f_\bfdelta$ and $f_{-\bfdelta}$ is
polynomial, i.e.\ that
\bas
\bfn &\,\mapsto\, [z_1^{k-j}] \prod_{i=1}^t ((-1)^{m_i} e^{\delta_i n_i z_1} -1)
[y_1^{f_\bfdelta(\bfn)}][z_1^j] D_{[1]} \\
&\phantom{\,\mapsto\,} + [z_1^{k-j}] \prod_{i=1}^t ((-1)^{m_i} e^{-\delta_i n_i z_1} -1)
[y_1^{f_{-\bfdelta}(\bfn)}][z_1^j] D_{[1]}
\eas
is the restriction of a polynomial to any collection of natural numbers
$n_i \equiv m_i \mod (2)$ is the fixed coset. If we denote
$f_\bfdelta^+ = \max(0,f_\bfdelta)$, the claim follows from the observation that
$Q_j(\tfrac12 f_\bfdelta^+(\bfn)) + (-1)^j Q_j(\tfrac12 f_{-\bfdelta}^+(\bfn))
= Q_j(\tfrac12 f_\bfdelta(\bfn))$
is globally a polynomial for $\bfn$ in a fixed congruence class and for~$Q_j$
with the parity as in Lemma~\ref{eq:1VarPoly}.
The polynomiality for $s=1$ and the $\widetilde{K_n}$-version follows by
summing up these expressions. 
\par
Second, we argue that $A_2'(\bfw,\bar{p}_k)$ without the additional
factors~$n$ in $\widetilde{K_n}$ is a piece-wise polynomial,
i.e.\ that the polynomial expression obtained previously using $\widetilde{K_n}$
is indeed divisible by~$n$. The divisibility by $n_i$ follows from
adding the contribution of $f_\delta$ and $f_{\delta'}$ where $\delta'$
differs from $\delta$ precisely in the $i$-th digit, since
$(n_i+n_j)^k + (n_i-n_j)^k$ is divisible by~$n_i$ independently of the
parity of~$k$.
\par
For the general case $s \geq 1$ follows along the same lines.
We first prove a generalization of Lemma~\ref{eq:1VarPoly}, stating
that for $\bfk = (k_1,\ldots,k_s)$ there is a polynomial~$Q_{\bfk}$
in~$s$ variables~$g_i$ such that 
\bes
Q_{\bfk}(g_1,\ldots,g_s) \= [y_1^{-2g_1}\cdots y_s^{-2g_s}][z_1^{k_1} \cdots z_s^{k_s}]
\Bigl({\widetilde{A} \, \prod_{i<j}D_{[ij]} \prod_j D_{[j]}}  \Bigr)
\ees
if all $g_i \geq 0$. Moreover this polynomial has the parity
\be \label{eq:Qparity}
Q_{\bfk}(g_1,\ldots, -g_j, \ldots,g_s) \= (-1)^{k_j}
Q_{\bfk}(g_1,\ldots, -g_j, \ldots,g_s) \,.
\ee
To see this, we write
\bas
&\phantom{\=} \partial_{z_1}^{k_1}\cdots\partial_{z_s}^{k_s} \Bigl({\widetilde{A}
\,\prod_{i<j}D_{[ij]} \prod_j D_{[j]}}  \Bigr)|_{z_1=\cdots=z_s=0} \\
& \=
\sum_{\bfl, \bfm,\bfn} R_{\bfk,\bfl,\bfm,\bfn}(y_1,\ldots,y_s) \cdot \prod_i
\frac{1}{(1-y_i^2)^{\ell_i}} \cdot
\prod_{i<j} \frac{1}{(1-y_iy_j)^{m_{ij}}} \cdot
\prod_{i<j}  \frac{1}{(1-y_j/y_i)^{n_{ij}}} 
\eas
for some polynomial $R_{\bfk,\bfl,\bfm,\bfn}$ with the components of
$\bfl,\bfm$ and $\bfn$ bounded in terms of~$\bfk$. Using the binomial
expansion of this expression we see that the coefficient 
$[y_1^{-2g_1}\cdots y_s^{-2g_s}]$ of this expression is a sum of polynomial
expressions in non-negative integers $u_i$, $v_{ij}$ and $w_{ij}$ over the
bounded simplex defined by 
\bes
u_i + \sum_j v_{ij} + \sum_{j<i} n_{ij} -\sum_{j>i} n_{ij} \= g_i, \quad
(i=1,\ldots, s)\,.
\ees
This sum is again a polynomial and this implies the polynomiality claim.
Moreover,
the argument ${\widetilde{A} \, \prod_{i<j}D_{[ij]} \prod_j D_{[j]}} $
in the definition of~$Q_\bfk$ is unchanged under the transformation
$(y_i, z_i) \mapsto (1/y_i, -z_i)$ since each $D_{[i]}$ has this property
and since this transformation swaps $D_{[ij]}$ with the $(i,j)$-factor
of $\widetilde{A}$. As in the case $s=1$ this implies palindromic
numerators and the parity statement~\eqref{eq:Qparity}.
\par
We show similarly that for each tuple~$(j_1,\ldots,j_s)$ separately the
function 
\bes \bfn = (n_1,\ldots,n_t) \,\mapsto\,
[y_1^0] \Bigl([z_1^{k_1-j_1}\cdots z_s^{k_s-j_s}]
\prod_{i=1}^t \widetilde{K}_{n_i}\,\cdot \,
[z_1^{j_1}\cdots z_s^{j_s}] \widetilde{A} \, \prod_{i<j}D_{[ij]} \prod_j D_{[j]}
\Bigr)
\ees
is a polynomial. Evaluation of $\prod_{i=1}^t \widetilde{K}_{n_i}$
now leads to $s$ linear forms $f_1(\bfn),\ldots,f_s(\bfn)$ with
each $n_i$ appearing in exactly one of the $f_j$, and with coefficient~$\pm 1$.
Given one such tuple, the contributions of $\pm f_1,\ldots,\pm f_s$ add
up to a global polynomial thanks to~\eqref{eq:Qparity}, as in the case $s=1$.
Finally adding the contribution of $f_1,\ldots,f_s$ and the linear
form with precisely the sign of~$n_i$ flipped gives the divisibility
by $n_i$ that was still left to prove.
\end{proof}
\par

%variables $\bfn = (n_1,\ldots,n_t)$
%for any collection of linear forms $\bff = (f_1,\ldots,f_s)$ in~$t$

%%%%%%%%%%%%%%%%%%%
\subsection{Quasimodularity of the number of pillowcase covers}
\label{sec:QMcovers}
%%%%%%%%%%%%%%%%%%%

Recall that we introduced in Section~\ref{sec:countingPillow} the generating
functions $N^0(\Hmu)$ and $N'(\Hmu)$ of the number of pillowcase covers
that are connected resp.\ without unramified components.
\par
\begin{Cor} \label{cor:noname}
For any ramification profile~$\Hmu$ the counting 
function $N^0(\Hmu)$ for connected pillowcase covers of profile~$\Hmu$ 
is a quasimodular form for the group~$\Gamma_0(2)$ of mixed weight less or
equal to ${\rm wt}(\Hmu) = |\Hmu| + \ell(\Hmu)$.
\end{Cor}
\par
\begin{proof} In~\eqref{eq:Covwbracket} we recalled that $N'(\Hmu)$
is the $w$-bracket of some element in~$\ol{\Lambda}$. By
Theorem~\ref{thm:bracketgraphsum} this series is thus a linear combination
of auxiliary brackets with entries~$p_\ell$ in the first arguments
and a product of $\olp_k$'s as last argument. The classical polynomiality
for triple Hurwitz numbers with $p_k$-arguments (summarized as
\cite[Theorem~4.1]{GMab}) and polynomiality in Theorem~\ref{thm:pbarpoly}
imply that both auxiliary functions $A'(\cdot)$ and $A_2'(\cdot)$ appearing in
the definition~\eqref{eq:defauxbrack} of auxiliary brackets are indeed
polynomials. That is, the auxiliary bracket is the sum over all subsets
$E^+(\Gamma)$ and all parity conditions $\pari \in \PC(\Gamma)$ of
graph sums of the form defined
by~\eqref{eq:SGammam} and~\eqref{eq:SGm}. By Theorem~\ref{thm:QMforS} such a
graph sum defines a quasimodular form of the weight as claimed. This gives
the result for $N'(\Hmu)$ and the claim for  $N^0(\Hmu)$ follows from
inclusion-exclusion, see e.g.\ \cite[Proposition~2.1]{GMab}.
\end{proof}

%%%%%%%%%%%%%%%%%%%%%%%%%%%%%%%%%%%%%%%%%%%%%%%%%%%%%%%%%%%%%%%%%%%%%%%%
\section{Application to Siegel-Veech constants} 
\label{sec:appSV}
%%%%%%%%%%%%%%%%%%%%%%%%%%%%%%%%%%%%%%%%%%%%%%%%%%%%%%%%%%%%%%%%%%%%%%%%

In this section we show that counting pillowcase covers with certain weight functions
also fall in the scope of the quasimodularity theorems and we prove
Theorem~\ref{intro:SV}. 
\par
Let $\lambda = (\lambda_1 \geq \lambda_2 \geq \cdots \geq \lambda_k)$ be 
a partition. For $p \in \ZZ$ we define the {\em $p$-th Siegel-Veech weight} 
of~$\lambda$ to be $%\be %\label{eq:def_pSV}
S_p(\lambda) = \sum_{j=1}^{k} \lambda_j^p\,.
$%\ee
With the conventions of Section~\ref{sec:countingPillow}, in particular the
definition of Hurwitz tuples in \eqref{eq:HT}, the core curves of the
horizontal cylinders have the monodromies
\bas \usi_0 &\=\ual_1\ual_4=\ual_2\ual_3(\uga_1\dots \uga_n)^{-1} \\
\usi_1 &\=\ual_1\ual_4\uga_1=\ual_2\ual_3(\uga_2\dots \uga_n)^{-1} \\
&\,\dots \\
\usi_n &\=\ual_1\ual_4\uga_1\dots \uga_n=\ual_2\ual_3\,.
\eas
Motivated by the relation to area-Siegel-Veech constants in
Proposition~\ref{prop:SVSV} below, we define the Siegel-Veech weighted 
Hurwitz numbers of a Hurwitz tuple~$h$ to be
\be \label{eq:SVofHT}
S_p(h) \= \sum_{i=0}^{n}S_p(\usi_i(h))\,.
\ee
Next, for $* \in \{', 0, \emptyset\}$ we package them into
the generating series 
\begin{equation} \label{eq:cpsimple}
c^*_{p}(d,\Hmu) \=   \sum_{j=1}^{|\Hur^*_d(\Hmu)|} S_p(\alpha^{(j)})\,,
\quad \text{and} \quad c_p^*(\Hmu) = \sum_{d \geq 0} c_{p}^*(d,\Hmu) q^d\,.
\end{equation}
These series admit the following graph sum decomposition. Let
$\widetilde{\NN}_\reg^{E(G)}$ be the special case of the height space
defined in~\eqref{eq:heightspace} with all the horizontal cylinders
on the base pillow of the same height, i.e.\ with $\ve_{4+i} = i/2(n-1)$
for $i=1,\ldots,n$.
\par
\begin{Prop} \label{prop:svgraphsum}
The generating series $c_{p}'(\Hmu)$ can be expressed in terms of
graph sums of triple Hurwitz numbers as
\bes %\label{eq:cpfmu}
c_{p}'(\Hmu) \=  \sum_{\Gamma} \frac1{|\Aut(\Gamma)|}c_p'(\Hmu,\Gamma)
\quad \text{where} \quad c_p'(\Hmu,\Gamma) \= \sum_{G \in \Gamma}c_p'(\Hmu,G)
\ees
and where
\bes
c_p'(\Hmu,G)=\!\!\!\sum_{{h\in \widetilde{\NN}_\reg^{E(G)}}, \atop w\in \Z_+^{E(G)}} \!\!
\Bigl( \sum_{e \in E(G)} \!\!\!h_e w_e^p \Bigr)
\prod_{e\in E(G)} w_e 
q^{h_e w_e} A_2'(\bfw_0,\nu) \!\!\prod_{v\in V(G)}A'(\bfw_v^-,\bfw_v^+,\mu_v)\,\, \delta(v)\,.
\ees
\end{Prop}
\par
\begin{proof}
The definition~\eqref{eq:SVofHT} together with the definition of the
Siegel-Veech weight is made such that a covering defined by a Hurwitz
tuples is counted with the weight given by the sum over all
horizontal cylinders~$C$ of $h(C)w(C)^p$. This results in the extra
factor in the formula for $c_p'(\Hmu,G)$ in comparison with the
formula for $N'(\Hmu,G)$ in~\eqref{eq:NpGamma2}. The summation over
$\widetilde{\NN}_\reg^{E(G)}$ is needed to ensure that each strip
of each cylinder is counted with the same weight. The whole formula is a direct
consequence of the correspondence theorem, i.e.\ of
Proposition~\ref{prop:corr}.
\end{proof}
\par
\begin{Cor} \label{Cor:SVmain}
For any ramification profile $\Hmu$ and any odd $p \geq -1$ the generating
series $c_p'(\Hmu)$ for counting pillowcase covers without
unramified components and with $p$-Siegel-Veech weight as well as the
generating 
series $c_p^0(\Hmu)$ for connected counting with $p$-Siegel-Veech weight
are quasimodular forms for the group~$\Gamma_0(2)$ of mixed
weight $ \leq {\rm wt}(\Hmu) +p+1$.
\end{Cor}
\par
Roughly, this follows from the polynomiality of $A'(\cdot)$ and $A_2'(\cdot)$
(i.e.\ from Theorem~\ref{thm:pbarpoly}) in a similar way as
Corollary~\eqref{cor:noname}. The extra factor $\sum_{e \in E(G)} h_e w_e^p$
raises the degree of the polynomial by $p+1$ and this results in the shifted
weight. We explain the procedure in detail in Section~\ref{sec:QMSV}

%%%%%%%%%%%%%%%%
\subsection{Relation to area Siegel-Veech constants} \label{sec:reltoSV}
%%%%%%%%%%%%%%%%

Siegel-Veech constants measure the growth rates of the number of saddle
connections or closed geodesics or equivalently embedded cylinders. Among
the various possibilities of weighting the count, the area weight is the
most important due to its connection to the sum of Lyapunov exponents (\cite{ekz}). 
In detail, 
\[c_{\textrm{area}}(X) \ =\lim\limits_{L\to\infty}\frac{N_{\textrm{area}}(T,L)}{\pi L^2},
\quad \text{where} \quad N_{\textrm{area}}(T,L) \= \sum_{Z\subset X\;\textrm{ cylinder}, \atop w(Z)\geq L}
\frac{\mbox{Area}(Z)}{\mbox{Area}(X)}.\]
is called the {\em (area) Siegel-Veech constant} of the flat surface~$X$. This constants 
are interesting both for generic flat surfaces of a given singularity type and
for pillowcase covers.
\par
\begin{Prop} \label{prop:SVSV}
  The area Siegel-Veech constant is related to Siegel-Veech weighted Hurwitz
  numbers by
  \[c_{\textrm{area}}(d,\Hmu) \= \frac{3}{\pi^2}\frac{c_{-1}^0(d,\Hmu)}{N_d^0(\Hmu)}\,.\]
\end{Prop}
\par
In particular, knowing the numerator and the denominator of the right hand
side to be quasimodular forms, and thus knowing the asymptotic behavior of
both $c_{-1}^0(d,\Hmu)$ and $N_d^0(\Hmu)$ as $d \to \infty$ allows to compute the
area Siegel-Veech constant of a generic surface with a given singularity type.
\par
\begin{proof}
The proof of \cite[Theorem 4]{ekz} or \cite[Theorem 3.1]{CMZ}
is easily adapted from torus covers to pillowcase covers.
\end{proof}

%%%%%%%%%%%%%%%%
\subsection{Quasimodularity of Siegel-Veech weighted graph sums} \label{sec:QMSV}
%%%%%%%%%%%%%%%%

The goal of this section is to prove Corollary~\ref{Cor:SVmain}. This will follow
from the following proposition. Recall from Section~\ref{sec:GraphSumisQM} the definition
of the distinguished edges $E^+(\Gamma)$ and the parity conditions $\pari \in \PC(\Gamma)$.
\par
\begin{Prop}  \label{prop:SVsumQM}
If $\mm=(m_1, \dots , m_{|E(\Gamma)|})$ is a tuple of even integers, then
for each $e_0 \in E(G)$ the graph sum ${S}^{SV}_{e_0}(\Gamma,E^+,\mm,\pari)
=\sum_{G\in(\Gamma,E^+)} {S}^{SV}_{e_0}(G,\mm,\pari)$, where
\bes
{S}^{SV}_{e_0}(G,\mm,\pari) \=\sum_{h\in\widehat\NN^{E(G)}, \atop w\in \NN_+^{E(G)^*}} 
\sum_{\bfw_0 \in \NN_{>0}^{E^0(\Gamma)} \atop  \bfw_0 \cong \pari\,{\rm mod}\, 2}
\frac{h_{e_0}}{w_{e_0}}\prod_{e\in E(G)}w_i^{m_i+1}q^{h_iw_i}\prod_{v\in V(G)^*}\delta(v)\,,
\ees
is a quasimodular form of mixed weight at most $k(\mm)=\sum_i(m_i+2)$.
\end{Prop}
\par
\begin{proof}
We may reduced to the reduced graph $\OG$ by computing the loop
contributions separately, compare Lemma~\ref{le:loopsplit} or rather \cite[Lemma~7.5]{GMab}.
If the $e_0$ is not a loop, then the loop contributions are quasimodular
by Lemma~\ref{le:loopOK}. If $e_0$ is a loop, we also need to take the extra
factor $h/w$ into account and note that $\sum_{w,h=1}^\infty hw^{m}q^{hw} = D_qS_{m}$ is
quasimodular (for $m \geq 2$ even) and similarly for the odd and even variants appearing
the proof of Lemma~\ref{le:loopOK}.
\par
To deal with $\OG$, we combine the construction of
Section~\ref{sec:ContourInt} and the Siegel-Veech weight in the proof
of \cite[Theorem~7.3]{GMab}. More precisely, we define if $e_0 \in E^*(\Gamma)$
\bas %\label{eq:defSVPG1}
P^{SV}_{\OG, E^+,\mm, \pari}(\zz)\=
\frac{D_q P^{(m_{e_0}-2)}(z_{v_1(e_0)} \pm z_{v_2(e_0)})}
{P^{(m_{e_0})}(z_{v_1(e_0)} \pm z_{v_2(e_0)})}
\cdot P_{\OG, E^+,\mm, \pari}(\zz)
\eas
if $m_{e_0}\geq 2$ and in the remaining case $m_{e_0} = 0$ we let
\bas %\label{eq:defSVPG2}
P^{SV}_{\OG, E^+,\mm, \pari}(\zz)\=
\frac{L(z_{v_1(e_0)} \pm z_{v_2(e_0)})}
{P(z_{v_1(e_0)} \pm z_{v_2(e_0)})}
\cdot P_{\OG, E^+,\mm, \pari}(\zz)
\eas
with~$L$ as in~\eqref{eq:LZ2}. In both cases the sign is chosen  according to
$e_0 \in E^\pm(\Gamma)$. 
This definition replaces the factor in~$P_{\OG, E^+,\mm, \pari}(\zz)$ corresponding
to the edge~$e_0$ in $P_{\OG, E^+,\mm, \pari}(\zz)$ is replaced by one with the extra
factor $h_{e_0}/w_{e_0}$. This follows from the power series expansion of~$P$
and~$L$ given in \cite[Equation~(41)]{GMab}, compare also the proof of Theorem~7.3
in loc.~cit. If $e_0 \in E^0(\Gamma)$, we define similarly 
\bas % \label{eq:defSVPG3}
P^{SV}_{\OG, E^+,\mm, \pari}(\zz) &\=
\frac{D_q P_{\pari_{e_o}}^{(m_{e_0}-2)}(2z_{v_2(e_0)})}
{P_{\pari_{e_o}}^{(m_{e_0})}(2z_{v_2(e_0)})}
\cdot P_{\OG, E^+,\mm, \pari}(\zz) \quad \text{or}\\
P^{SV}_{\OG, E^+,\mm, \pari}(\zz) &\=
\frac{L_{\pari_{e_0}}(2z_{v_2(e_0)})}
{P_{\pari_{e_0}}(2z_{v_2(e_0)})}
\cdot P_{\OG, E^+,\mm, \pari}(\zz) \\
\eas
according to $m_{e_0} \geq 2$ or $m_{e_0} =0$ respectively, where
$L_{\even}(z,\tau) = 2L(2z,2\tau)$ and $L_{\odd} = L - L_\even$. We this modified prefactor,
the same proof as in Proposition~\ref{prop:S1} shows that
\ba \label{eq:SVSasContInt}
{S}^{SV}(\OG,E^+,\mm, \pari)\= [\zeta_n^0, \dots, \zeta_1^0]\,
P^{SV}_{\OG, E^+,\mm, \pari}(\zz; \tau)\,.
\ea
Each of the factors in the definition of $P^{SV}_{\OG, E^+,\mm, \pari}(\zz)$ is a quasi-elliptic
quasimodular form by Proposition~\ref{prop:PZinQring}  and Proposition~\ref{prop:Qring} 
As in the case without Siegel-Veech weight, the claim follows from Theorem~\ref{thm:coeff0QM}
and the upgrade Lemma~\ref{le:EisG2} to get quasimodularity by the group $\Gamma_0(2)$.
\end{proof}
\par
\begin{proof}[Proof of Corollary~\ref{Cor:SVmain}]
We want to apply Proposition~\ref{prop:SVsumQM}. First, we pretend
for the moment that the summation in~\ref{prop:svgraphsum} is over
the normalized height space~$\widehat\NN^{E(G)}$ rather than
over~$\widetilde\NN_\reg^{E(G)}$ on prove quasimodularity of the
corresponding sum. Second,
to reduce the graph sum expression for Siegel-Veech constants in
Proposition~\ref{prop:svgraphsum} to those with polynomial entries, we have to
mimic the argument leading to Theorem~\ref{thm:bracketgraphsum}. Let
\bes 
[F_1,\ldots,F_n;F_0]^{p-SV} =  \sum_{\Gamma}\sum_{G \in \Gamma }\, [F_1,\ldots,F_n;F_0]^{p-SV}_G
\ees
and 
\bas
&\phantom{\=} [F_1,\ldots,F_n;F_0]_G^{p-SV} \\
&\=\!\!\!\sum_{h\in \widetilde{\NN}^{E(G)}, \atop w\in \Z_+^{E(G)}} \!\!
\Bigl( \sum_{e \in E(G)} \!\!\!h_e w_e^p \Bigr)
\prod_{e\in E(G)} w_e 
q^{h_e w_e} A_2'(\bfw_0,F_0) \!\!\prod_{v\in V(G)}A'(\bfw_v^-,\bfw_v^+,F_{\#_v})\,\, \delta(v)\,.
\eas
Then Proposition~\ref{prop:svgraphsum} can be generalized using
the correspondence theorem in the form Proposition~\ref{prop:corrvariant} to
\be \label{eq:SVcorrformula}
c_{p}'(\Hmu) \= [\prod_{i \not\in S} f_{\mu_i};
  \prod_{i \in S} f_{\mu_i}g_{\nu}]^{p-SV}
\ee
for any subset $S \subseteq \{1,\ldots,n\}$. To prove quasimodularity we
start with $S=\emptyset$ and decompose $g_\nu$ as linear combination of
$P_{I,J} = \prod_{j in J} p_{b_j } \, \prod_{i \in I} \olp_{a_i}$ for some~$a_i$ and~$b_j$.
For the summands where $J = \emptyset$ we write the first argument of the
bracket as linear combination of products of $p_k$ and use the polynomiality
for triple Hurwitz numbers with $p_k$-argument and Theorem~\ref{thm:pbarpoly}
to conclude thanks to Proposition~\ref{prop:SVsumQM}.
\par
To deal with the summands where $J = \emptyset$, we write $P_{I,J}$ a linear
combination with each time a product of $f_k$'s and one $g_\nu$. This is possible
by Theorem~\ref{thm:gnugens}. Since $J = \emptyset$, each term involves at least
one $f_k$. For each term we now apply~\eqref{eq:SVcorrformula} twice to move the $f_k$-product
to the first argument of the auxiliary bracket. By this procedure we have
reduced the weight to the~$g_\nu$ in the second argument and we conclude by induction.
\par
We have to justify the first simplification concerning height spaces.
Note that the shift to~$\widehat\NN^{E(G)}$ does not change the $q$-exponents
by Lemma~\ref{lem:heightspace}, so we may focus on the Fourier coefficients.
Since the expression for $c_p'(\Hmu)$ in Proposition~\ref{prop:svgraphsum}
involves a summation over all orientations we may combine the contributions
of~$G$ and the reverse orientation~$-G$, where the arrows of all
edges except for those emmanating from~$v_0$ have been inverted.
We thus obtain a pre-factor of
\bes
\sum_{h \in \NN_{\geq a_e}} \frac{h - \Delta(e)}w +
\sum_{h \in \NN_{\geq (1-a_e)}}\frac{h +  \Delta(e)}w
  \quad
\text{instead of}
\quad
\sum_{h \in \NN_{\geq a_e}} \frac{h}w +
\sum_{h \in \NN_{\geq (1-a_e)}}\frac{h}w
\ees
when using~$\widehat\NN^{E(G)}$. The difference is only the term $h=0$,
i.e.\ without an $h/w$-prefactor, and we recursively know these graph
sums to be quasimodular (in fact of smaller weight).
\par
Finally, the quasimodularity for $c_{p}^0(\Hmu)$ follows from the usual
inclusion-exclusion formulas, 
see \cite[Proposition~6.2]{CMZ} for the version with Siegel-Veech weight.
\end{proof}

%%%%%%%%%%%%%%%%
\subsection{Siegel-Veech weight and representation theory} \label{sec:SVandRep}
%%%%%%%%%%%%%%%%

The reader familiar with \cite{CMZ} will recall that counting function
for Hurwitz tuples, even with Siegel-Veech weight, can be expressed efficiently
using the representation theory of the symmetric group. More precisely,
for $\Hmu$ the profile of a torus cover
\be \label{eq:Burnside}
 N_d(\Hmu) \= \sum_{\lambda \in \Part(d)} \prod_{i=1}^n f_{\mu^{(i)}}(\lambda)\quad 
\text{and}\quad
 c_{p} (d, \Hmu) \=   
\sum_{\lambda \in \Part(d)} \prod_{i=1}^nf_{\mu_i}(\lambda)   T_p(\lambda)\,,
%\frac 1{d!}\sum_{\tau\in \Part(d)} z_{\tau}  S_p(\tau) \chi^\lambda(\tau)^2. 
\ee
where $T_p(\lambda) \= \sum_{\xi \in Y_\lambda} h(\xi)^{p-1}$ and where
$h(\xi)$ is the hook-length of the cell~$\xi$ of the Young diagram $Y_\lambda$. 
\par
It would be very useful to have a similar formula for the Siegel-Veech weighted
counting of pillowcase covers. We are only aware of the following much more
complicated formula.
\par
\begin{Prop} \label{prop:SVviaRep}
The number of all covers of degree~$d$ with profile~$\Hmu$
counted with $p$-Siegel-Veech weight is 
\bas
c_p(d, \Hmu) &\=\frac{1}{l(\mu)+1}\sum_{k=0}^{l(\mu)}\sum_C {S_p(C)}\sum_{\lambda, \lambda'}\sqrt{w(\lambda)w(\lambda')}|C|\chi^{\lambda}(C)\chi^{\lambda'}(C) \cdot \\
&\qquad \qquad \qquad\qquad \qquad \qquad \cdot
g_{\nu}(\lambda')\prod_{i=1}^k f_{\mu_i}(\lambda')\prod_{i=k+1}^{l(\mu)}f_{\mu_i}(\lambda) \,.
\eas
\end{Prop}
\par
In particular we are not aware of an operator on Fock space whose $q$-trace computes
the generating series with Siegel-Veech weight. Note that the $\mathfrak{W}$-operator
of \cite[Theorem~4]{eopillow} has the property $\langle v_\lambda|\mathfrak{W}
|v_\lambda \rangle = w(\lambda)$,
but it is not true that  $\langle v_\lambda|\mathfrak{W}
|v_\nu \rangle = \sqrt{w(\lambda)w(\nu)}$ for $\lambda \neq \nu$. Finding a vertex
operator with this property would be a way to use Proposition~\ref{prop:SVviaRep}
to express Siegel-Veech weighted generating series as $q$-traces.
\par
\begin{proof} The monodromy of the core curves of the cylinders 
of a Hurwitz tuple $h \in \Hur_d(\Hmu)$ is given by
\bas \usi_0 &\=\ual_1\ual_4=\ual_2\ual_3(\uga_1\dots \uga_n)^{-1} \\
\usi_1&\=\ual_1\ual_4\uga_1=\ual_2\ual_3(\uga_2\dots \uga_n)^{-1} \\
\dots \\
\usi_n&\=\ual_1\ual_4\uga_1\dots \uga_n=\ual_2\ual_3\,.
\eas
To count Hurwitz tuples with Siegel-Veech weight, say for the $k$-th cylinder,
we split the defining equation as 
\be \label{eq:csplit}
\ual_1\ual_4\uga_1, \dots, \uga_k \= c
\=\ual_2\ual_3\uga_n^{-1}\dots \uga_{k+1}^{-1} \ee
and count the solutions of each side separately. That is, we denote by $C_1,
C_2, C_3, C_4$, $C, C^{\mu}_1, \dots, C^{\mu}_n$, respectively,  the conjugacy
classes  of permutations of type
\[(\nu, 2^{|\nu|/2-d}), (2^d), (2^d), (2^d), (k^{c_k}), (\mu_1, 1^{d-\mu_1}),
\dots, (\mu_n, 1^{d-\mu_n}).\]
We denote by $c_p(S_{2d};C_1, C_4,C_1^{\mu}, \dots, C_k^{\mu}, C )$ the number of
solutions of 
\be \label{eq:partoffundeq}
\ual_1\ual_4\uga_1, \dots, \uga_k\= c
\ee 
with $\ual_i$ of conjugacy class $C_i$ and $\uga_i$ of conjugacy class
$C^{\nu}_i$), $c$ of conjugacy class $C$, counted with weight $S_p(c)=S_p(C)$.
We claim that
\bas c_p(S_{2d};C_1, C_4,C_1^{\mu}, \dots, C_k^{\mu}, C )=\frac{|C_1||C_4||C_1^{\mu}|\dots|C_k^{\mu}||C|}{|G|}S_p(C)\cdot \\ \cdot\sum_{\chi} \frac{\chi(C_1)\chi(C_4)\chi(C_1^{\mu})\dots\chi(C_k^{\mu})\chi(C)}{\chi(1)^{k+1}}\,.
\eas
To see this, we revisit the proof of the orthogonality relations, see
\cite[Theorem 7.2.1]{serre08}. We introduce the class function
\[\phi(x)=S_p(C){1}_{\{x\in C\}}=\sum_{\chi}c_{\chi}\chi\] with
\[c_{\chi}=\int_G S_p(C){1}_{\{x\in C\}}\overline\chi(x)dx=\frac{|C|}{|G|}S_p(C)\chi(C).\] 
We have
\bas
I(\phi)\=\int_{G^{k+2}} \phi(t_1\ual_1t_1^{-1}t_4\ual_4t_4^{-1}s_1\uga_1s_1^{-1}\dots s_k\uga_ks_k^{-1}y)dt_1dt_4ds_1\dots ds_k\\=\sum_{\chi}c_{\chi}\frac{\chi(C_1)\chi(C_4)\chi(C_1^{\mu})\dots
  \chi(C_k^{\mu})\chi(y)}{\chi(1)^{k+2}}
\eas
and the left hand side is $I(\phi)=N_p(S_{2d};C_1, C_4,C_1^{\mu}, \dots, C_k^{\mu}, C )/|G|^{k+2}$. Taking $y=1$ we get the formula of the claim.
\par
In the second step we count the solution of the right equality
of~\eqref{eq:csplit} with weight $c_p(S_{2d};C_1, C_4,C_1^{\mu}, \dots, C_k^{\mu},
C)$. The class function is now
\[\phi(x) \= \frac{1}{|C|} c_p(S_{2d};C_1, C_4,C_1^{\mu}, \dots, C_k^{\mu}, x )
1_{x\in C}\,,\] and its coefficients are
\[c_{\chi} \=\frac{1}{|G|}N_p(S_{2d};C_1, C_4,C_1^{\mu}, \dots, C_k^{\mu}, C )
\overline{\chi}(C)\,.\]
With the argument as above we conclude that the number of solutions counted
with weight is
\[|C_2||C_3||C_n^{\mu}|\dots |C_{k+1}^{\mu}|\sum_{\chi}c_{\chi}
\frac{\chi(C_2)\chi(C_3)\chi(C_n^{\mu})\dots \chi(C_{k+1}^{\mu})}
{\chi(1)^{n-k+1}}\,.\]
We sum now on all conjugacy classes $C$ and use the definition of $f_\mu$
and $g_\nu$ to obtain the result.
\end{proof}

%%%%%%%%%%%%%%%%%%%%%%%%%%%%%%%%%%%%%%%%%%%%%%%%%%%%%%%%%%%%%%%%%%%%%%%%%%%%%
\section{Example: $\cQ(2,1,-1^3)$} 
\label{sec:Example}
%%%%%%%%%%%%%%%%%%%%%%%%%%%%%%%%%%%%%%%%%%%%%%%%%%%%%%%%%%%%%%%%%%%%%%%%%%%%%

In this section, we treat the example of the stratum $\cQ(2,1,-1^3)$ from A
to Z to illustrate all sections of the paper. This stratum is the lowest
dimensional example exhibiting all the relevant aspects.
\par
Integer points in the stratum $\cQ(2,1,-1^3)$ correspond in our setting to
covers of the pillow ramified over five points: the four corners
$P_1, \dots P_4$ and an additional point $P_5$ (see Figure~\ref{cap:PP1gens}),
with the ramification profile
$\Pi=(\mu^{(1)},\mu^{(2)}, \mu^{(3)}, \mu^{(4)}, \mu^{(5)})$ where
$\mu^{(1)}=(3,1,1,1,2^{d-3})$ over $P_1$, $\mu^{(2)}=\mu^{(3)}=\mu^{(4)}=(2^d)$
over $P_2, P_3, P_4$ and $\mu^{(5)}=(2, 1^{2d-2})$ over $P_5$. Here $2d$ is the
degree of the cover. In this particular case, we cannot have covers with at
least two ramified connected components, so $N'(\Pi)=N^0(\Pi)$.

%%%%%%%%%%%%%%%%%%%%%%%%%%%%%%%%
\subsection{Counting covers}
%%%%%%%%%%%%%%%%%%%%%%%%%%%%%%%

By \cite{eopillow} and Theorem~\ref{intro:count} the generating series
$N^0(\Pi)$ is a quasimodular form for $\Gamma_0(2)$ of mixed weight less
or equal to $6$. Computing the first coefficients of the series we get
that\footnote{We choose in this paper to consider series in $q^{d}$ rather
than as series in $q^{2d}$ as in \cite{eopillow}. The integer $d$ can be
seen as the area of the cover while $2d$ is the degree of the cover.}
\bas
N^{0}(\Pi)&\= 360G_2(q^2)^3 -360G_2(q)G_2(q^2)^2  + 72G_2(q)^2G_2(q^2)  - 30G_4(q^2)G_2(q^2)\\&  - \frac{5}{4}G_4(q^2)+ 3G_2(q)^2+ 15G_2(q^2)^2- 15G_2(q)G_2(q^2).
\eas
Our goal is to retrieve this result by considering all graph contributions.
Standard Hurwitz theory (see~\eqref{eq:A2withF} and~\cite{eopillow}) gives that 
\[N^{0}(\Pi)=\frac{1}{N(\Pi_\emptyset)}\sum_{\lambda}g_{3,1,1,1}(\lambda)f_2(\lambda)\textrm{w}^2(\lambda)q^{|\lambda|/2}\]
where $g_{3,1,1,1}(\lambda)= f_{3,1,1,1,2,\dots,2}(\lambda) /{f_{2, \dots, 2}(\lambda)}$
and $w(\lambda)$ as in~\eqref{eq:defw}.
%where $g_{3,1,1,1}(\lambda)=\cfrac{f_{3,1,1,1,2,\dots,2}(\lambda)}{f_{2, \dots, 2}(\lambda)}$ and $\textrm{w}(\lambda)=\left(\cfrac{\dim(\lambda)}{|\lambda|!}\right)f_{2, \dots, 2}(\lambda)^2$.
\par
Following the strategy in Theorem~\ref{thm:bracketgraphsum} we need to
express the graph sums with local contributions $g_{3,1,1,1}$ and~$f_2$
by graphs sums using the functions $p_k$ and $\olp_k$ for which we have nice
polynomiality results. We compute that 
\bes g_{(3,1,1,1)} \=-\frac{1}{4}\olp_1p_1+\frac{1}{108}\olp_1^3-\frac{1}{36}\olp_2\olp_1+\frac{3}{8}\olp_1+\frac{2}{27}\olp_3, \quad 
f_2 \=\frac{1}{2}p_2
\ees
The definition
\bas \ol{g}_{(3,1,1,1)}&:=\frac{1}{108}\olp_1^3-\frac{1}{36}\olp_2\olp_1+\frac{3}{8}\olp_1+\frac{2}{27}\olp_3\\
{g}^{deg}_{(3,1,1,1)}&:=-\frac{1}{4}\olp_1
\eas
and the resulting decomposition
\[g_{(3,1,1,1)}=\ol{g}_{(3,1,1,1)}+g^{deg}_{(3,1,1,1)}p_1.\]
isolates the products of $\olp_i$'s so we can use the polynomiality results
of Section~\ref{sec:QPoly}. Writing
$\ol{g}_{3111}f_2 = \frac{1}{2}f_2f_1g_{11} + f_2g_{3111}$
and $g_{3111}^{deg}p_1f_2 = -\frac{1}{2}g_{11}f_2f_2+\frac{1}{48}g_{11}f_2$ we can make geometric sense
of the formal decomposition, as counting degenerate covers, e.g.\
in the stratum $\cQ(2,0,-1^2)$. We determine the simple Hurwitz numbers with
2-stabilization $A'_2$  for products of $\olp_i$'s, proven to be
quasi-polynomials in Theorem~\ref{thm:pbarpoly} by computing the first few
terms. As a result 
\bas A_2'\big((w),\ol{g}_{(3,1,1,1)}\big)&=\frac{1}{24}w^2+\frac{1}{3} \\
A_2'\big((w_1,w_2,w_3),\ol{g}_{(3,1,1,1)}\big)&=\begin{cases}\frac{1}{2} \mbox{ if two }w_i\mbox{ are odd}\\
\frac{3}{2} \mbox{ if all }w_i\mbox{ are even}\end{cases}\\
A'_2\big((w),g^{deg}_{(3,1,1,1)}\big)&=-\frac{1}{4}
\eas 
We recall that these polynomials are only defined for $\sum w_i$ even.
Similarly, the double Hurwitz numbers $A'$ for products of $p_i$'s
are polynomials, in fact
\bes
A'\big((w),(w),p_1\big)\=1\,, \quad A'\big((w_1),(w_2,w_3),f_2\big) \= 1\,,
\ees
We have thus collected all local polynomials.
\par
\begin{figure}
\begin{center}
\begin{tikzpicture}[line cap=round,line join=round,>=triangle 45,x=1.0cm,y=1.0cm]
\clip(-0.86,0.12) rectangle (7.8,3.98);
\draw [line width=1pt] (0.,1.)-- (0.,2.);
\draw [line width=1pt] (0.,2.5) circle (0.5cm);
\draw (0.1,2.1) node[anchor=north west] {$f_2$};
\draw (0.24,1.1) node[anchor=north west] {$\overline{g}_{(3,1,1,1)}$};
\draw [line width=1pt] (3.,3.)-- (3.,1.);
\draw [shift={(3.5333333333333328,2.)},line width=1pt]  plot[domain=2.0607536530486246:4.222431654130961,variable=\t]({1.*1.133333333333333*cos(\t r)+0.*1.133333333333333*sin(\t r)},{0.*1.133333333333333*cos(\t r)+1.*1.133333333333333*sin(\t r)});
\draw [shift={(2.4666666666666672,2.)},line width=1pt]  plot[domain=-1.080839000541169:1.0808390005411688,variable=\t]({1.*1.133333333333333*cos(\t r)+0.*1.133333333333333*sin(\t r)},{0.*1.133333333333333*cos(\t r)+1.*1.133333333333333*sin(\t r)});
\draw (3.28,3.3) node[anchor=north west] {$f_2$};
\draw (3.22,1.1) node[anchor=north west] {$\overline{g}_{(3,1,1,1)}$};
\draw [line width=1pt] (6.,1.5)-- (6.,2.5);
\draw [line width=1pt] (6.,3.) circle (0.5cm);
\draw [line width=1pt] (6.,1.) circle (0.5cm);
\draw (6.16,2.7) node[anchor=north west] {$f_2$};
\draw (6.3,1.7) node[anchor=north west] {$\overline{g}_{(3,1,1,1)}$};
\draw (-0.6,3.9) node[anchor=north west] {$A$};
\draw (2.0,3.9) node[anchor=north west] {$B$};
\draw (4.9,3.9) node[anchor=north west] {$C$};
\begin{scriptsize}
\draw [color=black] (0.,1.)-- ++(-2.0pt,-2.0pt) -- ++(4.0pt,4.0pt) ++(-4.0pt,0) -- ++(4.0pt,-4.0pt);
\draw [fill=black] (0.,2.) circle (2.0pt);
\draw [fill=black] (3.,3.) circle (2.0pt);
\draw [color=black] (3.,1.)-- ++(-2.0pt,-2.0pt) -- ++(4.0pt,4.0pt) ++(-4.0pt,0) -- ++(4.0pt,-4.0pt);
\draw [color=black] (6.,1.5)-- ++(-2.0pt,-2.0pt) -- ++(4.0pt,4.0pt) ++(-4.0pt,0) -- ++(4.0pt,-4.0pt);
\draw [fill=black] (6.,2.5) circle (2.0pt);
\end{scriptsize}
\end{tikzpicture}
\end{center}
\caption{Admissible graphs for $\ol{g}_{(3,1,1,1)}\cdot f_2$}\label{fig:adm1}
\end{figure}
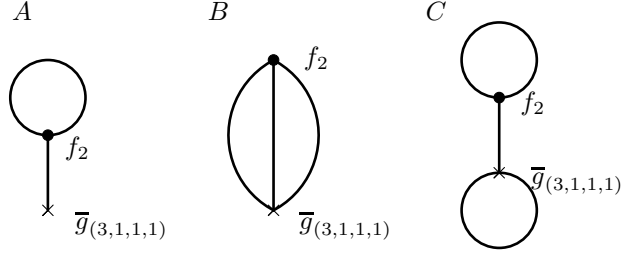

Now we glue the local surfaces together, encoding the gluings by the various possible
global graphs. The contribution of $\ol{g}_{(3,1,1,1)}\cdot f_2$ is encoded
in graphs with two vertices,  one special vertex $v_0$ of valency one or
three, and another trivalent vertex $v_1$ since all other valencies
result in a zero local polynomial and thus in a zero contribution.
By convention, we represent $v_0$ as the bottom of the graph, marked
with a cross. Disregarding orientations, we obtain three admissible
graphs, shown in Figure~\ref{fig:adm1}.
\par
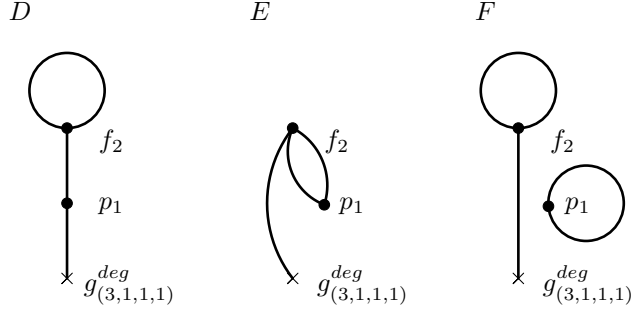
\begin{figure}
\begin{center}
\begin{tikzpicture}[line cap=round,line join=round,>=triangle 45,x=1.0cm,y=1.0cm]
\clip(-4.3,2.18) rectangle (5.08,6.94);
\draw [line width=1.pt] (-3.,3.)-- (-3.,5.);
\draw [line width=1.pt] (3.,3.)-- (3.,5.);
\draw [line width=1.pt] (-3.,5.5) circle (0.5cm);
\draw [line width=1.pt] (3.,5.5) circle (0.5cm);
\draw [line width=1.pt] (3.9,4.) circle (0.5015974481593781cm);
\draw [shift={(1.298235294117647,4.)},line width=1.pt]  plot[domain=2.4852404429949058:3.7979448641846805,variable=\t]({1.*1.638723551698923*cos(\t r)+0.*1.638723551698923*sin(\t r)},{0.*1.638723551698923*cos(\t r)+1.*1.638723551698923*sin(\t r)});
\draw [shift={(0.6977570093457943,4.690841121495327)},line width=1.pt]  plot[domain=2.724512232197865:4.339887162377095,variable=\t]({1.*0.7631802252741182*cos(\t r)+0.*0.7631802252741182*sin(\t r)},{0.*0.7631802252741182*cos(\t r)+1.*0.7631802252741182*sin(\t r)});
\draw [shift={(-0.3971428571428568,4.24)},line width=1.pt]  plot[domain=-0.3080527810237754:1.0892668684191507,variable=\t]({1.*0.8575094454171293*cos(\t r)+0.*0.8575094454171293*sin(\t r)},{0.*0.8575094454171293*cos(\t r)+1.*0.8575094454171293*sin(\t r)});
\draw (-3.9,6.8) node[anchor=north west] {$D$};
\draw (-0.7,6.8) node[anchor=north west] {$E$};
\draw (2.3,6.8) node[anchor=north west] {$F$};
\draw (-2.9,3.3) node[anchor=north west] {$g^{deg}_{(3,1,1,1)}$};
\draw (0.2,3.3) node[anchor=north west] {$g^{deg}_{(3,1,1,1)}$};
\draw (3.2,3.3) node[anchor=north west] {$g^{deg}_{(3,1,1,1)}$};
\draw (-2.7,4.2) node[anchor=north west] {$p_1$};
\draw (0.5,4.2) node[anchor=north west] {$p_1$};
\draw (3.5,4.2) node[anchor=north west] {$p_1$};
\draw (-2.7,5.1) node[anchor=north west] {$f_2$};
\draw (0.3,5.1) node[anchor=north west] {$f_2$};
\draw (3.3,5.1) node[anchor=north west] {$f_2$};
\begin{scriptsize}
\draw [color=black] (-3.,3.)-- ++(-2.0pt,-2.0pt) -- ++(4.0pt,4.0pt) ++(-4.0pt,0) -- ++(4.0pt,-4.0pt);
\draw [fill=black] (-3.,4.) circle (2.0pt);
\draw [fill=black] (-3.,5.) circle (2.0pt);
\draw [color=black] (0.,3.)-- ++(-2.0pt,-2.0pt) -- ++(4.0pt,4.0pt) ++(-4.0pt,0) -- ++(4.0pt,-4.0pt);
\draw [fill=black] (0.42,3.98) circle (2.0pt);
\draw [fill=black] (0.,5.) circle (2.0pt);
\draw [color=black] (3.,3.)-- ++(-2.0pt,-2.0pt) -- ++(4.0pt,4.0pt) ++(-4.0pt,0) -- ++(4.0pt,-4.0pt);
\draw [fill=black] (3.4,3.96) circle (2.0pt);
\draw [fill=black] (3.,5.) circle (2.0pt);
\end{scriptsize}
\end{tikzpicture}
\end{center}
\caption{Admissible graphs for $\ol{g}_{(3,1,1,1)}\cdot f_2$} \label{fig:adm2}
\end{figure}

Similarly, the contribution of $g^{deg}_{(3,1,1,1)}\cdot p_1\cdot f_2$ is encoded
in graphs with three vertices, one special vertex~$v_0$ of valency~$1$, a
vertex~$v_1$ of valency~$2$ and a vertex~$v_2$ of valency~$3$. The graphs
are listed in Figure~\ref{fig:adm2}.
\par
Next, we are supposed to sum over all possible orientations of these
graphs. We sort the orientations by the subset of coherently
oriented edges, i.e.\ by those distinguished with a~$+$ in the notation
of Section~\ref{sec:GraphSumisQM}. Note that certain decorations give trivial
contributions, as the graph $A$ with the loop decorated with $+$. In
fact, considering all possible orientations of half-edges compatible with
this decoration, we see that we get incompatible width conditions for the
vertex~$v_1$. This implies that integrating the corresponding propagator we
will get terms like $\delta(w_2=w_1+w_2)$ that are always trivial. 
\par
The next step is to associate to each decorated graph its propagator, and then
to integrate this propagator (get its $\zeta^0$-coefficient) to obtain the
contributions of each individual graph. In the following table we consider
only decorations with non trivial contributions.  For each vertex, we indicate
the corresponding integration variable $z_i$ (on the left of the vertex) and
the corresponding local polynomial (on the right). For each decorated graph
we give the associated propagator and the contribution to the volume.
The contributions were computed using three independent methods: the reduction algorithm used in the proof of Theorem~\ref{thm:coeff0QM}, the computation of the first terms of the q-series using graph sums \eqref{eq:SGm}, and, independently, using extraction of $[\zeta^0]$ coefficients, combined with a numerical test
of quasimodularity (test of linear dependency with the basis of quasimodular forms). 
\par
\begin{figure}
\label{table:example}
\begin{tabular}{|c|c|p{3.3cm}|p{4.1cm}|}
\hline
Graph $\Gamma$ &$\cfrac{1}{|\Aut(\Gamma)|}$& Propagator & Contribution\\
\hline
 \raisebox{-0.5\totalheight}{
\begin{tikzpicture}[line cap=round,line join=round,>=triangle 45,x=1.0cm,y=1.0cm]
\clip(-0.86,0.12) rectangle (2,3.98);
\draw [line width=1pt] (0.,1.)-- (0.,2.);
\draw [line width=1pt] (0.,2.5) circle (0.5cm);
\draw (0.5,2.6) node[anchor=north west] {$-$};
\draw (0.1,2.1) node[anchor=north west] {$1$};
\draw (-0.5,2.05) node[anchor=north west] {$z$};
\draw (0.1,1.1) node[anchor=north west] {$\frac{w_1^2}{24}+\frac{1}{3}$};
\draw (-0.6,3.9) node[anchor=north west] {$A$};
\begin{scriptsize}
\draw [color=black] (0.,1.)-- ++(-2.0pt,-2.0pt) -- ++(4.0pt,4.0pt) ++(-4.0pt,0) -- ++(4.0pt,-4.0pt);
\draw [fill=black] (0.,2.) circle (2.0pt);
\end{scriptsize}
\end{tikzpicture}}

&$\cfrac{1}{2}$& $\left(\frac{1}{24}P^{(2)}(z)+\frac{1}{3}P(z)\right)\cdot \newline \cdot  P(2z)$ 
& $\frac{7}{30}G_6-\frac{8}{3}G_4G_2+\frac{5}{9}G_4-\frac{4}{3}G_2^2$ \\
\hline
 \raisebox{-0.5\totalheight}{
\begin{tikzpicture}[line cap=round,line join=round,>=triangle 45,x=1.0cm,y=1.0cm]
\clip(2.,0.12) rectangle (4.5,3.98);
\draw [line width=1pt] (3.,3.)-- (3.,1.);
\draw [shift={(3.5333333333333328,2.)},line width=1pt]  plot[domain=2.0607536530486246:4.222431654130961,variable=\t]({1.*1.133333333333333*cos(\t r)+0.*1.133333333333333*sin(\t r)},{0.*1.133333333333333*cos(\t r)+1.*1.133333333333333*sin(\t r)});
\draw [shift={(2.4666666666666672,2.)},line width=1pt]  plot[domain=-1.080839000541169:1.0808390005411688,variable=\t]({1.*1.133333333333333*cos(\t r)+0.*1.133333333333333*sin(\t r)},{0.*1.133333333333333*cos(\t r)+1.*1.133333333333333*sin(\t r)});
\draw (3.28,3.3) node[anchor=north west] {$1$};
\draw (2.4,3.25) node[anchor=north west] {$z$};
\draw (3.22,1.1) node[anchor=north west] {$\frac{1}{2}$};
\draw (2.0,3.9) node[anchor=north west] {$B_1$};
\draw (1.85,1.7) node[anchor=north west] {$\odd$};
\draw (2.9,2.2) node[anchor=north west] {$\odd$};
\draw (3.6,1.9) node[anchor=north west] {$\even$};
\begin{scriptsize}
\draw [fill=black] (3.,3.) circle (2.0pt);
\draw [color=black] (3.,1.)-- ++(-2.0pt,-2.0pt) -- ++(4.0pt,4.0pt) ++(-4.0pt,0) -- ++(4.0pt,-4.0pt);
\end{scriptsize}
\end{tikzpicture}}
&$\cfrac{1}{2}$& $\frac{1}{2}P_{\odd}(z)^2P_{\even}(z)$ &$-\frac{1}{3}(16G_{22}G_{42}+8G_2G_{42})+\frac{1}{5}(-192G_{22}^3+608G_2G_{22}^2-384G_2^2G_{22}+64G_2^3) $\\
\hline
 \raisebox{-0.5\totalheight}{
\begin{tikzpicture}[line cap=round,line join=round,>=triangle 45,x=1.0cm,y=1.0cm]
\clip(2.,0.12) rectangle (4.5,3.98);
\draw [line width=1pt] (3.,3.)-- (3.,1.);
\draw [shift={(3.5333333333333328,2.)},line width=1pt]  plot[domain=2.0607536530486246:4.222431654130961,variable=\t]({1.*1.133333333333333*cos(\t r)+0.*1.133333333333333*sin(\t r)},{0.*1.133333333333333*cos(\t r)+1.*1.133333333333333*sin(\t r)});
\draw [shift={(2.4666666666666672,2.)},line width=1pt]  plot[domain=-1.080839000541169:1.0808390005411688,variable=\t]({1.*1.133333333333333*cos(\t r)+0.*1.133333333333333*sin(\t r)},{0.*1.133333333333333*cos(\t r)+1.*1.133333333333333*sin(\t r)});
\draw (3.28,3.3) node[anchor=north west] {$1$};
\draw (2.4,3.25) node[anchor=north west] {$z$};
\draw (3.22,1.1) node[anchor=north west] {$\frac{3}{2}$};
\draw (2.0,3.9) node[anchor=north west] {$B_2$};
\draw (2,1.8) node[anchor=north west] {$\even$};
\draw (2.7,2.2) node[anchor=north west] {$\even$};
\draw (3.5,1.8) node[anchor=north west] {$\even$};
\begin{scriptsize}
\draw [fill=black] (3.,3.) circle (2.0pt);
\draw [color=black] (3.,1.)-- ++(-2.0pt,-2.0pt) -- ++(4.0pt,4.0pt) ++(-4.0pt,0) -- ++(4.0pt,-4.0pt);
\end{scriptsize}
\end{tikzpicture}}
&$\cfrac{1}{6}$& $\frac{3}{2}P_{\even}(z)^3$ & $48G_{22}G_{42}-16G_2G_{42}+\frac{1}{5}(-832G_{22}^3+768G_2G_{22}^2-384G_2^2G_{22}+64G_2^3)$\\
\hline
 \raisebox{-0.5\totalheight}{
\begin{tikzpicture}[line cap=round,line join=round,>=triangle 45,x=1.0cm,y=1.0cm]
\clip(4.5,0.12) rectangle (7.3,3.98);
\draw [line width=1pt] (6.,1.5)-- (6.,2.5);
\draw [line width=1pt] (6.,3.) circle (0.5cm);
\draw [line width=1pt] (6.,1.) circle (0.5cm);
\draw (6.16,2.7) node[anchor=north west] {$1$};
\draw (5.45,2.65) node[anchor=north west] {$z$};
\draw (6.4,1.75) node[anchor=north west] {$\frac{1}{2}$};
\draw (4.9,3.9) node[anchor=north west] {$C_1$};
\draw (6.5,3.5) node[anchor=north west] {$-$};
\draw (6.5,1.2) node[anchor=north west] {$\odd$};
\draw (4.7,1.2) node[anchor=north west] {$\odd$};
\draw (6,2.1) node[anchor=north west] {$\even$};
\begin{scriptsize}
\draw [color=black] (6.,1.5)-- ++(-2.0pt,-2.0pt) -- ++(4.0pt,4.0pt) ++(-4.0pt,0) -- ++(4.0pt,-4.0pt);
\draw [fill=black] (6.,2.5) circle (2.0pt);
\end{scriptsize}
\end{tikzpicture}}
&$\cfrac{1}{4}$& $\frac{1}{2}S_{0,\odd}P_{\even}(z)P(2z)$ & $\left(G_2-2G_{22}-\frac{1}{24}\right)\cdot\newline\left(\frac{5}{6}G_4-2G_2^2\right)$\\
\hline
 \raisebox{-0.5\totalheight}{
\begin{tikzpicture}[line cap=round,line join=round,>=triangle 45,x=1.0cm,y=1.0cm]
\clip(4.5,0.12) rectangle (7.3,3.98);
\draw [line width=1pt] (6.,1.5)-- (6.,2.5);
\draw [line width=1pt] (6.,3.) circle (0.5cm);
\draw [line width=1pt] (6.,1.) circle (0.5cm);
\draw (6.16,2.7) node[anchor=north west] {$1$};
\draw (5.45,2.65) node[anchor=north west] {$z$};
\draw (6.3,1.75) node[anchor=north west] {$\frac{3}{2}$};
\draw (4.9,3.9) node[anchor=north west] {$C_2$};
\draw (6.5,3.5) node[anchor=north west] {$-$};
\draw (6.45,1.2) node[anchor=north west] {$\even$};
\draw (4.6,1.2) node[anchor=north west] {$\even$};
\draw (6,2.1) node[anchor=north west] {$\even$};
\begin{scriptsize}
\draw [color=black] (6.,1.5)-- ++(-2.0pt,-2.0pt) -- ++(4.0pt,4.0pt) ++(-4.0pt,0) -- ++(4.0pt,-4.0pt);
\draw [fill=black] (6.,2.5) circle (2.0pt);
\end{scriptsize}
\end{tikzpicture}}
&$\cfrac{1}{4}$& $\frac{3}{2}S_{0,\even}P_{\even}(z)P(2z)$ &$3\left(2G_{22}+\frac{1}{12}\right)\left(\frac{5}{6}G_4-2G_2^2\right)$\\
\hline
\end{tabular}
\caption{Graphs for $\cQ(2,1, -1^3)$: Contribution of $\ol{g}_{(3,1,1,1)}\cdot f_2$}
\end{figure}

\begin{figure}
\label{table:example2}
\begin{tabular}{|c|c|p{3.3cm}|p{4.1cm}|}
\hline
Graph $\Gamma$ &$\cfrac{1}{|\Aut(\Gamma)|}$& Propagator & Contribution\\
\hline
 \raisebox{-0.5\totalheight}{\begin{tikzpicture}[line cap=round,line join=round,>=triangle 45,x=1.0cm,y=1.0cm]
\clip(-4.3,2.18) rectangle (-1.5,6.94);
\draw [line width=1.pt] (-3.,3.)-- (-3.,5.);
\draw [line width=1.pt] (-3.,5.5) circle (0.5cm);
\draw (-3.9,6.8) node[anchor=north west] {$D_a$};
\draw (-2.9,3.3) node[anchor=north west] {$-\frac{1}{4}$};
\draw (-2.7,4.2) node[anchor=north west] {$1$};
\draw (-2.7,5.1) node[anchor=north west] {$1$};
\draw (-3.1,4.7) node[anchor=north west] {$+$};
\draw (-2.6,5.8) node[anchor=north west] {$-$};
\draw (-3.7,5.1) node[anchor=north west] {$z_2$};
\draw (-3.7,4.2) node[anchor=north west] {$z_1$};
\begin{scriptsize}
\draw [color=black] (-3.,3.)-- ++(-2.0pt,-2.0pt) -- ++(4.0pt,4.0pt) ++(-4.0pt,0) -- ++(4.0pt,-4.0pt);
\draw [fill=black] (-3.,4.) circle (2.0pt);
\draw [fill=black] (-3.,5.) circle (2.0pt);
\end{scriptsize}
\end{tikzpicture}}
&$\cfrac{1}{2}$& $-\frac{1}{4}P(z_1)P(z_1-z_2) \cdot \newline\cdot  P(2z_2)$&   $A-3B$ \\
\hline
 \raisebox{-0.5\totalheight}{\begin{tikzpicture}[line cap=round,line join=round,>=triangle 45,x=1.0cm,y=1.0cm]
\clip(-4.3,2.18) rectangle (-1.5,6.94);
\draw [line width=1.pt] (-3.,3.)-- (-3.,5.);
\draw [line width=1.pt] (-3.,5.5) circle (0.5cm);
\draw (-3.9,6.8) node[anchor=north west] {$D_b$};
\draw (-2.9,3.3) node[anchor=north west] {$-\frac{1}{4}$};
\draw (-2.7,4.2) node[anchor=north west] {$1$};
\draw (-2.7,5.1) node[anchor=north west] {$1$};
\draw (-3.1,4.7) node[anchor=north west] {$-$};
\draw (-2.6,5.8) node[anchor=north west] {$-$};
\draw (-3.7,5.1) node[anchor=north west] {$z_2$};
\draw (-3.7,4.2) node[anchor=north west] {$z_1$};
\begin{scriptsize}
\draw [color=black] (-3.,3.)-- ++(-2.0pt,-2.0pt) -- ++(4.0pt,4.0pt) ++(-4.0pt,0) -- ++(4.0pt,-4.0pt);
\draw [fill=black] (-3.,4.) circle (2.0pt);
\draw [fill=black] (-3.,5.) circle (2.0pt);
\end{scriptsize}
\end{tikzpicture}}
&$\cfrac{1}{2}$& $-\frac{1}{4}P(z_1)P(z_1+z_2)\cdot \newline\cdot P(2z_2)$&  $A+B$\\
\hline
 \raisebox{-0.5\totalheight}{\begin{tikzpicture}[line cap=round,line join=round,>=triangle 45,x=1.0cm,y=1.0cm]
\clip(-1,2.18) rectangle (1.5,6.94);
\draw [shift={(1.298235294117647,4.)},line width=1.pt]  plot[domain=2.4852404429949058:3.7979448641846805,variable=\t]({1.*1.638723551698923*cos(\t r)+0.*1.638723551698923*sin(\t r)},{0.*1.638723551698923*cos(\t r)+1.*1.638723551698923*sin(\t r)});
\draw [shift={(0.6977570093457943,4.690841121495327)},line width=1.pt]  plot[domain=2.724512232197865:4.339887162377095,variable=\t]({1.*0.7631802252741182*cos(\t r)+0.*0.7631802252741182*sin(\t r)},{0.*0.7631802252741182*cos(\t r)+1.*0.7631802252741182*sin(\t r)});
\draw [shift={(-0.3971428571428568,4.24)},line width=1.pt]  plot[domain=-0.3080527810237754:1.0892668684191507,variable=\t]({1.*0.8575094454171293*cos(\t r)+0.*0.8575094454171293*sin(\t r)},{0.*0.8575094454171293*cos(\t r)+1.*0.8575094454171293*sin(\t r)});
\draw (-0.7,6.8) node[anchor=north west] {$E$};
\draw (0.2,3.3) node[anchor=north west] {$-\frac{1}{4}$};
\draw (0.5,4.2) node[anchor=north west] {$1$};
\draw (0.2,5.3) node[anchor=north west] {$1$};
\draw (-0.6,5.3) node[anchor=north west] {$z_2$};
\draw (-0.2,4.2) node[anchor=north west] {$z_1$};
\draw (-0.1,4.7) node[anchor=north west] {$-$};
\draw (0.3,4.7) node[anchor=north west] {$+$};
\begin{scriptsize}
\draw [color=black] (0.,3.)-- ++(-2.0pt,-2.0pt) -- ++(4.0pt,4.0pt) ++(-4.0pt,0) -- ++(4.0pt,-4.0pt);
\draw [fill=black] (0.42,3.98) circle (2.0pt);
\draw [fill=black] (0.,5.) circle (2.0pt);
\end{scriptsize}
\end{tikzpicture}}
&$1$& $-\frac{1}{4}P(z_2)P(z_1-z_2)\cdot \newline\cdot P(z_1+z_2)$&  $A+B$\\
\hline
 \raisebox{-0.5\totalheight}{\begin{tikzpicture}[line cap=round,line join=round,>=triangle 45,x=1.0cm,y=1.0cm]
\clip(2,2.18) rectangle (4.5,6.94);
\draw [line width=1.pt] (3.,3.)-- (3.,5.);
\draw [line width=1.pt] (3.,5.5) circle (0.5cm);
\draw [line width=1.pt] (3.9,4.) circle (0.5015974481593781cm);
\draw (2.3,6.8) node[anchor=north west] {$F$};
\draw (3.2,3.3) node[anchor=north west] {$-\frac{1}{4}$};
\draw (3.5,4.2) node[anchor=north west] {$1$};
\draw (3.3,5.1) node[anchor=north west] {$1$};
\draw (2.4,5.1) node[anchor=north west] {$z_2$};
\draw (2.9,4.2) node[anchor=north west] {$z_1$};
\draw (3.9,4.2) node[anchor=north west] {$+$};
\draw (3.5,5.7) node[anchor=north west] {$-$};
\begin{scriptsize}
\draw [color=black] (3.,3.)-- ++(-2.0pt,-2.0pt) -- ++(4.0pt,4.0pt) ++(-4.0pt,0) -- ++(4.0pt,-4.0pt);
\draw [fill=black] (3.4,3.96) circle (2.0pt);
\draw [fill=black] (3.,5.) circle (2.0pt);
\end{scriptsize}
\end{tikzpicture}}
&$\cfrac{1}{2}$& $-\frac{1}{4}S_0 P(z_2)P(2z_2)$&$-\frac{1}{4}\left(G_2+\frac{1}{24}\right)\left(-\frac{20}{3}G_{42}+\right. \newline\left.80G_{22}^2 -80G_2G_{22}+16G_2^2\right)$\\
\hline
\end{tabular}
\caption{Graphs for $\cQ(2,1, -1^3)$: Contribution of ${g}^{deg}_{(3,1,1,1)}\cdot  p_1\cdot f_2$}
\end{figure}
\par
We provide details of the method of Theorem~\ref{thm:coeff0QM} for the
graphs A and D, to illustrate the algorithm in the proof.
We denote by $T= \{0, \tfrac{1}{2}, \tfrac{\tau}{2}, \tfrac{1+\tau}{2} \}$
the set of 2-torsion points and let $T^*= T\setminus\{0\}$. 
Note that, considering the residues at the poles at $T$, we have
\[
 %Z(2z;\tau)&\=\frac{1}{2}\left(\sum_{a\in T}Z(z-a;\tau)+1\right)\\
%Z(2z;2\tau)&\=\frac{1}{2}\left(Z(z;\tau)+Z(z-1/2;\tau)\right)\\
P(2z;\tau) \=\frac{1}{4}\sum_{a\in T}P(z-a;\tau)\,.
\]
Consequently, $P(2z)P(z)$ is an elliptic function with a pole of order~$4$
at~$0$, and poles of order~$2$ at~$T^*$, whose residues are easily compensated by $\frac{1}{4}\frac{1}{6}P''(z)$, $\frac{1}{4}P(z-a)P(a)$ for $a\in T^*$ respectively. Then, compensating the remaining pole of order 2 at 0, we get
\bas P(2z)P(z)&=\frac{1}{24}P''(z)+G_2P(z)+\frac{1}{4}\sum_{a\in T^*}P(a)\left(P(z)+P(z-a)\right) +\frac{5}{3}G_4-4G_2^2.\eas 
Since $P$ is the derivative of a $1-$periodic function $Z$, the contour integral is then reduced to 
\[[\zeta^0]P(2z)P(z)\=\frac{5}{3}G_4-4G_2^2\,.\]
Proceeding similarly with the term $P(2z)P''(z)$ we get 
\bas P(2z)P''(z) & = \frac{1}{80}P^{(4)}(z)+\frac{G_2}{2}P''(z)+\frac{49}{2}G_4P(z)\\
& \+ \frac{1}{4}\sum_{a\in T^*}\left(P(a)P''(z)+P''(a)P(z-a)\right) +\frac{28}{5}G_6-64G_4G_2\,,
\eas
hence
\[[\zeta^0]P(2z)P''(z) \=\frac{28}{5}G_6-64 G_4G_2\,,\]
which gives the contribution of the graph $A$.  The contributions of the
graphs $B$ and $C$ are computed similarly using the decomposition 
\[P(2z;2\tau)\=\frac{1}{4}\left(P(z;\tau)+P(z-1/2;\tau)\right).\]
\par
Computing the contribution of the graph $D_a$ we will see quasi-elliptic
functions appearing. We first decompose $P(z_1-z_2)P(2z_2)$ in the additive
basis with respect to~$z_2$. For this purpose decompose as usual $P(2z_2)$
into the sum of the four contributions of the 2-torsion points.
The term $P(z_1-z_2)P(z_2)$ has a pole of order 2 at $z_2=z_1$ (which is
compensated by $P(z_1-z_2)P(z_1)$),  a pole of order 2 at $z_2=0$ (which is
compensated by $P(z_1)P(z_2)$), and it also has a pole of order one at
$z_2=z_1$ (which is compensated by $Z(z_1-z_2)P'(z_1)$), and a pole of
order 1 at $z_2=0$ (finally compensated by $Z(z_2)P'(z_1)$) Proceeding
similarly for the three other two-torsion points, we get
\bas P(z_1-z_2)P(2z_2)&= \frac{1}{4}\sum_{a\in T}\Big( [P(z_2-a)+P(z_2-z_1)]P(z_1-a)\\&+[Z(z_2-a)-Z(z_2-z_1)]P'(z_1-a)\Big) +R(z_1)
\eas 
where %\[R(z_1)=\frac{4}{3}P''(2z1)-\frac{1}{4}\sum_a P'(z1-a)Z(z1-a)+\frac{1}{6} P(2z1)E2+\frac{1}{144}(E2^2-E4)\]
\[R(z_1)\=\frac{4}{3}P''(2z_1)-\frac{1}{4}\sum_{a\in T}\Big( P'(z_1-a)Z(z_1-a)\Big)-4 P(2z_1)G_2+\frac{1}{6}G_2^2-\frac{5}{3}G_4.\]
As a result,
\[[\zeta_2^0]P(z_1)P(z_1-z_2)P(2z_2) \= P(z_1)R(z_1)\,.\]
We already showed how to treat terms like $P(z_1)P''(2z_1)$, $P(z_1)P(2z_1)$,
so we focus on  $P(z_1)P'(z_1-a)Z(z_1-a)$ in the sequel. The product
$S_2 = P(z_1)P'(z_1-a)$ is an elliptic function, so examining the pole orders
we get 
\[S_2 \= \begin{cases} 
%\frac{1}{12}P^{(3)}(z_1)-\frac{1}{12}P'(z_1)E_2 \mbox{ if }a=0\\
\frac{1}{12}P^{(3)}(z_1)+2P'(z_1)G_2 & \mbox{ if }a=0\\
P''(a)[Z(z_1-a)-Z(z_1)]+P(a)P'(z_1-a) &\mbox{ if }a=1/2\\
P''(a)[Z(z_1-a)-Z(z_1)+\frac{1}{2}]+P(a)P'(z_1-a)
&\mbox{ if }a \in \{\tau/2, (1+\tau)/2 \}\,.\end{cases}
\]
When multiplying the right hand side by $Z(z_1-a)$ the right hand
side belongs to the additive basis, except for the case $Z(z_1)Z(z_1-a)$
for $a\in T^*$. But since 
\[\Delta \Big( Z(z_1)Z(z_1-a)\Big) \=
\frac{1}{2}\Delta\Big(Z(z_1)^2+Z(z_1-a)^2\Big)\,,\]
the function $Z(z_1)Z(z_1-a)-\frac{1}{2}(Z(z_1)^2+Z(z_1-a)^2)$ is elliptic
with poles of order less or equal to 2 at $0$ and $a$, so we get
\bas Z(z_1)Z(z_1-a)&=\frac{1}{2}\Big(Z(z_1)^2+Z(z_1-a)^2 -P(z_1)-P(z_1-a)\Big)+3G_2-\frac{1}{2}P(a)\\&+ \begin{cases}0 &\mbox{  if }a=1/2\\
\frac{1}{2}[Z(z-a)-Z(a)]+\frac{1}{8} &\mbox{  if }a=\tau/2, (1+\tau)/2\end{cases}\eas
and finally
\[ [\zeta_1^0]Z(z_1)Z(z_1-a) \=
G_2-\frac{1}{2}P(a)+\begin{cases}\frac{1}{6} &\mbox{if} a=1/2
\\ -\frac{5}{24} &\mbox{  if }a=\tau/2, (1+\tau)/2\,. \end{cases}\]
In total the contribution of graph $D_a$ is a combination of the
quasimodular forms
\bas A &=(-40G_{22}+10G_2)G_{42}+352G_{22}^3-408G_2G_{22}^2+144G_2^2G_{22}-16G_2^3 \\
B & = -5/4G_{42}+12G_{22}^2-12G_2G_{22}+3G_2^2\eas
as indicated in the table. Note that each contribution is a quasimodular
form for $\Gamma(2)$, but also a series in $q$, so in fact it is a quasimodular
form for $\Gamma_0(2)$ (as we remarked already Lemma~\ref{le:EisG2}). The
sum of all contributions in the table is finally the quasimodular form
given at the beginning of this subsection.
%
%The generated series for pillowcases covers with this ramification profile is then:
%\bas N^\circ(\Pi)&\=(72G_{2}(q^2) + 3)G_2(q)^2 + (-360G_{2}(q^2)^2 - 15G_{2}(q^2))G_2(q) + 360G_{2}(q^2)^3\\& + 15G_{2}(q^2)^2 - 30G_{4}(q^2)G_{2}(q^2) - 5/4G_{4}(q^2),\eas
%as it can be checked directly by computing the first terms of the series. Note that this is a quasimodular form of mixed weight $4$ and $6$.

%%%%%%%%%%%%%%%%%%%%%%%%%%%%%%%%
\subsection{Siegel-Veech weight}
%%%%%%%%%%%%%%%%%%%%%%%%%%%%%%%
Everything is ready to compute the contributions of these graphs with Siegel-Veech weight. We have the same graphs with the same local polynomials, we just associate a slightly modified propagator to take care of the weight. The recipe to get this propagator from the old one is simple, for example if they are no distinguished loops: for each edge of the graph replace the corresponding factor $P$ by $L$ and $P^{m}$ by $D_qP^{m-2}$ if $m\geq 2$, and then sum over all edges of the graph (see Section~\ref{sec:QMSV} for precise statement). The last step of integration is then similar to the previous case. We give the results in the Table~\ref{fig:examplesv} (we do not copy the factors $1/|\Aut(\Gamma)|$ which are the same; we group the graphs of same type). In the column contribution, $lwt$ stands for lower weight terms.
\begin{figure}
\label{table:examplesv}
\begin{tabular}{|c|p{5.1cm}|p{5cm}|}
\hline
Graph & Propagator & Contribution\\
\hline
$A$ & $\frac{1}{2}\Big[\Big(\frac{1}{24}D_qP(z)+\frac{1}{3}L(z)\Big)P(2z) \newline +\Big(\frac{1}{24}P^{(2)}(z)+\frac{1}{3}P(z)\Big)L(2z)\Big] $
& $\frac{1}{2}\Big[-\frac{ 880}{3}G_{22}^3 + 340G_{2}G_{22}^2 -120G_{2}^2G_{22}  + \frac{ 100}{3}G_{42}G_{22} + \frac{ 40}{3}G_{2}^3  + -\frac{ 25}{3}G_{42}G_{2}\Big] + lwt$\\
\hline
$B$ & $\frac{1}{2}\cdot\frac{1}{2}\Big(2L_{odd}(z)P_{odd}(z)P_{even}(z)\newline +P_{odd}(z)^2L_{even}(z)\Big)\newline + \frac{1}{6}\cdot\frac{3}{2}\Big(3P_{even}(z)^2L_{even}(z)\Big)$
& $\frac{1}{2}\Big[-\frac{ 560}{3}G_{22}^3 + 280G_{2}G_{22}^2 -160G_{2}^2G_{22} + \frac{ 100}{3}G_{42}G_{22} + \frac{ 80}{3}G_{2}^3 - \frac{ 50}{3}G_{42}G_{2}\Big]$\\
\hline
$C$ & $\frac{1}{4}\Big[\Big(\frac{1}{2} S_{0,odd}+\frac{3}{2}S_{0,even}\Big)\cdot \newline \Big(L_{even}(z)P(2z) + P_{even}(z)L(2z)\Big)\newline + \Big(\frac{1}{2}S_{0,odd}^{SV}+\frac{3}{2}S_{0,even}^{SV}\Big)\cdot\newline P_{even}(z)P(2z)\Big] $
&  $\frac{1}{4}\Big[180G_{22}^3 -115G_{2}G_{22}^2 -29G_{2}^2  G_{22} -15G_{42}G_{22} + 13G_{2}^3  -\frac{ 65}{12}G_{42} G_{2} \Big]+ lwt$\\
\hline
$D$ & $\frac{1}{2}\cdot(-\frac{1}{4})\cdot\Big( L(z_1)P(z_1-z_2)P(2z_2) \newline + P(z_1L(z_1-z_2)P(2z_2)\newline + P(z_1)P(z_1-z_2)L(2z_2) \newline 
+L(z_1)P(z_1+z_2)P(2z_2) \newline + P(z_1L(z_1+z_2)P(2z_2)\newline + P(z_1)P(z_1+z_2)L(2z_2)\Big) $
&$\frac{1}{2}\Big[ 352G_{22}^3 -408G_{2}G_{22}^2 + 144G_{2}^2G_{22} - 40G_{42}G_{22}  -16G_{2}^3  + 10G_{42}G_{2} \Big] + lwt$
\\
\hline
$E$ & $-\frac{1}{4}\Big( L(z_1)P(z_1-z_2)P(z_1+z_2) \newline + P(z_1L(z_1-z_2)P(z_1+z_2)\newline + P(z_1)P(z_1-z_2)L(z_1+z_2)\Big) $
&$264G_{22}^3 -306G_{2}G_{22}^2 + 108G_{2}^2G_{22}  - 30G_{42}G_{22} -12G_{2}^3  + \frac{ 15}{2}G_{42}G_{2} + lwt$\\
\hline
$F$ & $\frac{1}{2}\cdot(-\frac{1}{4})\cdot\Big(P(z_2)P(2z_2) + \newline L(z_2)P(2z_2)+P(z_2)L(2z_2)\Big)$ &$\frac{1}{2}\Big[-\frac{ 65}{2}G_{2}G_{22}^2 + \frac{ 65}{2}G_{2}^2G_{22} -\frac{ 13}{2}G_{2}^3  + \frac{ 65}{24}G_{42}G_{2} \Big] + lwt$\\
\hline
\end{tabular}
\caption{Graphs for $\cQ(2,1, -1^3)$: Siegel-Veech contribution}
\label{fig:examplesv}
\end{figure}
In the compilation of the table we used 
\bas S_{0}^{SV}&\=\sum_{w,h=1}^{\infty} hq^{wh}\=S_0\=G_2+\frac{1}{24}\\
S_{0, even}^{SV}&\=\sum_{w,h=1}^{\infty} h q^{(2w)h}\= G_{22}+\frac{1}{24}
\eas
and $S_{0_odd}^{SV} \= S_0^{SV}-S_{0, even}^{SV}$. Summing up, we get the
generating series
\bas c_{-1}^0(\Pi)&\= 245G_{22}^3 -245G_2G_{22}^2 + 49G_2^2G_{22}  -\frac{245}{12}G_{42}G_{22} + lwt.
\eas
The lower weight terms are weight $2$ and $4$ terms, as we can expect from the weight of $N^{0}(\Pi)$ (as $L$ contains lower weight terms). This series is proportional to $N^{0}(\Pi)$, since this stratum is non-varying (see \cite{chenmoeller} for more explanation). We will see in the next section that the coefficient of proportionality is the Siegel-Veech constant of the stratum.

\subsection{Contributions to volumes of strata and Siegel-Veech constants}

Evaluation of volumes and Siegel-Veech constants of strata are closely related to the asymptotics of the generating series $N^0(\Pi)$ and $c_{-1}(\Pi)$ as $q$ tends to 1 (or equivalently $\tau$ tends to 0), as stated in \cite{eopillow} (see also \cite[Proposition 7]{GoujardExplVal}). This asymptotics can be easily obtained thanks to the quasimodularity property for Eisenstein series: the transformation $\tau\to -1/\tau$ relates the asymptotics as $\tau\rightarrow 0$ to the asymptotics as $\tau\rightarrow i\infty$.

We overview briefly the results of \cite[Section 9]{CMZ} here and define  two polynomials describing the growth of a quasimodular form (for $\Gamma_0(2)$) near $\tau=0$, so at the same time the average growth of its Fourier coefficients. 

We recall that $QM(\Gamma_0(2))$ is the space of even weight quasimodular forms for $\Gamma_0(2)$, and it is generated as a polynomial ring by $G_2, G_{22}, G_{42}$, see \eqref{eq:G02gens}.

\begin{Defi}We define the map $\Ev$ as the unique algebra homomorphism from $QM(\Gamma_0(2))$ to $\Q[X]$ sending $G_2$ to $-X/24-1/2$, $G_{22}$ to $-X/96-1/4$ and $G_{42}$ to $X^2/3840$.
\end{Defi}

Setting $h=-2\pi i\tau$, we define 
\bas \ev[F](h)\=\frac{1}{h^k}\Ev[F]\Big(-\frac{4\pi^2}{h}\Big)\; \in \Q[\pi^2][1/h]
\eas
for $F\in QM_{2k}(\Gamma_0(2))$ (weight $2k$ quasimodular form). This polynomial describes the growth of $F(\tau)$ near $\tau=0$ directly (also for mixed weight forms), as proved in the following Proposition. 

\begin{Prop} For $F\in QM(\Gamma_0(2))$ we have
\[F(i\v)=\ev[F](2\pi\v) + (\textrm{small})\quad (\v\searrow 0) \] where ``small'' means terms that tends exponentially quickly to 0.
\end{Prop}

\begin{proof}This is directly derived from the modularity properties $G_2(-1/\tau)=\tau^2G_2(\tau)-\tau/4\pi i$ and $G_4(-1/\tau)=\tau^4G_4(\tau)$.
\end{proof}

In Figure~\ref{table:contribexample} we give the $h$-evaluation of all individual graphs. Note that the $h$-evaluation of lower weight terms is $O\big(\frac{1}{h^4}\big)$.
\begin{figure}
%  \begin{table}    
\renewcommand{\arraystretch}{1.5}

\begin{tabular}{|c|c|c|}
\hline
Graph $\Gamma$ & $\ev[N^0(\Gamma, \Pi)](h)$ & $\ev[c_{-1}(\Gamma,\Pi)](h)$\\
\hline
$A$ & $\frac{4}{45}\frac{\pi^4}{h^5}+O\big(\frac{1}{h^4}\big)$ & $O\big(\frac{1}{h^4}\big)$\\
\hline
$B$ & $\frac{2}{15}\frac{\pi^4}{h^5}+O\big(\frac{1}{h^4}\big)$ & $\frac{5}{36}\frac{\pi^4}{h^5}+O\big(\frac{1}{h^4}\big)$\\
\hline
$C$ & $\frac{1}{9}\frac{\pi^4}{h^5}+O\big(\frac{1}{h^4}\big)$ & $\frac{11}{144}\frac{\pi^4}{h^5}+O\big(\frac{1}{h^4}\big)$\\
\hline
$D$ & $O\big(\frac{1}{h^4}\big)$ & $O\big(\frac{1}{h^4}\big)$\\
\hline
$E$ & $O\big(\frac{1}{h^4}\big)$ & $O\big(\frac{1}{h^4}\big)$\\
\hline
$F$ & $-\frac{1}{36}\frac{\pi^4}{h^5}+O\big(\frac{1}{h^4}\big)$ & $-\frac{13}{288}\frac{\pi^4}{h^5}+O\big(\frac{1}{h^4}\big)$\\
\hline
Total & $\vol= \pi^4/3072$ &  $\frac{\pi^2}{3}c_{area}=\frac{49}{72}$\\
\hline
\end{tabular}
\caption{Graphs for $\cQ(2,1, -1^3)$: growth polynomials and contribution to the volume and the Siegel-Veech constant of the stratum}\label{table:contribexample}
%\end{table}
\end{figure}
In this table \[\vol=\frac{2\dim}{2^{\dim} \dim !}\cdot \lim\limits_{h\to 0} \big(\ev[N^{0}(\Pi)](h)\cdot h^{\dim}\big)\] is the volume of $\cQ(2,1, -1^3)$ (with respect to the Eskin-Okounkov convention ; the volume for the Athreya-Eskin-Zorich convention has an additional factor $4^5/2\cdot 3!=3072$, see \cite[Lemma 2]{GoujardExplVal} for discussion concerning volume normalizations), and $\dim=\dim_\C\cQ(2, 1, -1^3)=5$. 

In the general case, the Siegel-Veech constant $c_{area}$ is computed using
\[\frac{\pi^2}{3}c_{area}\=\lim\limits_{D\to\infty}\cfrac{\sum_{d=1}^{D}c_{-1}^0(d,\Pi)}{\sum_{d=1}^{D}N^{0}(d,\Pi)}\=\cfrac{\lim\limits_{h\to 0} \big(\ev[c_{-1}^{0}(\Pi)](h)\cdot h^{\dim}\big)}{\lim\limits_{h\to 0} \big(\ev[N^0(\Pi)](h)\cdot h^{\dim}\big)},\] since both limits are finite.

In our example, the stratum  is non varying so the series $c_{-1}^0(\Pi)$ and $N^{0}(\Pi)$ are proportional (not the individual graph contributions though). The Siegel-Veech constant is just the proportionality factor. We get the following constant,
\[\frac{\pi^2}{3}c_{area}(\cQ(2,1, -1^3))=\frac{49}{72}.\]
in agreement with the value computed in \cite{chenmoeller}.
%Interpretation
%\bas
%\ol{g}_{3111}&\=\frac{1}{2}f_1g_{11} + g_{3111}\\
%g_{3111}^{deg}&\=-\frac{1}{2}g_{11}\\
%p_1&\= f_1-1/24
%\eas

\bibliography{my}
\end{document}